\documentclass[10pt]{article}
\usepackage{makeidx}
\usepackage{multirow}
\usepackage{multicol}
\usepackage{graphicx}
\usepackage{epstopdf}
\usepackage{ulem}
\usepackage{hyperref}
\usepackage{amsmath}
\usepackage{amssymb}
\usepackage{setspace}
\usepackage{color}
\usepackage{float}
\author{Aishwarjya Gogoi}
\title{}
\usepackage[paperwidth=612pt,paperheight=792pt,top=72pt,right=72pt,bottom=72pt,left=72pt]{geometry}

\makeatletter
	{\par\setlength{\parindent}{#3}
	\setlength{\leftmargin}{#1}       \setlength{\rightmargin}{#2}%
	\advance\linewidth -\leftmargin       \advance\linewidth -\rightmargin%
	\advance\@totalleftmargin\leftmargin  \@setpar{{\@@par}}%
	\parshape 1\@totalleftmargin \linewidth\ignorespaces}{\par}%
\makeatother 

\begin{document}
\begin{doublespace}
\begin{center}

\textbf{{\LARGE Enhanced HLLEM and HLL-CPS schemes for all Mach number flows based using anti-diffusion coefficients}}

\end{center}
\begin{center}
A. Gogoi*, J. C. Mandal* and A. Saraf**\\
*Department of Aerospace Engineering, Indian Institute of Technology, Bombay, Mumbai, 400076, India\\
**Aeronautical Development Agency, Vimanapura, Bengaluru 560017, India
\end{center}
\end{doublespace}
\begin{center}
\uppercase{ABSTRACT}
\end{center}
\begin{doublespace}
This paper compares the HLLEM and HLL-CPS schemes for Euler equations and proposes improvements for all Mach number flows. Enhancements to the HLLEM scheme involve adding anti-diffusion terms in the face normal direction and modifying anti-diffusion coefficients for linearly degenerate waves near shocks. The HLL-CPS scheme is improved by adjusting anti-diffusion coefficients for the face normal direction and linearly degenerate waves. Matrix stability, linear perturbation, and asymptotic analyses demonstrate the robustness of the proposed schemes and their ability to capture low Mach flow features. Numerical tests confirm that the schemes are free from shock instabilities at high speeds and accurately resolve low Mach number flow features.

\section{Introduction }

Godunov-type approximate Riemann solvers \cite{godu, toro-book} are widely used for computing convective fluxes in the Euler and Navier-Stokes equations. Notable solvers include the Roe scheme \cite{roe}, HLL \cite{hll}, HLLE \cite{einf1}, Osher and Solomon scheme \cite{osher}, HLLEM \cite{einf2}, and HLLC \cite{toro1, batt}. The Roe scheme \cite{roe} linearizes the Euler equations, accurately resolving isolated shocks, contact, and shear discontinuities. The HLL scheme \cite{hll} also resolves isolated shocks but neglects intermediate contact and shear waves, limiting its practical use. The HLLE scheme \cite{einf1}, an extension of HLL with Einfeldt's wave-speed estimates, satisfies entropy conditions and ensures positivity. The HLLEM scheme \cite{einf2} enhances HLLE by adding anti-diffusion terms to capture intermediate waves, offering accuracy comparable to Roe's scheme while retaining HLLE's entropy and positivity features. The HLLC scheme \cite{toro1} further restores contact and shear waves in the HLLE framework by enforcing Rankine-Hugoniot jump conditions, maintaining both entropy and positivity. Mandal and Panwar's HLL-CPS scheme \cite{mandal} splits flux into convective and pressure components, with the convective part treated by an upwind method and the pressure flux evaluated using HLL principles, replacing density jumps with pressure jumps based on isentropic flow. While this method resolves contact discontinuities exactly, it struggles with boundary layer accuracy \cite{kita1, kita-slau2}.

It is observed that the approximate Riemann solvers \cite{roe, osher, einf2, toro1} that resolve the linearly degenerate intermediate waves present in the Riemann problem are susceptible to numerical shock instabilities like carbuncle phenomenon, kinked Mach stem, odd-even decoupling and low-frequency post-shock oscillation for slowly moving shock. On the other hand, the HLLE scheme \cite{hll, einf1} which does not resolve the linearly degenerate intermediate waves, is free from numerical shock instabilities. The approximate Riemann solvers are incapable of resolving flow features at low Mach numbers, especially on structured grids. Some of the Riemann solvers that have been proposed to resolve the low Mach flow features have been found to suffer from checkerboard-type decoupling of the pressure and velocity fields. The cause of numerical instabilities in the complete approximate Riemann solvers has been widely investigated ever since the numerical instabilities were first presented by Peery and Imlay \cite{peery}. Liou \cite{liou1} has conjectured that the schemes that suffer from numerical shock instability have a pressure jump term in the numerical mass flux while those free from shock instability do not contain a pressure jump term in the mass flux. Based on the conjecture of Liou, Kim et al. \cite{kim} have proposed a shock-stable Roe scheme by reducing the contribution of the pressure jump terms in the mass flux. Quirk \cite{quirk} has carried out linearized perturbation analysis and observed that the schemes in which the pressure perturbations feed the density perturbations will be afflicted by the odd-even decoupling problem. Quirk has proposed a hybrid Roe-HLLE scheme \cite{quirk} to overcome the numerical instability problems of the Roe scheme. Pandolfi and D'Ambrosio \cite{pand} have carried out a more refined perturbation analysis by including the velocity perturbations and have observed that the schemes in which the density perturbations are fed by pressure perturbations and the shear velocity perturbation are undamped exhibit severe carbuncle phenomenon while the schemes in which the density perturbations are damped but the shear velocity perturbations are undamped are prone to light carbuncle. Dumbser et al. \cite{dumb-matrix} have developed a matrix stability method to evaluate the stability of upwind schemes and have observed that the origin of the carbuncle phenomenon is localized in the upstream region of the shock. Several rotated Roe solvers \cite{ren, nishi}, where the flux evaluation direction gets rotated when the velocity difference exceeds a threshold value, have been proposed to overcome the numerical shock instability problems of the Roe scheme. However, the rotated Riemann solvers require additional computational time as the Riemann problem needs to be solved twice. Park and Kwon \cite{park} have proposed a modification in the anti-diffusion coefficients of the linearly degenerate intermediate waves to overcome the numerical shock instability problems in the HLLEM scheme. Shen et al. \cite{shen} have identified the exact resolution of the shear wave as the cause of the numerical instabilities in the HLLC schemes and have proposed a shock-stable, hybrid HLLC-HLLCM scheme which does not resolve the shear wave in the vicinity of the shock. Similarly, Xie et al. \cite{xie-hllem} have identified the shear wave as the cause of instability in the HLLEM  Riemann solver and introduced a pressure function in the HLLEM scheme such that the anti-diffusion terms corresponding to the shear wave is drastically reduced in the vicinity of the shock. Chen et al. \cite{chen-roe} have proposed a carbuncle-free Roe scheme by introducing shear viscosity into the momentum flux in the vicinity of the shock through a pressure-based sensing function while preserving the shear wave elsewhere. Chen et al. \cite{chen-all-scheme} have also proposed the addition of shear viscosity terms into the Godunov-type schemes to cure the schemes from the carbuncle phenomenon along with a multi-dimensional pressure-based function to restore the shear wave. However, the diagnosis of the shear wave as the cause of the numerical shock instabilities in the HLL-type scheme by Shen et al. \cite{shen}, Xie et al. \cite{xie-hllem}, and Chen et al. \cite{chen-roe, chen-all-scheme} is yet to be universally accepted. Simon and Mandal \cite{san3} have shown that the HLLEMS scheme of Xie et al. \cite{xie-hllem} is not fully free from numerical shock instability problems. Simon and Mandal \cite{san3} have shown that suppression of the anti-diffusion terms of both the contact and shear waves is more effective in curing the numerical shock instability problems in the HLLEM Riemann solver. Simon and Mandal \cite{san1, san2, san3} have also proposed strategies like anti-diffusion control (ADC) and Selective Wave-speed Modification (SWM) to cure the numerical instability problems in the HLLC and HLLEM Riemann solver. Kemm \cite{kemm} has shown that the resolution of the contact wave in the HLLEM Riemann solver also contributed to the numerical shock instabilities, although its contribution to the numerical shock instabilities was less as compared to the shear wave. Fleischmann et al. \cite{fleis1, fleis2} have identified the inappropriate numerical dissipation in the transverse direction of a grid-aligned shock to be the cause of the numerical shock instability in the Godunov-type approximate Riemann solvers. The inappropriate numerical dissipation in the transverse direction is attributed to the low Mach number in the transverse direction of the grid-aligned shock. Fleischmann et al. \cite{fleis1} have proposed a reduction of the overall numerical dissipation in the classical Roe scheme and the component-wise local Lax-Fredrichs (cLLF) scheme to suppress the shock instability. The decrease in the overall numerical dissipation is achieved by modifying the non-linear acoustic eigenvalues. Fleischmann et al. \cite{fleis2} have observed that  similar modification to the non-linear eigenvalues in the HLLC scheme leads to a dissipative upwind scheme. Fleischmann et al. \cite{fleis2} have proposed a centred HLLC scheme comprising central-difference terms and numerical dissipation, and proposed a modification of the non-linear eigenvalues in the numerical dissipation terms to suppress shock instability. Fleischmann et al. \cite{fleis2} have demonstrated reduced overall numerical dissipation with the modification of the non-linear eigenvalues and have shown accurate results over a cylinder at low Mach numbers.

Approximate Riemann solvers often produce unphysical results in low Mach number flows. To address this, various methods have been developed, including preconditioning \cite{guillard1, park-pre}, non-linear eigenvalue modification \cite{li-roe}, velocity jump scaling \cite{rieper1, rieper2}, velocity reconstruction \cite{th1, th3}, central differencing of momentum flux \cite{della1}, and the addition of anti-diffusion terms \cite{della2}, among others, to enhance their performance in such flows. It was shown by Guillard and Viozat \cite{guillard1} and Rieper \cite{rieper1, rieper2} that, in the limit as Mach number tends to zero, the pressure fluctuations in the Roe scheme for the discrete Euler equations are of the order of Mach number, while for the continuous Euler equations, the pressure fluctuations are of the order of Mach number squared. Guillard and Viozat \cite{guillard1} have proposed a preconditioned Roe scheme that had a correct pressure scaling with Mach number. Park et al. \cite{park-pre} have proposed a preconditioned HLLE+ scheme for low-speed flows. Li and Gu \cite {li-roe} have observed that the preconditioned Roe schemes suffer from the global cut-off problem. They have proposed an all-speed Roe scheme by modifying the non-linear eigenvalues of the numerical dissipation terms \cite{li-roe} in the Roe scheme. The modification of the non-linear eigenvalues comprises scaling of the acoustic speeds by a local Mach number function. They \cite{li-roe} have shown that the zero-order velocity field in the all-speed Roe scheme is subject to a divergence constraint and the second-order pressure satisfies a Poisson-type equation in the case of constant entropy. They have also observed that the zero-order velocity field computed by the traditional preconditioned Roe scheme does not satisfy the divergence constraint as the Mach number vanishes. Rieper \cite{rieper2} has carried out an asymptotic analysis of the Roe scheme and identified the jump in the normal velocity component as the cause of the loss of accuracy at low Mach numbers. He \cite{rieper2} has proposed a low Mach number fix for the Roe scheme called LMRoe comprising scaling of the normal velocity jump in the numerical dissipation terms by a local Mach number function. He has demonstrated accurate results for inviscid flow past a cylinder at very low Mach numbers with the modified LMRoe scheme, without any checkerboard-type velocity-pressure decoupling. Thornber et al. \cite{th1} have proposed a modification of the reconstructed velocities at the cell interface of the Godunov-type schemes like HLLC \cite{toro1} such that the velocity difference reduces with decreasing Mach number and the arithmetic mean is reached at Mach zero. They have demonstrated analytically and numerically that the numerical dissipation of the scheme became constant in the limit of zero Mach, while the numerical dissipation of the traditional schemes tended to infinity in the limit of zero Mach. The velocity reconstruction method of Thornber et al.\cite{th1} has been used in the HLLC scheme to obtain accurate results for a wide range of Mach numbers \cite{th3}. The low Mach number fix of Rieper has been applied to both the normal and tangential velocity jumps in the numerical dissipation terms of the Roe scheme by Osswald et al. \cite{ossw}. Accuracy similar to the method of Thornber et al. \cite{th1} has been obtained for turbulent flows by Osswald et al. Osswald et al. \cite{ossw} have also conjectured that Thornber et al.'s method may have undesirable non-linear side-effects as the method involved modification of the physical flux function. Dellacherie \cite{della1} has observed that the Godunov-type schemes cannot be accurate in low Mach numbers due to loss of invariance. He has proposed central differencing to discretize the momentum flux or to discretize the pressure gradient in momentum flux to recover invariance property and improve accuracy in low Mach number flows. Dellacherie et al. \cite{della2} have also proposed stable, all-Mach Godunov-type schemes by adding Mach number-based anti-diffusion terms to the Godunov-type schemes such that central-differencing is obtained in the limit of Mach zero. The theoretical framework of Dellacherie has justified the existence of all-Mach number Godunov-type schemes. Yu et al. \cite{yu-hllem2} have implemented Dellacherie et al.'s low Mach correction in the HLLEM scheme. They have also proposed a scaling of the terms responsible for the inaccurate behaviour of the HLLEM scheme at low Mach numbers with a Mach number function. 

It is observed that most of the analyses and fixes for the high Mach number numerical shock instabilities and the low Mach number inaccuracies have been proposed for the Roe \cite{xie-roe}, AUSM \cite{shima, chen-ausm}, HLLC \cite{xie-hllc, chen-hllc}, and HLLEM \cite{qu-hllem} schemes. The HLL-CPS scheme \cite{mandal} has been numerically shown to be positivity preserving, entropy satisfying, and contact capturing. However, the HLL-CPS scheme exhibit numerical shock instability under certain conditions and is incapable is resolving shear layers and low Mach flow features. In this paper, an attempt is made to eliminate the numerical shock instabilities in the HLL-CPS scheme and also improve its ability to resolve shear layers and low Mach flow features. A set of anti-diffusion coefficients are proposed for the HLL-CPS scheme which helps suppress the numerical shock instabilities for high-speed flows and also enables the resolution of shear layers and the low Mach flow features. A similar set of anti-diffusion coefficients is also proposed for the HLLEM and HLLC schemes for suppressing numerical shock instability and resolving low Mach flow features. 

Most of the analyses and solutions for numerical shock instabilities at high Mach numbers and inaccuracies at low Mach numbers have been developed for the Roe \cite{xie-roe}, AUSM \cite{shima, chen-ausm}, HLLC \cite{xie-hllc, chen-hllc}, and HLLEM \cite{qu-hllem} schemes. The HLL-CPS scheme \cite{mandal} has been shown to preserve positivity, satisfy entropy conditions, and accurately capture contact discontinuities. However, the HLL-CPS scheme exhibits numerical shock instabilities under certain conditions and struggles to resolve shear layers and low Mach number flow features. In this paper, we aim to eliminate these numerical shock instabilities in the HLL-CPS scheme and enhance its ability to resolve shear layers and low Mach number flow features. We propose a set of anti-diffusion coefficients for the HLL-CPS scheme that suppresses numerical shock instabilities in high-speed flows and improves the resolution of shear layers and low Mach features. A similar set of anti-diffusion coefficients is also proposed for the HLLEM scheme to address shock instabilities and improve low Mach flow resolution.

The paper is structured into six sections. Section 2 outlines the governing equations, followed by a description of the finite volume discretization in Section 3. Section 4 covers the flux evaluation method, while Section 5 presents numerical results for various 2D test cases, spanning a wide range of Mach numbers from very low to very high. Finally, the conclusions are discussed in Section 6.

\section{Governing Equations}
The governing equations for the two-dimensional inviscid compressible flow can be expressed in their conservative form as
\begin{equation} \label{2dge}
\dfrac{\partial{}\boldsymbol{\acute{U}}}{\partial{}t}+\dfrac{\partial{}\boldsymbol{\acute{F}(\acute{U})}}{\partial{}x}+\dfrac{\partial{}\boldsymbol{\acute{G}(\acute{U})}}{\partial{}y}=0
\end{equation}
where $\boldsymbol{\acute{U}}$, $\boldsymbol{\acute{F}(\acute{U})}$, and $\boldsymbol{\acute{G}(\acute{U})}$ are the vector of conserved variables, x-directional and y-directional fluxes respectively, and can be written as
\begin{equation}
\boldsymbol{\acute{U}}=\left[\begin{array}{c} \rho \\ \rho{}u \\ \rho{}v \\ \rho{}E \end{array}\right] \hspace{1cm} \boldsymbol{\acute{F}(\acute{U})}=\left[\begin{array}{c} \rho{}u \\ \rho{}u^2+p \\ \rho{}uv \\ (\rho{}E+p)u \end{array}\right]  \hspace{1cm} \boldsymbol{\acute{G}(\acute{U})}=\left[\begin{array}{c} \rho{}v \\ \rho{}uv \\ \rho{}v^2+p \\ (\rho{}E+p)v \end{array}\right] 
\end{equation}
where $\rho$, $u$, $v$, $p$ and $E$ stand for density, x-directional and y-directional velocities in global coordinates, pressure and specific total energy. The system of equations is closed by the equation of state
\begin{equation}
p=(\gamma-1)\left(\rho{}E-\dfrac{1}{2}\rho(u^2+v^2)\right)
\end{equation} jj
where $\gamma$ is the ratio of specific heat. In this paper, a calorically perfect gas with $\gamma=1.4$ is considered. The integral form of the governing equations is
\begin{equation}\label{2dns}
\frac{d}{dt}\int_{\Omega{}}\boldsymbol{\acute{U}}d\Omega{}+\oint_{\partial{}\Omega{}}\boldsymbol{(\acute{F},\acute{G}).n}ds=0
\end{equation}
where $\partial{}\Omega$ is the boundary of the control volume  $\Omega$ (area in case of two dimensions) and  $\boldsymbol{n}$ denotes the outward pointing unit normal vector to the surface $\partial\Omega$.

\section{Finite Volume Discretization}
The finite volume discretization of equation (4) for a structured, quadrilateral mesh can be expressed as
\begin{equation}
\dfrac{d\boldsymbol{\acute{U}}_i}{dt}=-\dfrac{1}{|\Omega|_i}\sum_{l=1}^4[\boldsymbol{(\acute{F}, \acute{G})}_l.\boldsymbol{n}_l]\Delta{}s_l
\end{equation}
where $\boldsymbol{\bar{U}}_i$ is the cell averaged conserved state vector, $|\Omega|_i$ is the area of the cell, $\boldsymbol{n}_l$ denotes the unit normal vector and $\Delta{}s_l$ denotes the length of each interface. The rotational invariance property of the Euler equations is utilized to express the flux $\boldsymbol{(F,G)}_l.\boldsymbol{n}_l$ as 
\begin{equation} \label{fv}
\dfrac{d\boldsymbol{\acute{U}}_i}{dt}=-\dfrac{1}{|\Omega|_i}\sum_{l=1}^{4}(\mathcal{T}_{il})^{-1}\boldsymbol{F}(\mathcal{T}_{il}\boldsymbol{U}_i,\mathcal{T}_{il}\boldsymbol{U}_l)\Delta{}s_l
\end{equation}
where  $\boldsymbol{F}(\mathcal{T}_{il}\boldsymbol{U}_i,\mathcal{T}_{il}\boldsymbol{U}_l)$  is the inviscid  face normal flux vector, $\mathcal{T}_{il}$ is the rotation matrix and ${\mathcal{T}_{il}}^{-1}$ is its inverse at the edge between cells $i$ and $l$. 
 \begin{equation}
\mathcal{T}_{il}=\left[\begin{array}{ccccc}1& 0& 0& 0\\ 0& n_x& n_y& 0 \\  0& -n_y& n_x& 0 \\ 0& 0& 0& 1\end{array}\right] \hspace{1cm} \text{and} \hspace{1cm} {\mathcal{T}_{il}}^{-1}=\left[\begin{array}{ccccc}1& 0& 0& 0\\ 0& n_x&- n_y& 0 \\  0& n_y& n_x& 0 \\ 0& 0& 0& 1\end{array}\right]
\end{equation}
where $\boldsymbol{n}=(n_x,n_y)$ is the outward unit normal vector on the edge between cells $i$ and $l$. 

In the next section, the HLLEM and HLL-CPS schemes for the face normal flux $\boldsymbol{F}(\mathcal{T}_{il}\boldsymbol{U}_i,\mathcal{T}_{il}\boldsymbol{U}_l)$ shown in equation (\ref{fv}) shall be described and improvements to these schemes shall be proposed. 

\section{Flux Description and Improvement}
\subsection{HLLEM Flux }\label{flux-comparison}
The flux of various Godunov-type, HLL family schemes in the direction normal to the cell interface can be expressed as a sum of the HLL flux and additional anti-diffusion terms as per Einfeldt et al. \cite{einf2}, Kim et al. \cite{kim}, and Park and Kwon \cite{park}. Thus the flux of the HLL-family scheme can be written as
\begin{equation} \label{hll-type-flux}
\boldsymbol{F}(\boldsymbol{U}_R,\boldsymbol{U}_L)=\dfrac{S_R\boldsymbol{F}_L-S_L\boldsymbol{F}_R}{S_R-S_L}+\dfrac{S_RS_L}{S_R-S_L}(\Delta{}\boldsymbol{U}-B\Delta{}\boldsymbol{U})
\end{equation} 
where $S_R$ and $S_L$ are the fastest right and left running wave speeds, $\boldsymbol{U}_R=\mathcal{T}_{il}\boldsymbol{U}_i$ and $\boldsymbol{U}_L=\mathcal{T}_{il}\boldsymbol{U}_k$ are the face normal state variables to the right and left of the interface, $\boldsymbol{F}_R=\boldsymbol{F}(\boldsymbol{U}_R)$ and $\boldsymbol{F}_L=\boldsymbol{F}(\boldsymbol{U}_L)$ are the face normal fluxes to the right and left of the interface, $B\Delta{}\boldsymbol{U}$ is the anti-diffusion term, and $\Delta(.)=(.)_R-(.)_L$.
For the HLLEM scheme, the wave speeds have been defined by Einfeldt et al. \cite{einf2} as
\begin{equation} \label{wavespeed-hll}
S_L=min(0,\;u_{nL}-a_L,\;\tilde{u}_n-\tilde{a}) \hspace{1cm} S_R=max(0,\;u_{nR}+a_R,\;\tilde{u}_n+\tilde{a})
\end{equation}
Einfeldt et al. \cite{einf2} have shown that with the above choice of wave speeds, the HLLEM scheme is entropy satisfying and positively conservative\cite{einf2}. The anti-diffusion term  $B\Delta{}\boldsymbol{U}$ of the HLLEM scheme is \cite{einf2}
\begin{equation} \label{bdq-roe-hllem}
B\Delta{}\boldsymbol{U}=\delta_2\tilde{\alpha}_2\tilde{R}_2+\delta_3\tilde{\alpha}_3\tilde{R}_3
\end{equation}
where $\delta_2$ and $\delta_3$ are the anti-diffusion coefficients for the contact and shear waves, $\tilde{\alpha}_2=\Delta{}\rho-\dfrac{\Delta{}p}{\tilde{a}^2}$ and $\tilde{\alpha}_3=\tilde{\rho}\Delta{}u_t$ are the wave strength of the contact and shear waves, $R_2=\left[1,\;\tilde{u}_n,\;\tilde{u}_t,\; \frac{1}{2}(\tilde{u}_n^2+\tilde{u}_t^2)\right]^T$ and $R_3=[0,\;0,\; 1,\; \tilde{u}_t]^T$ are the right eigenvectors of the contact and shear waves in the flux Jacobian matrix. The anti-diffusion coefficients are $\delta{}_2=\delta{}_3=\frac{\tilde{a}}{\tilde{a}+|\bar{u}|}$ where $|\bar{u}|$ is the approximate speed of the contact discontinuity and is defined as $|\bar{u}|=\frac{S_R+S_L}{2}$ \cite{einf2}. Park and Kwon \cite{park} have observed that the HLLEM scheme of Einfeldt \cite{einf2} is not as accurate as the Roe scheme in presence of contact discontinuity and proposed the Roe-averaged velocity as the velocity of the contact discontinuity for improving the accuracy of the HLLEM scheme. Thus, for the HLLEM scheme with the modification of Park and Kwon,  $\delta{}_2=\delta{}_3=\frac{\tilde{a}}{\tilde{a}+|\tilde{u}_n|}$. 
\subsection{HLL-CPS Flux and its HLLEM-type Expression}
In the HLL-CPS scheme, the convective and pressure fluxes are split following either the Zha-Bilgen \cite{zha1} or AUSM splitting \cite{ausm1}. The convective flux is evaluated by an upwind method \cite{mandal} and the pressure flux is evaluated using the HLL method \cite{hll} with suitable modification for the resolution of the contact discontinuity. 
For the Zha-Bilgen splitting, the HLL-CPS flux at the cell interface is \cite{mandal}
\begin{equation} \label{hll-cps-flux}
\boldsymbol{F}(\boldsymbol{U}_R,\boldsymbol{U}_L )= M_{nK}\left[\begin{array}{c} \rho \\ \rho{}u_n \\ \rho{}u_t\\ \rho{}E\end{array}\right]_Ka_K+\dfrac{S_R\boldsymbol{F}_{2L}-S_L\boldsymbol{F}_{2R}}{S_R-S_L}+\dfrac{S_RS_L}{S_R-S_L}\left[\begin{array}{c}\dfrac{\Delta{}p}{\bar{a}^2} \\ \dfrac{\Delta{}(pu_n)}{\bar{a}^2} \\ \dfrac{\Delta{}(pu_t)}{\bar{a}^2} \\ \dfrac{\Delta{}p}{\gamma-1}+\dfrac{\Delta(pq^2)}{2\bar{a}^2} \end{array}\right]
\end{equation}
where 
\begin{equation}
K=\left\lbrace\begin{array}{c} L \:\;\; \text{if} \:\;\; \bar{u}_n\geq0.0\\ R \:\;\; \text{if} \:\;\; \bar{u}_n < 0.0\end{array}\right., \hspace{2mm}
M_{nK}=\left\lbrace\begin{array}{c}\dfrac{\bar{u}_n}{\bar{u}_n-S_L}\;\; \text{if}\;\; \bar{u}_n\geq0.0\\
\dfrac{\bar{u}_n}{\bar{u}_n-S_R}\;\; \text{if}\;\; \bar{u}_n<0.0\end{array}\right., \hspace{2mm}
a_K=\left\lbrace\begin{array}{c} \ u_{nL}-S_L \:\; \text{if} \:\; \bar{u}_n\geq0.0\\\ u_{nR}-S_R \:\; \text{if} \:\; \bar{u}_n < 0.0\end{array}\right.
\end{equation}
\begin{equation}
\bar{u}_n=\frac{1}{2}(u_{nL}+u_{nR}), \hspace{1cm} \boldsymbol{F}_{2K}=[0,\; p_K,\; 0,\; (pu_n)_K]^T
\end{equation}
The HLL-CPS pressure flux can be expressed in the form of the HLLEM scheme as a sum of the HLL flux and additional anti-diffusion terms. Thus, the HLL-CPS flux shown in equation (\ref{hll-cps-flux}) can be written as
\begin{equation} \label{hll-cps-flux2}
\boldsymbol{F}(\boldsymbol{U}_R, \boldsymbol{U}_L)= M_{nK}\left[\begin{array}{c} \rho \\ \rho{}u_n \\ \rho{}u_t\\ \rho{}E\end{array}\right]_Ka_K+\dfrac{S_R\boldsymbol{F}_{2L}-S_L\boldsymbol{F}_{2R}}{S_R-S_L}+\dfrac{S_RS_L}{S_R-S_L}\left(\Delta{}\boldsymbol{U}-B\Delta{}\boldsymbol{U}\right)
\end{equation} 
where
\begin{equation}
\Delta{}\boldsymbol{U}=\left[\begin{array}{c} \Delta{}\rho \\ \Delta{}(\rho{}u_n) \\ \Delta{}(\rho{}u_t) \\ \Delta{}(\rho{}E)\end{array}\right] \hspace{1cm}B\Delta{}\boldsymbol{U}=\left[\begin{array}{c}\Delta{}\rho{}-\dfrac{\Delta{}p}{\bar{a}^2} \\ \Delta{}(\rho{}u_n)-\dfrac{\Delta{}(pu_n)}{\bar{a}^2} \\ \Delta{}\rho{}u_t-\dfrac{\Delta{}(pu_t)}{\bar{a}^2} \\ \dfrac{\Delta{}(\rho{}q^2)}{2}-\dfrac{\Delta(pq^2)}{2\bar{a}^2} \end{array}\right]
\end{equation}
Rearranging $B\Delta{}\boldsymbol{U}$ of the HLL-CPS scheme in terms of the contact and shear waves, we obtain
\begin{equation}
B\Delta{}\boldsymbol{U}=\left(\Delta{}\rho-\dfrac{\Delta{}p}{\bar{a}^2}\right) \left[\begin{array}{c}1 \\ \bar{u_n} \\  \bar{u_t}\\ \dfrac{\bar{q^2}}{2}\end{array}\right]+\left(\dfrac{\gamma-1}{\gamma}\right)\bar{\rho}\left[\begin{array}{c} 0 \\  \Delta{}u_n\\ \Delta{}u_t \\ \Delta{}(q^2)/2\end{array}\right]
\end{equation}
This can be further expressed as
\begin{equation} \label{bdq-hllcps-3}
	B\Delta{}\boldsymbol{U}=\left(\Delta{}\rho-\dfrac{\Delta{}p}{\bar{a}^2}\right) \left[\begin{array}{c}1 \\ \bar{u}_n \\  \bar{u}_t\\ \dfrac{\bar{q}^2}{2}\end{array}\right]+\left(\dfrac{\gamma-1}{\gamma}\right)\bar{\rho}\Delta{}u_t\left[\begin{array}{c} 0 \\  0 \\ 1  \\ \bar{u}_t\end{array}\right] +\left(\dfrac{\gamma-1}{\gamma}\right)\bar{\rho}\Delta{}u_n\left[\begin{array}{c} 0 \\  1  \\ 0 \\ \bar{u_n}\end{array}\right]
\end{equation}
It can be seen from equation (\ref{bdq-hllcps-3}) that the anti-diffusion coefficients of the contact and shear wave for the HLL-CPS scheme are $\delta{}_2=1.0$ and $\delta{}_3=\frac{\gamma-1}{\gamma}$. The HLL-CPS scheme has additional anti-diffusion in the direction normal to the cell interface and the coefficient (which can be denoted by $\delta_n$) is equal to $\delta_3$. For the HLLEM scheme, the normal anti-diffusion coefficient $\delta_n$ that is present in the HLL-CPS scheme can be considered zero. The normal anti-diffusion present in the HLL-CPS scheme may appear unphysical. However, it will be shown later that anti-diffusion in the normal direction helps achieve accuracy for low-speed flows. Further, according to Dellacherie et al. \cite{della2}, anti-diffusion terms in the face normal direction are essential in all Godunov-type schemes to improve accuracy in low Mach number flows. It is also interesting to note that the anti-diffusion coefficient for the shear wave in the HLL-CPS scheme is $\delta_3=0.2875$ for a perfect gas with $\gamma=1.4$. Therefore, the ability of the HLL-CPS scheme to resolve the shear layers shall be better than the classical HLLE scheme but inferior to the HLLEM scheme. Overall, the HLL-CPS scheme can therefore be described as a perfect contact resolving but an approximate shear resolving scheme.

The numerical dissipation of the convective flux evaluation method of Mandal and Panwar \cite{mandal} can be assessed against the un-split HLL/HLLEM method. The HLL-CPS scheme can be expressed as a sum of central-difference terms and numerical dissipation as follows
\begin{equation}
\boldsymbol{F}_1(\boldsymbol{U}_R,\boldsymbol{U}_L)+\boldsymbol{F}_2(\boldsymbol{U}_R,\boldsymbol{U}_L) =\frac{1}{2}(\boldsymbol{F}_{1R}+\boldsymbol{F}_{1L}) + D_1 + \frac{1}{2}(\boldsymbol{F}_{2R}+\boldsymbol{F}_{2L}) + D_2
\end{equation}
The numerical dissipation of the convective flux of the HLL-CPS scheme with Zha-Bilgen splitting is
\begin{equation}
D_{1}=\boldsymbol{F}_1(\boldsymbol{U}_R,\boldsymbol{U}_L)-\frac{1}{2}(\boldsymbol{F}_{1R}+\boldsymbol{F}_{1L})
\end{equation}
After some mathematical manipulations, the numerical dissipation of the HLL-CPS convective flux can be written as
\begin{equation}\label{d1-hllcps}
D_1=\left\lbrace\begin{array}{l} \frac{1}{2}u_{nR}(\boldsymbol{U}_L-\boldsymbol{U}_R)+\frac{1}{2}M_{nL}\boldsymbol{U}_L(u_{nL}-u_{nR}) \hspace{1cm} \text{for} \hspace{1cm} u_n\ge 0 \\
\frac{1}{2}u_{nL}(\boldsymbol{U}_R-\boldsymbol{U}_L)+\frac{1}{2}M_{nR}\boldsymbol{U}_R(u_{nR}-u_{nL}) \hspace{1cm} \text{for} \hspace{1cm} u_n< 0 
\end{array} \right.
\end{equation}
The overall numerical dissipation of the HLLEM scheme, after expanding the flux into convective and pressure parts can be written as
\begin{equation}
D_{HLLEM}=\frac{1}{2}\dfrac{S_R+S_L}{S_R-S_L}(\boldsymbol{F}_{1L}-\boldsymbol{F}_{1R}+\boldsymbol{F}_{2L}-\boldsymbol{F}_{2R})+\frac{S_RS_L}{S_R-S_L}(\Delta{}\boldsymbol{U}-B\Delta{}\boldsymbol{U})
\end{equation}
The numerical dissipation of the HLLEM convective flux is
\begin{equation}
D_{1,HLLEM}=\frac{1}{2}\dfrac{S_R+S_L}{S_R-S_L}(\boldsymbol{F}_{1L}-\boldsymbol{F}_{1R})
\end{equation}
Considering $\dfrac{S_R+S_L}{S_R-S_L}\approx M_n$ for low-speed flows, the numerical dissipation of the HLLEM convective flux is 
\begin{equation}\label{d1-hllem}
D_{1,HLLEM}=\frac{1}{2}\tilde{M}_n\tilde{u}_n(\boldsymbol{U}_{L}-\boldsymbol{U}_{R})+\frac{1}{2}\tilde{M}_n\tilde{\boldsymbol{U}}(u_{nL}-u_{nR})
\end{equation}
Comparing equations (\ref{d1-hllcps}) and \ref{d1-hllem}), it can be seen that the numerical dissipation of the convective flux of the HLL-CPS scheme is marginally higher than that of the HLLEM scheme, which scales with the face normal Mach number. Nevertheless, the numerical dissipation of the convective flux of the HLL-CPS scheme approaches zero as velocity and Mach number approach zero and hence is suitable for low-speed flows. 
\subsection{Comparison of anti-diffusion coefficients}
The HLLEM scheme with the modification of Park and Kwon \cite{park} and the HLL-CPS scheme \cite{mandal} are capable of  resolving stationary contact discontinuity exactly but these schemes are incapable of resolving the flow features at low Mach numbers. The HLLEM scheme exhibits severe numerical shock instabilities while the HLL-CPS scheme is free from the carbuncle phenomenon at very high speeds \cite{mandal}, but exhibits mild instability for the supersonic corner test case. The HLLEM scheme is capable of resolving boundary layers accurately while the HLL-CPS scheme is incapable of accurate boundary layer resolution \cite{kita1, kita-slau2}. Thus, it is observed that the HLLEM and HLL-CPS schemes are similar in some aspects and different in some aspects. The difference in the behaviour of these schemes can be explained through the differences in their anti-diffusion coefficients. Thus, in the present work, a comparison of the anti-diffusion coefficients of these schemes is made for various flow problems and optimum anti-diffusion coefficients are evolved to make these schemes applicable for all-Mach number flows.
\subsubsection{One-dimensional flow problems}
\begin{figure}[H]
\begin{center}
\includegraphics[width=225pt]{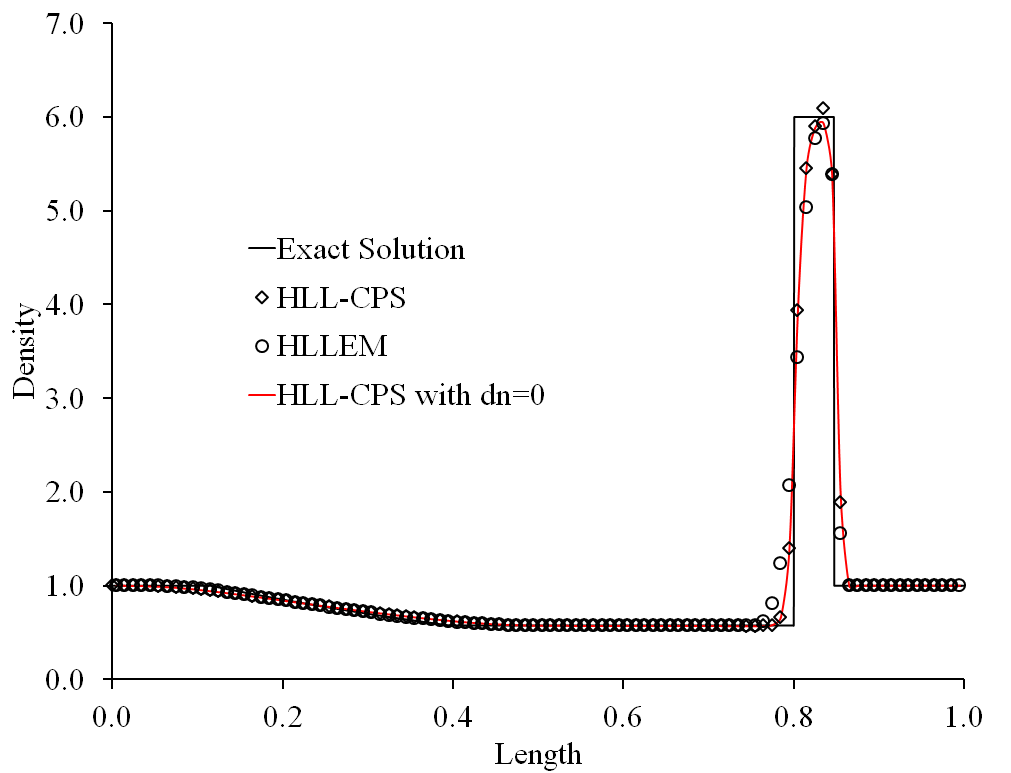} \includegraphics[width=225pt]{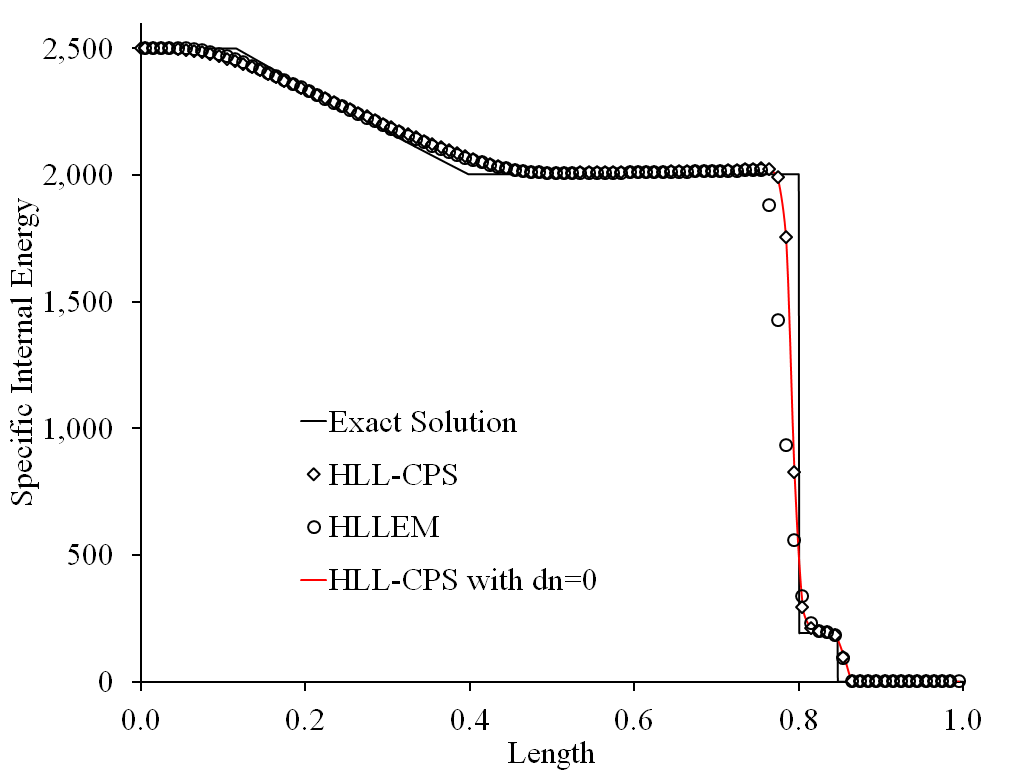}
\caption{Results of Sod shock tube problem}
\label{hllcps-hllem-testcase5}
\end{center}
\end{figure}

For one-dimensional flows, the shear wave is absent, and hence $\delta_3=0$. For the HLLEM scheme, the anti-diffusion coefficients are $\delta_2=\frac{\tilde{a}}{\tilde{a}+|\tilde{u}_n|}$ and $\delta_n=0$, while for the HLL-CPS scheme the anti-diffusion coefficients are $\delta_2=1$ and $\delta_n=\frac{\gamma-1}{\gamma}$. For isolated stationary contact discontinuity, the anti-diffusion coefficient for the contact wave $\delta_2$ becomes equal to unity for both the HLLEM and HLL-CPS schemes and hence both the schemes are capable of resolving stationary contact discontinuity. However, for a moving contact discontinuity, the HLL-CPS scheme  with $\delta_2=1$ will have higher anti-diffusion than the HLLEM scheme with $\delta_2=\frac{\tilde{a}}{\tilde{a}+|\tilde{u}_n|}$ and hence HLL-CPS scheme may resolve a moving contact discontinuity better than the HLLEM scheme. The HLLEM scheme shows monotone behaviour across a strong shock while the HLL-CPS scheme does not show perfectly monotone behaviour for all the shock tube problems \cite{mandal}. The results of the HLLEM and HLL-CPS schemes for a very severe shock tube problem with a strong shock running rightwards, a stationary contact and a left-running expansion wave are shown in Fig. \ref{hllcps-hllem-testcase5}. The initial conditions for the problem are $(\rho,\;u,\; p)_L=(1.0,\;-19.59745,\;1000.0)$ and $(\rho,\; u,\; p)_R=(1.0,\;-19.59475,\;0.01)$. The pre-shock Mach number here is about 165. It can be seen from the figure that  the HLLEM scheme is perfectly monotone while the HLL-CPS scheme is marginally non-monotone with a slight overshoot in density. The results of a modified HLL-CPS scheme with $\delta_n=0$ are also shown in the figure. It can be seen from the figure that the modified HLL-CPS scheme is perfectly monotone like the HLLEM scheme. It can be concluded that for obtaining a monotone solution across an infinitely strong shock, the value of $\delta_n$ must be equal to zero or nearly equal to zero across the shock. The figure also indicates that the moving contact discontinuity resolving ability of the HLL-CPS scheme is superior to the HLLEM scheme. 
\subsubsection{High-speed flows}
\begin{figure}[H]
\begin{center}
\includegraphics[width=250pt] {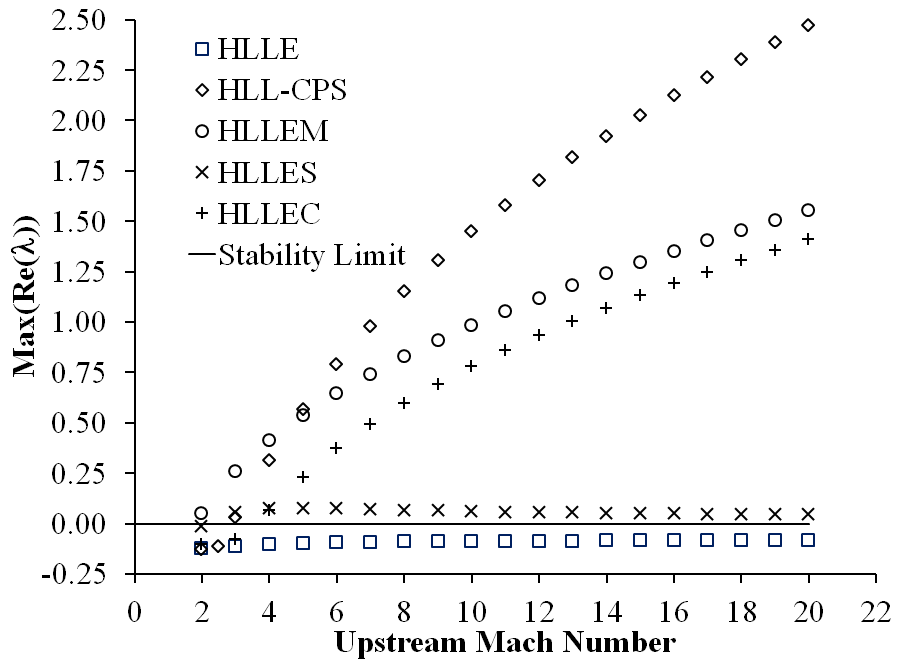}
\caption{Plot of maximum real eigenvalues vs upstream Mach number for the HLL-family schemes}
\label{eigenvalues-hll-family-schemes}
\end{center}
\end{figure}

The HLLEM scheme exhibit numerical shock instabilities like carbuncle phenomenon, odd-even decoupling, kinked Mach stem, etc., for high-speed flows. On the other hand, the HLLE scheme, which does not resolve the linearly degenerate waves, is free from the numerical shock instability problems. The HLL-CPS scheme \cite{mandal}, which resolves the contact discontinuity exactly, exhibits mild instability for the supersonic corner test case. The HLL-CPS scheme is susceptible to the carbuncle phenomena under conditions like highly stretched grids. The HLLES scheme proposed by Xie et al. \cite{xie-hllem}, which resolves the shear wave but excludes the contact wave, exhibits numerical shock instabilities. The HLLEC scheme, which resolves the contact wave but excludes the shear wave, has been shown to be free from numerical shock instabilities by Xie et al. \cite{xie-hllem}. However, Simon and Mandal \cite{san3} have shown that the HLLEC scheme is also prone to numerical shock instabilities. The HLLE+ scheme proposed by Park and Kwon \cite{park} has been shown to be free from numerical shock instabilities for viscous flows. However, Tramel et al. \cite{tramel} have observed that the HLLE+ scheme is also susceptible to numerical shock instabilities. The anti-diffusion coefficients of the various HLL-family schemes are shown in Table \ref{anti-diffusion-coeff}. It can be seen from the table that all the schemes that have non-zero anti-diffusion coefficients are susceptible to numerical shock instabilities. 
\begin{table}[H]
\caption{Anti-diffusion coefficients of the HLL-family schemes}
\label{anti-diffusion-coeff}
\begin{center}
\begin{doublespace}
\begin{tabular}{|c|c|c|c|c|}
\hline Sl No & Scheme & $\delta_2$& $\delta_3$ &  $\delta_n$\\
\hline 1 & HLLE & 0 & 0 & 0\\ 
\hline 2 &HLL-CPS & 1 & $\dfrac{\gamma}{\gamma-1}$ & $\dfrac{\gamma}{\gamma-1}$\\ 
\hline 3 & HLLEM & $\dfrac{\tilde{a}}{\hat{a}+|\tilde{u}_n|}$ & $\dfrac{\tilde{a}}{\hat{a}+|\tilde{u}_n|}$ & 0\\ 
\hline 4 &HLLEC & $\dfrac{\tilde{a}}{\tilde{a}+|\tilde{u}_n|}$ &  0  & 0\\ 
\hline 5 &HLLES & 0 &$\dfrac{\tilde{a}}{\tilde{a}+|\tilde{u}_n|}$  & 0\\ 
\hline 6 &HLLE+ & 0.5 (around  shock)&  0.5 (around shock) & 0\\ 
\hline
\end{tabular}
\end{doublespace}
\end{center}
\end{table}
The matrix stability analysis method of Dumbser et al. \cite{dumb-matrix} is a popular method of evaluating the stability of numerical schemes. The results of the matrix stability analysis of the various HLL-family schemes are shown in Fig. \ref{eigenvalues-hll-family-schemes}. The matrix stability analysis is carried out like Dumbser et al. \cite{dumb-matrix} and a `thin shock' is considered in the analysis. It can be seen from the figure that the maximum real eigenvalues of all the schemes with non-zero anti-diffusion coefficients are positive and suffer from numerical shock instability problems, while the HLLE scheme with anti-diffusion coefficients of zero has negative maximum real eigenvalues for all-Mach numbers and is free from the numerical shock instabilities. Therefore, it is felt that in the HLL-family schemes, the numerical shock instabilities can be reduced or avoided if the anti-diffusion coefficients for the shear and the contact waves are close to zero ($\delta_{2,3}\to 0$) around the shock. The anti-diffusion coefficient in the face normal direction $\delta_n$ is also required to be close to zero for obtaining a monotone solution across the shock, as shown in the previous section.
\subsubsection{Low Mach number flows}
For the HLLEM scheme, the coefficients $\delta_2$ and $\delta_3$ approach unity as the Mach number tends to zero and while $\delta_n$ is equal to zero irrespective of the Mach number. For the HLL-CPS scheme, the coefficients are $\delta_2=1$, $\delta_3=\delta_n=\frac{\gamma-1}{\gamma} $ irrespective of Mach number. The HLL-CPS and HLLEM schemes are incapable of resolving the flow features at low Mach numbers. The cause of the unphysical behaviour of the approximate Riemann solvers at low Mach numbers can be ascertained from the asymptotic analysis pioneered by Guillard and Viozat \cite{guillard1} and Rieper \cite{rieper1, rieper2}. It has been shown by Guillard and Viozat \cite{guillard1} that for the continuous Euler equations, the pressure is constant in space up to fluctuations of the order of the reference Mach number squared ($M_{*}^2$). Asymptotic analysis of the Roe scheme \cite{roe} carried out by Guillard and Viozat \cite{guillard1} reveals that the pressure fluctuation at very low Mach numbers is of the order of the reference Mach number $M_*$. The unphysical behaviour of the Roe scheme \cite{roe} at low Mach numbers has been attributed to the wrong order of pressure fluctuations. It has been shown by Rieper \cite{rieper1, rieper2} that for the Roe scheme the incorrect order of pressure fluctuation is due to the normal velocity jump terms in the normal momentum flux equation. Rieper has proposed a scaling of the normal velocity jump terms in the normal momentum flux equation of the Roe scheme by a local Mach number function \cite{rieper1, rieper2} to improve the accuracy of the scheme at low Mach numbers. This modification of Rieper can be expressed as a sum of the original term and anti-diffusion term as $f(M)\Delta{}u_n=\Delta{}u_n-(1-f(M))\Delta{}u_n$ and hence the anti-diffusion coefficient is $\delta_n=1-f(M)$. It is observed that this anti-diffusion coefficient in the normal momentum equation $\delta_n$ tends to unity as Mach tends to zero. Accurate results have been obtained by Rieper for inviscid flow around a cylinder with the above modification to the Roe scheme \cite{rieper1, rieper2}. Similarly, the low Mach modification of Godunov-type flux proposed by Dellacherie et al. \cite{della2} is
\begin{equation}
\boldsymbol{F}(\boldsymbol{U}_R,\boldsymbol{U}_L)=\boldsymbol{F}(\boldsymbol{U}_R,\boldsymbol{U}_L)^{G} -\dfrac{1}{2}(1-\theta) \rho{}a \left[0,\;\Delta{}u_n,\;0,\;0\right]^T 
\end{equation}
where $\boldsymbol{F}(\boldsymbol{U}_R,\boldsymbol{U}_L)^{G}$ is the original Godunov-type flux and $\theta$ is a Mach number-based function. Here, the anti-diffusion coefficient in the face normal momentum flux equation ($\delta_n=1-\theta$) approaches unity as the Mach number tends to zero.

Asymptotic analysis of the HLL-CPS scheme is therefore carried out like Guillard and Viozat \cite{guillard1} to identify the cause of non-physical behaviour of the scheme in low Mach number flows. The HLL-CPS scheme for the discrete Euler equations (\ref{fv},\ref{hll-cps-flux}) at low Mach numbers can be written as 
\begin{equation}\label{hllcps-lowspeed}
\begin{split}
& |\Omega|_i\dfrac{\partial{}U_i}{\partial{}t}+ \sum_{l\epsilon\upsilon(i)}\mathcal{T}^{-1}\left[u_{n,il}\left[\begin{array}{c} \rho\\ \rho{}u_n\\ \rho{}u_t \\ \rho{}e\end{array}\right]_K+\left[\begin{array}{c} 0 \\ p\\ 0 \\ pu_n\end{array}\right]_{il}+ \left(\dfrac{u_n}{2a}\right)_{il}\left[\begin{array}{c} 0 \\\Delta_{il}(p)\\ 0 \\ \Delta_{il}(pu_n)\end{array}\right]\right]\Delta{}s_{il}-\\ &\mathcal{T}^{-1}\left(\dfrac{u_n^2-a^2}{2a^3}\right)_{il}\left[\begin{array}{c}\Delta_{il}(p) \\ \Delta_{il}(pu_n) \\ \Delta_{il}(pu_t) \\ \dfrac{a^2_{il}\Delta_{il}(p)}{\gamma-1}+\dfrac{\Delta_{il}(pq^2)}{2} \\\end{array}\right]\Delta{}s_{il}=0
\end{split}
\end{equation}
The HLL-CPS scheme shown in (\ref{hllcps-lowspeed}) can be written in non-dimensional form, after rotation, as 
\begin{equation}
\begin{split}
& |\bar{\Omega}|_i\dfrac{\partial{}\bar{U}_i}{\partial{}\bar{t}}+  \sum_{l\epsilon\upsilon(i)}\left[\bar{u}_{n,il}\left[\begin{array}{c} \bar{\rho}\\ \bar{\rho}\bar{u}\\ \bar{\rho}\bar{v} \\ \bar{\rho}\bar{e}\end{array}\right]_K+\left[\begin{array}{c} 0 \\ \dfrac{\bar{p}}{M_*^2}n_x\\  \dfrac{\bar{p}}{M_*^2}n_y \\ \bar{p}\bar{u_n}\end{array}\right]_{il}+ \left(\dfrac{\bar{u}_nM_*}{2\bar{a}}\right)_{il}\left[\begin{array}{c} 0 \\ \left(\dfrac{n_x}{M_*^2}\right)_{il}\Delta_{il}(\bar{p})\\  \left(\dfrac{n_y}{M_*^2}\right)_{il} \Delta_{il}(\bar{p})\\ \Delta_{il}(\bar{p}\bar{u}_n)\end{array}\right]\right]\Delta{}s_{il}-  \\ &  \sum_{l\epsilon\upsilon(i)}\left[\left(\dfrac{\bar{u}_n^2M_*}{2\bar{a}^3}-\dfrac{1}{2\bar{a}M_*}\right)_{il}\left[\begin{array}{c}\Delta_{il}(\bar{p}) \\ \Delta_{il}(\bar{p}\bar{u}) \\ \Delta_{il}(\bar{p}\bar{v}) \\ \dfrac{\bar{a}^2_{il}\Delta_{il}(\bar{p})}{\gamma-1}+\dfrac{\Delta_{il}(\bar{p}\bar{q}^2)M_*^2}{2} \\\end{array}\right]\right]\Delta{}s_{il}=0
\end{split}
\end{equation}
Expanding the terms asymptotically and arranging the terms in the equal power of $M_*$, we obtain 
\begin{enumerate}
	\item{Order of $M_*^0$}
	\begin{enumerate}
		\item{Continuity Equation}
		\begin{equation}
		|\bar{\Omega}|\dfrac{\partial\bar{\rho}_{0i}}{\partial{}\bar{t}}+\sum_{l\epsilon\upsilon(i)}\left(\bar{u}_{n0,il}\bar{\rho}_{0K}+\dfrac{1}{2\bar{a}_{0il}}\Delta_{il}\bar{p}_1\right)\Delta{}s_{il}=0
		\end{equation}
		\item{ x-momentum Equation}
		\begin{equation}
		\begin{split}
		& |\bar{\Omega}|\dfrac{\partial(\bar{\rho}_0\bar{u}_0)_i}{\partial{}\bar{t}}+ \sum_{l\epsilon\upsilon(i)}\left(\bar{u}_{n0,il}(\bar{\rho}_0\bar{u}_0)_K+\left(\bar{p}_2n_x\right)_{il}+\left(\dfrac{\bar{u}_{n0}}{2\bar{a}_0}n_x\right)_{il}\Delta_{il}\bar{p}_1 \right)\Delta{}s_{il}\\ & +\sum_{l\epsilon\upsilon(i)}\left(\dfrac{1}{2\bar{a}_{0il}}\Delta_{il}(\bar{p}_1\bar{u}_0+\bar{p}_0\bar{u}_1)\right)\Delta{}s_{il}=0
		\end{split}
		\end{equation}
		\item{ y-momentum Equation}
		\begin{equation}
		\begin{split}
		& |\bar{\Omega}|\dfrac{\partial(\bar{\rho}_0\bar{v}_0)_i}{\partial{}\bar{t}}+ \sum_{l\epsilon\upsilon(i)}\left(\bar{u}_{n0,il}(\bar{\rho}_0\bar{v}_0)_K+\left(\bar{p}_2n_y\right)_{il}+\left(\dfrac{\bar{u}_{n0}}{2\bar{a}_0}n_y\right)_{il}\Delta_{il}\bar{p}_1 \right)\Delta{}s_{il}\\ & +\sum_{l\epsilon\upsilon(i)}\left(\dfrac{1}{2\bar{a}_{0il}}\Delta_{il}(\bar{p}_1\bar{v}_0+\bar{p}_0\bar{v}_1)\right)\Delta{}s_{il}=0
		\end{split}
		\end{equation}
	\end{enumerate}
	\item{Order of $M_*^{-1}$} 
	\begin{enumerate}
		\item{Continuity Equation}
		\begin{equation}\label{deltap0}
		\sum_{l\epsilon\upsilon(i)}\left(\dfrac{1}{2\bar{a}_{il}}\Delta_{il}\bar{p}_0\right)\Delta{}s_{il}=0
		\end{equation}
		From the above equation, we obtain $\sum_{l\epsilon\upsilon(i)}\Delta_{il}\bar{p}_0=0$ and hence $\bar{p}_0=$ constant for all $\mathbf{i}$.
		\item{ x-momentum Equation}
		\begin{equation}
		\sum_{l\epsilon\upsilon(i)}\left((\bar{p}_1n_x)_{il}+\left(\dfrac{\bar{u}_{n0}}{2\bar{a}_0}n_x\right)_{il}\Delta_{il}\bar{p}_0+\dfrac{1}{2\bar{a}_{0il}}\Delta_{il}(\bar{p}_0\bar{u}_0)\right)=0
		\end{equation}
		\begin{equation}
		\sum_{l\epsilon\upsilon(i)}\left((\bar{p}_1n_x)_{il}+\left(\dfrac{\bar{u}_{n0}}{2\bar{a}_0}n_x\right)_{il}\Delta_{il}\bar{p}_0+\dfrac{1}{2\bar{a}_{0il}}(\bar{p}_{0il}\Delta_{il}\bar{u}_{0}+\bar{u}_{0il}\Delta_{il}\bar{p}_0)\right)=0
		\end{equation}
		As per equation (\ref{deltap0}) $\sum_{l\epsilon\upsilon(i)}\Delta_{il}\bar{p}_0=0$. Therefore, 
		\begin{equation}\label{pnx}
		\sum_{l\epsilon\upsilon(i)}(\bar{p_1}n_x)_{il}=\sum_{l\epsilon\upsilon(i)}-\dfrac{1}{2}\dfrac{\bar{p}_{0,il}}{\bar{a}_{0il}}\Delta_{il}\bar{u}_0
		\end{equation}
		\item{ y-momentum Equation}
		\begin{equation}
		\sum_{l\epsilon\upsilon(i)}\left((\bar{p}_1n_y)_{il}+\left(\dfrac{\bar{u}_{n0}}{2\bar{a}_0}n_x\right)_{il}\Delta_{il}\bar{p}_0+\dfrac{1}{2\bar{a}_{0il}}\Delta_{il}(\bar{p}_0\bar{v}_0)\right)\Delta{}s_{il}=0
		\end{equation}
		As per equation (\ref{deltap0}) $\sum_{l\epsilon\upsilon(i)}\Delta_{il}\bar{p}_0=0$. Therefore,
		\begin{equation}\label{pny}
		\sum_{l\epsilon\upsilon(i)}\bar{p_1}n_y=-\sum_{l\epsilon\upsilon(i)}\dfrac{1}{2}\dfrac{\bar{p}_{0,il}}{\bar{a}_{0il}}\Delta_{il}\bar{v}_0
		\end{equation}	
	\end{enumerate}
\end{enumerate}

Rotating the x and y-momentum equation results of (\ref{pnx}, \ref{pny}) into the normal direction, we obtain
\begin{equation}
\sum_{l\epsilon\upsilon(i)}\bar{p}_1n_xn_x+\sum_{l\epsilon\upsilon(i)}\bar{p}_1n_yn_y=-\sum_{l\epsilon\upsilon(i)}\dfrac{1}{2}\dfrac{\bar{p}_{0,il}}{\bar{a}_{il}}\Delta_{il}\bar{u}_0n_x-\sum_{l\epsilon\upsilon(i)}\dfrac{1}{2}\dfrac{\bar{p}_{0,il}}{\bar{a}_{0il}}\Delta_{il}\bar{v}_0n_y
\end{equation}
\begin{equation} \label{hllcps-p1}
\sum_{l\epsilon\upsilon(i)}\bar{p}_1=-\sum_{l\epsilon\upsilon(i)}\dfrac{1}{2}\dfrac{\bar{p}_{0,il}}{\bar{a}_{il}}\Delta_{il}\bar{u}_{n0}
\end{equation}

Therefore, the HLL-CPS scheme permits  pressure fluctuation of the type $p(x,t)=p_0(t)+M_*p_1(x,t)$ and hence, the scheme will not be able to resolve the flow features at low Mach numbers. It can be seen from the above equation that that the normal velocity jump is the cause of the non-physical behaviour of the HLL-CPS scheme at low Mach numbers. Further, it is apparent that the normal velocity jump will vanish in the limit of Mach number zero if the anti-diffusion coefficient $\delta_n$ tends to unity as Mach number tends to zero. 

Therefore, for resolving the low Mach flow features, all the anti-diffusion coefficients ($\delta_2,\;\delta_3$ and $\delta_n$) of the HLLEM and HLL-CPS schemes should approach unity as the Mach number tends to zero ($\delta_2,\; \delta_3,\;\delta_n\to 1$ as $M \to 0$). It shall be proved later by asymptotic analysis that pressure fluctuations similar to the continuous Euler equations are obtained in the HLLEM and HLL-CPS schemes if the anti-diffusion coefficients, $\delta_2,\;\delta_3,\;\delta_n\to 1$ as $M \to 0$.
\subsubsection{Boundary layer resolution} \label{bl-profile}
\begin{figure}[H]
	\begin{center}	
	\includegraphics[width=220pt]{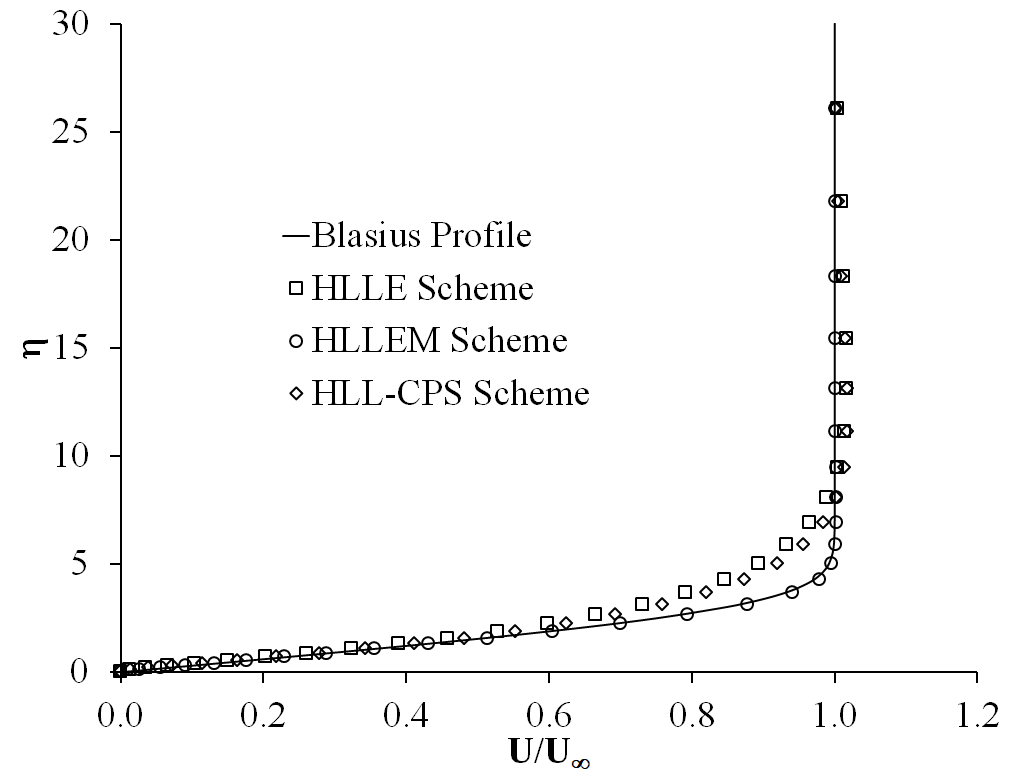}
	\caption{Boundary layer profile for $M_{\infty}=0.20$ laminar flow over a flat plate computed by the HLLE, HLL-CPS and HLLEM schemes}
	\label{bl-hll-family}
	\end{center}
\end{figure}

The Roe scheme, the HLLE+ scheme, and the HLLEM scheme with Park modification are capable of resolving the wall boundary layer accurately while the HLLE scheme is incapable of accurate boundary layer resolution. The flat plate boundary layer profile of the HLL-CPS scheme has been shown to lie between the HLLEM and HLLE schemes \cite{kita1, kita-slau2}. For the HLLEM scheme, the anti-diffusion coefficients are $\delta_2=\delta_3=\frac{\tilde{a}}{\tilde{a}+|\tilde{u}_n|}$ and $\delta_n=0$. Across the flat plate boundary layer, the Roe-averaged normal velocity $|\tilde{u}_n|$, which is also the approximate speed of the contact discontinuity, is almost equal to zero, and hence the anti-diffusion coefficients $\delta_2$ and $\delta_3$ are equal to (or almost equal to) unity in the HLLEM and Roe schemes. For the HLL-CPS scheme the anti-dissipation coefficients are $\delta_2=1$ and $\delta_3=\delta_n=0.2857$ for a fluid with $\gamma=1.4$ like air. The flat plate boundary layer profiles of the HLL-CPS and HLLEM schemes for a free-stream Mach number of 0.20 are shown in Fig. \ref{bl-hll-family}. The results are shown after 50,000 iterations with a CFL number of 0.50 on a grid size of $81\times33$. The third-order computations are carried out with the MUSCL approach \cite{van-leer-muscl} ($\kappa=1/3$) without any limiters. It can be seen from the figure that the HLLEM scheme resolves the boundary layer accurately as its anti-diffusion coefficient for the contact and shear waves ($\delta_2$ and $\delta_3$) are nearly equal to unity across the boundary layer. Therefore, it can be said that for accurate boundary layer resolution, the anti-diffusion coefficient for the contact and shear waves, $\delta_2$ and $\delta_3$ should be equal or close to unity inside the boundary layer like the Roe/HLLE+/HLLEM scheme.
\subsection{Improvement of HLLEM and HLL-CPS Schemes}
The following guidelines for the anti-diffusion coefficients can now be defined for obtaining all-speed HLL-CPS and HLLEM scheme:
\begin{enumerate}
 \item For obtaining a monotone solution across an infinitely strong shock, $\delta_n$ should be equal or close to zero as in the HLLE and HLLEM schemes.
 \item For avoiding  numerical shock instabilities at high Mach numbers, the anti-diffusion coefficients should be close to zero in the vicinity of the shock ($\delta_{2,3,n} \to 0$) like the classical HLLE scheme.
 \item For obtaining accurate solutions at low Mach numbers all three anti-diffusion coefficients should tend to unity as the Mach number tends to zero ($\delta_{2,3,n}\to 1$ as Mach $\to 0$).
 \item For the accurate resolution of boundary layers, the anti-diffusion coefficients of the contact and shear waves should be equal or close to unity ($\delta_2,\delta_3 \approx 1$), as in the Roe/HLLE+/HLLEM schemes.
\end{enumerate}
Based on the above guidelines, the anti-diffusion coefficients are improved to make the HLL-CPS scheme a genuinely all-Mach number scheme. The anti-diffusion coefficient in the HLL-CPS scheme for the contact wave ($\delta_2=1$) is excessive while the anti-diffusion coefficient of the HLLEM scheme for the contact wave is more appropriate. Kemm \cite{kemm} and Simon and Mandal \cite{san3} have shown that the exact resolution of the contact wave also contributes to numerical shock instabilities in the HLLEM scheme. Therefore, the anti-diffusion coefficients of both the contact and shear waves in the original HLLEM scheme are modified with a multi-dimensional pressure function like Simon and Mandal \cite{san3}. The anti-diffusion coefficients for the contact and shear waves in the improved HLL-CPS schemes can be written as
\begin{equation} \label{delta23-new}
\delta_{2,new}=\delta_{3,new}=\dfrac{\tilde{a}}{\tilde{a}+|\tilde{u}_n|}f_{p1}
\end{equation}
where $f_{p1}$ is a multi-dimensional pressure function. The pressure function $f_{p1}$ at the interfaces $(i+\frac{1}{2},j)$ and $(i,j+\frac{1}{2})$ are taken as 
\begin{equation} \label{fp-2d}
	\begin{array}{l} 
	f_{p1,i+\frac{1}{2},j}=min\left(f_{p,i+\frac{1}{2},j},f_{p,i,j+\frac{1}{2}},f_{p,i+1,j+\frac{1}{2}},f_{p,i,j-\frac{1}{2}},f_{p,i+1,j-\frac{1}{2}}\right) \\
	f_{p1,i,j+\frac{1}{2}}=min\left(f_{p,i,j+\frac{1}{2}},f_{p,i-\frac{1}{2},j},f_{p,i-\frac{1}{2},j+1},f_{p,i+\frac{1}{2},j},f_{p,i+\frac{1}{2},j+1}\right)
	\end{array}
\end{equation}
where
\begin{equation} \label{fp}
f_{p}=min\left(\left(\dfrac{p_L}{p_R}\right), \left(\dfrac{p_R}{p_L}\right)\right)^3
\end{equation}
where $p_L$ and $p_R$ are the static pressure to the right and left of the interface. 
Alternate, less dissipative anti-diffusion coefficients of the contact and shear wave for the improved HLL-CPS scheme are $\delta_{2,new, HLL-CPS}=\delta_{3,new, HLL-CPS}=fp$. However, these anti-diffusion coefficients may not meet the stability requirements of the HLLEM scheme \cite{einf2}. The anti-diffusion coefficient in the face normal direction is taken as a function of the local Mach number like Rieper \cite{rieper1, rieper2} and Dellacherie et al \cite{della2}. However, mild oscillations and non-monotone behaviour are observed with the anti-diffusion coefficient based only on the local Mach number \cite{ossw}. The anti-diffusion coefficient is scaled by the pressure-based function in order to obtain monotone solutions and suppress numerical shock instabilities. The face normal anti-diffusion coefficient can be defined as
\begin{equation} \label{deltan-new}
\delta_{n,new}=[1-f(M)]f_{p1}
\end{equation}
where $f_{p1}$ is the pressure function defined in equations (\ref{fp-2d}, \ref{fp}), $f(M)=min(max(M_R,\;M_L),\;1)$, $M_R$ and $M_L$ are the local Mach numbers to the right and left of the interface.

The improved anti-diffusion term $B\Delta{}\boldsymbol{U}$ for the HLLEM and HLL-CPS schemes can be written as 
\begin{equation} \label{bdq-new}
B\Delta{}\boldsymbol{U}_{new}=\delta_{2,new}\tilde{\alpha}_2\tilde{R}_2+\delta_{3,new}\tilde{\alpha}_3\tilde{R}_3 +\delta_{n,new}\tilde{\rho}\Delta{}u_n\left[0,\;1,\;0,\;\tilde{u}_n\right]^T
\end{equation}
where $\delta_{2,new}$ and $\delta_{3,new}$ are the anti-diffusion coefficients of the contact and shear waves defined in equation (\ref{delta23-new}), $\tilde{\alpha}_2=\Delta{}\rho-\dfrac{\Delta{}p}{\tilde{a}^2}$ and $\tilde{\alpha}_3=\tilde{\rho}\Delta{}u_t$ are the wave strength of the contact and shear waves, $\tilde{R}_2=\left[1,\; \tilde{u}_n,\;\tilde{u}_t,\;\frac{1}{2}({u}_u^2+\tilde{u}_t^2)\right]^T$ and $\tilde{R}_3=[0,\;0,\;1,\;\tilde{u}_t]^T$ are the right eigenvectors of the contact and shear wave in the flux Jacobian matrix, $\delta_{n,new}$ is the anti-diffusion coefficient for resolution of low Mach flow features defined in equation (\ref{deltan-new}). The new improved $B\Delta{}\boldsymbol{U}$ can now be used in equation (\ref{hll-cps-flux2}) for the HLL-CPS scheme and in equation (\ref{hll-type-flux}) for the HLLEM scheme. The improved schemes can be named HLL-CPS-FP and HLLEM-FP respectively. In the proposed HLL-CPS-FP scheme the anti-diffusion coefficients are modified from their original values of $\delta_2=1$, $\delta_3=\delta_n=\frac{\gamma-1}{\gamma}$ to $\delta_{2,new}=\frac{\tilde{a}}{\tilde{a}+|\tilde{u}_n|}f_{p1}$ or $\delta_{2,new}=f_{p1}$, $\delta_{3,new}=\frac{\tilde{a}}{\tilde{a}+|\tilde{u}_n|}f_{p1}$ or $\delta_{3,new}=f{p1}$, and $\delta_{n,new}=(1-f(M))f_{p1}$. In the proposed HLLEM-FP scheme, the anti-diffusion coefficients for the contact and shear wave are scaled with a pressure function and additional anti-diffusion is introduced in the normal direction. The additional anti-diffusion in the normal direction shall reduce dissipation due to the normal velocity jump and improve performance in low Mach number flows. It may be noted here that, in the HLL-CPS scheme, arithmetic averaging is used instead of Roe averaging and hence the Roe averaged terms in the above equation can be replaced by arithmetic average for the HLL-CPS scheme. The proposed HLLEM scheme can also be expressed as a sum of the original HLLEM scheme and additional terms. Since the HLLC scheme is similar to the HLLEM scheme, the additional terms can also be used with the original HLLC scheme to improve its shock stability and its ability to resolve low Mach flow features. Therefore, the all-Mach number HLLEM and HLLC schemes, in terms of the original schemes and additional terms, can be written as
\begin{equation} \label{hllc-allmach}
\boldsymbol{F}(\boldsymbol{U_R, U_L})_{All-Mach}=\boldsymbol{F}(\boldsymbol{U_R, U_L})+\dfrac{S_RS_L}{S_R-S_L}\left((f_{p1}-1)\sum_{k=2}^3(\delta_k\tilde{\alpha}_k\tilde{R}_k)+\delta_{n,new}\tilde{\rho}\Delta{}u_n(0,\;1,\;0,\;\tilde{u}_n)^T\right)
\end{equation}
where $\boldsymbol{F}(\boldsymbol{U_R, U_L})$ is the original HLLEM/HLLC flux, $f_{p1}$ is the pressure function defined in equations (\ref{fp} and \ref{fp-2d}), $\delta_2=\delta_3=\frac{\tilde{a}}{\tilde{a}+|\tilde{u}_n|}$ are the anti-diffusion coefficients of the contact and shear waves in the original HLLEM scheme, $\delta_{n,new}$ is the anti-diffusion coefficient for resolving low Mach flow features defined in equation (\ref{deltan-new}). It can be seen from the equation that across a strong shock, $f_{p1}\to 0$ and hence diffusion of the contact and shear wave is introduced in the vicinity of the shock to suppress numerical shock instabilities. In regions of smooth flows, $f_{p1}\to 1$ and hence the contact and shear wave resolution capability of the original HLLEM and HLLC schemes is retained.

In the following sections, the suitability of the proposed schemes for low-speed flows is demonstrated through asymptotic analysis while the stability of the proposed HLL-CPS-FP and HLLEM-FP schemes is demonstrated through linear perturbation and matrix stability analyses.
\subsection{Analysis of proposed schemes}
Asymptotic analysis, linear perturbation and matrix stability analyses of the proposed HLL-CPS and HLLEM schemes are carried out to demonstrate that the proposed schemes are capable of resolving low Mach flow features and are robust against numerical shock instabilities.
\subsubsection{Asymptotic analysis}
Asymptotic analysis of the proposed HLL-CPS and HLLEM schemes is carried out like Guillard and Viozat \cite{guillard1}. Asymptotic analysis of the proposed HLL-CPS-FP scheme is placed in Appendix I while asymptotic analysis of the proposed HLLEM-FP scheme is placed in Appendix II. It can be seen that the pressure fluctuation in the proposed HLL-CPS-FP and HLLEM-FP schemes are of the type $p(x,t)=p_0(t)+M_{*}^2p_2(x,t)$. Hence, these schemes shall be capable of resolving the flow features at low Mach numbers.
\subsubsection{Linear perturbation analysis}
Linear perturbation analysis of the proposed schemes is carried out like Quirk \cite{quirk} and Pandolfi-D'Ambrosio \cite{pand}. A saw-tooth type perturbation is introduced into the primitive variables like density, shear velocity and pressure. The evolution of the density, shear velocity, and pressure perturbations with time in  the proposed schemes is obtained like Quirk \cite{quirk} and Pandolfi-D'Ambrosio \cite{pand}. The normalized flow properties are described as follows
\begin{equation}
\rho=1\pm \hat{\rho}, \hspace{1cm} u=u_0\pm \hat{u}, \hspace{1 cm} v=0, \hspace{1cm}  p=1\pm\hat{p}
\end{equation}
where $\hat{\rho}$, $\hat{u}$ and $\hat{p}$ are the density, shear velocity, and pressure perturbations respectively. It may be noted that here $u$ is the shear velocity and $v$ is the face normal velocity.
The evolution of the density and pressure perturbations for the proposed schemes along with the original schemes is shown in Table \ref{pert-anal-new-schemes}. In the HLLEM scheme, the density perturbations are fed by the pressure perturbations while the shear velocity perturbations remain unchanged. In the HLL-CPS scheme, the density perturbations are fed by the pressure perturbations, just like the HLLEM scheme, but the shear velocity perturbations are damped in the HLL-CPS scheme. In the proposed HLL-CPS-FP and HLLEM-FP schemes, the pressure switch $f_{p1}$ along with the density and pressure perturbations influence the evolution of the density perturbations. Across a strong shock, the value of $f_{p1}$ becomes close to zero and the density perturbations are damped like the HLLE scheme. The shear velocity perturbations in the proposed HLL-CPS-FP and HLLEM-FP schemes are influenced by the pressure switch $f_{p1}$. Across a strong shock, the value of $f_{p1}$ becomes close to zero and the shear velocity perturbations are damped like the HLLE scheme. Since the density and shear velocity perturbations are damped in presence of strong shock in the HLL-CPS-FP and HLLEM-FP schemes, these schemes are likely to be free from numerical shock instabilities like the classical HLLE scheme. 
\begin{table}[H]
\caption{Result of linear perturbation analysis of the proposed HLL-CPS and HLLEM schemes}
\label{pert-anal-new-schemes}
\begin{center}
\begin{doublespace}
\begin{tabular}{|c|c|c|c|c|}
\hline Sl No & Scheme & $\hat{\rho}^{n+1}=$ & $\hat{u}^{n+1}=$ &  $\hat{p}^{n+1}=$\\
\hline 1 & HLL-CPS & $\hat{\rho}^n-\dfrac{2\nu}{\gamma}\hat{p}^n$ & $\hat{u}^n(1-\dfrac{2\nu}{\gamma})$ & $\hat{p}^n(1-2\nu)$\\
\hline 2 &HLL-CPS-FP & $\hat{\rho}^n(1-2\nu{}(1-f_{p1}))-\dfrac{2\nu}{\gamma}f_p\hat{p}^n$ & $\hat{u}^n(1-2\nu{}(1-f_{p1}))$ & $\hat{p}^n(1-2\nu)$\\
\hline 3 & HLLEM & $\hat{\rho}^n-\dfrac{2\nu}{\gamma}\hat{p}^n$ & $\hat{u}^n$ & $\hat{p}^n(1-2\nu)$\\
\hline 4 &HLLEM-FP & $\hat{\rho}^n(1-2\nu{}(1-f_{p1}))-\dfrac{2\nu}{\gamma}f_p\hat{p}^n$ & $\hat{u}^n(1-2\nu{}(1-f_{p1}))$ & $\hat{p}^n(1-2\nu)$\\
\hline
\end{tabular}
\end{doublespace}
\end{center}
\end{table}
\subsubsection{Matrix stability analysis}
Matrix stability analysis is carried out  for the proposed schemes like Dumbser et al. \cite{dumb-matrix}. A 2D computational domain $[0,1] \times [0,1]$ is considered and is discretized into $11 \times 11$ grid points. The raw state is a steady normal shock wave with some perturbations which can be expressed as
\begin{equation}
\mathbf{U}_\mathbf{i}=\mathbf{U}^0_\mathbf{i}+\delta{}\mathbf{U}_\mathbf{i}
\end{equation}
where $\mathbf{i}=(i,j)$ is the index of the cell, $\mathbf{U}^0_\mathbf{i}$ is the solution of a steady shock, and $\delta{}\mathbf{U}_\mathbf{i}$ is a small numerical random perturbation. Substituting the expression into finite volume formulation and after some evolution, the perturbations can be expressed as
\begin{equation}
\left(\begin{array}{c} \delta{}\mathbf{U}_1 \\ . \\ . \\ . \\ \delta{}\mathbf{U}_M\end{array}\right)=exp^{\mathbf{S}t}\left(\begin{array}{c} \delta{}\mathbf{U}_1 \\ . \\ . \\ . \\ \delta{}\mathbf{U}_M\end{array}\right)_0
\end{equation}
where $\mathbf{S}$ is the stability matrix based on the Riemann solver. The perturbations will remain bounded if the maximum of the real part of the eigenvalues of $\mathbf{S}$ is non-positive, i.e.,
\begin{equation}
max(Re(\lambda(\mathbf{S})))\le 0
\end{equation}
The initial data is provided by the exact Rankine-Hugoniot solution in the x-direction. The upstream and downstream states of the primitive variables are
\begin{equation}
\begin{array}{l}
W_L=(\rho, u, v, p)_L=\left(1, 1, 0, \dfrac{1}{\gamma{}M^2_a}\right), x< 0.5 \\ 
W_R=(\rho, u, v, p)_R=\left(f(M_a),\dfrac{1}{f(M_a)}, 0, \dfrac{g(M_a)}{\gamma{}M^2_a}\right), x>0.5
\end{array}
\end{equation}
where $M_a$ is the upstream Mach number,
$f(M_a)=\left(\dfrac{2}{\gamma+1}\dfrac{1}{M^2_a}+\dfrac{\gamma-1}{\gamma+1}\right)^{-1}$,  
$g(M_a)=\left(\dfrac{2\gamma}{\gamma+1}{M^2_a}-\dfrac{\gamma-1}{\gamma+1}\right)$.  

\begin{figure}[H]
	\begin{center}
	\includegraphics[width=220pt]{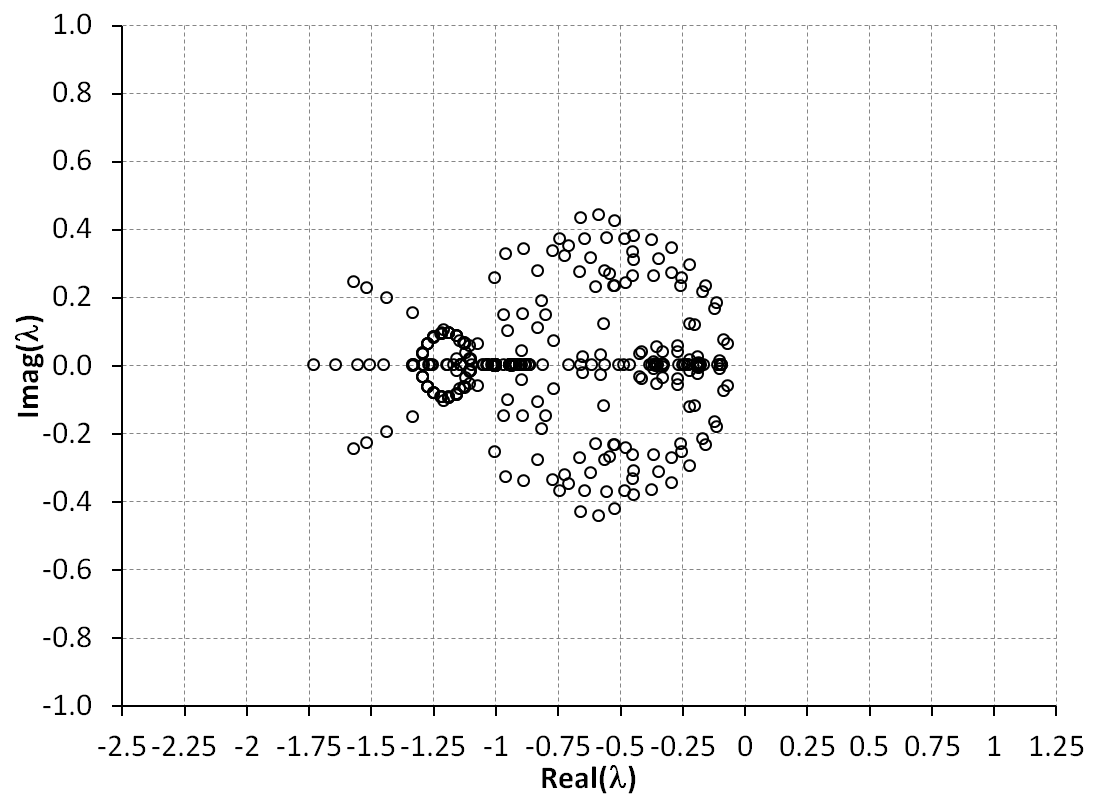}  \includegraphics[width=220pt]{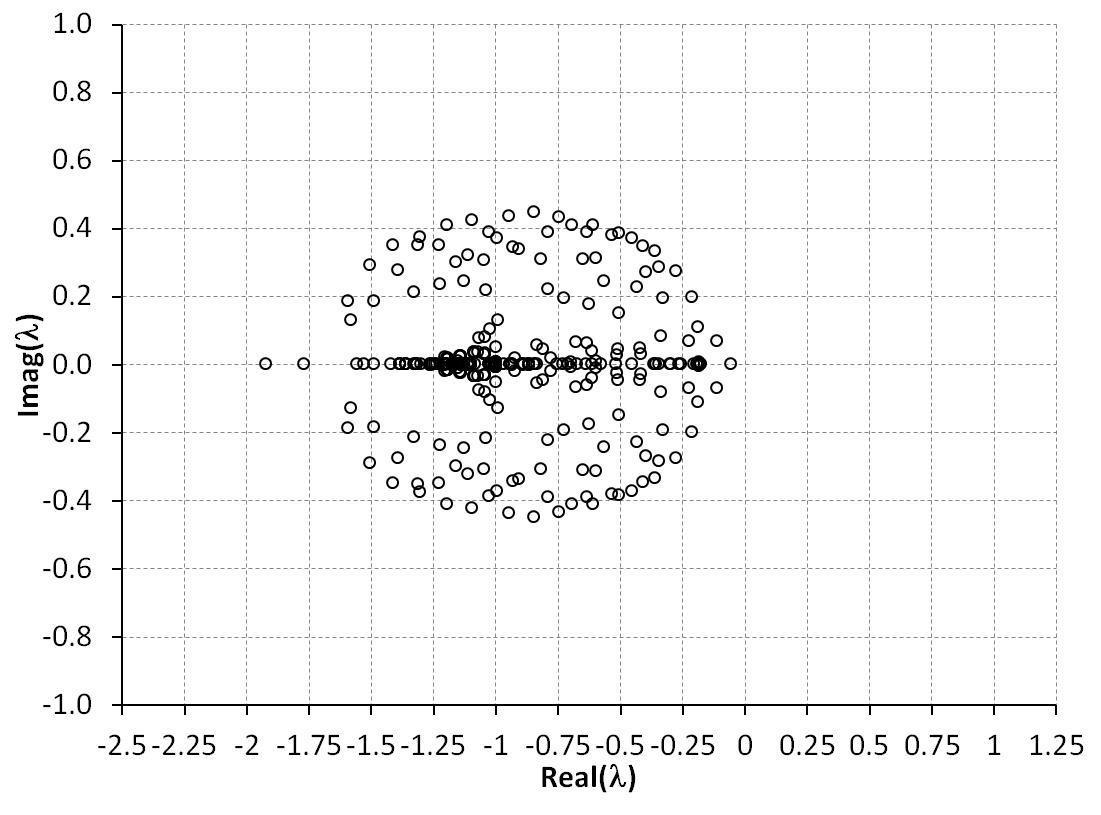} \\
	(a) HLLEM-FP Scheme \hspace{5cm} (b) HLL-CPS-FP Scheme
	\caption{Distribution of the eigenvalues of $\mathbf{S}$ in the complex plane for the proposed HLLEM-FP and HLL-CPS-FP Schemes for an upstream  Mach number of 7.0}
	\label{hllemfp-eigenvalue-m7}
	\end{center}
\end{figure}

\begin{figure}[H]
	\begin{center}
	\includegraphics[width=220pt]{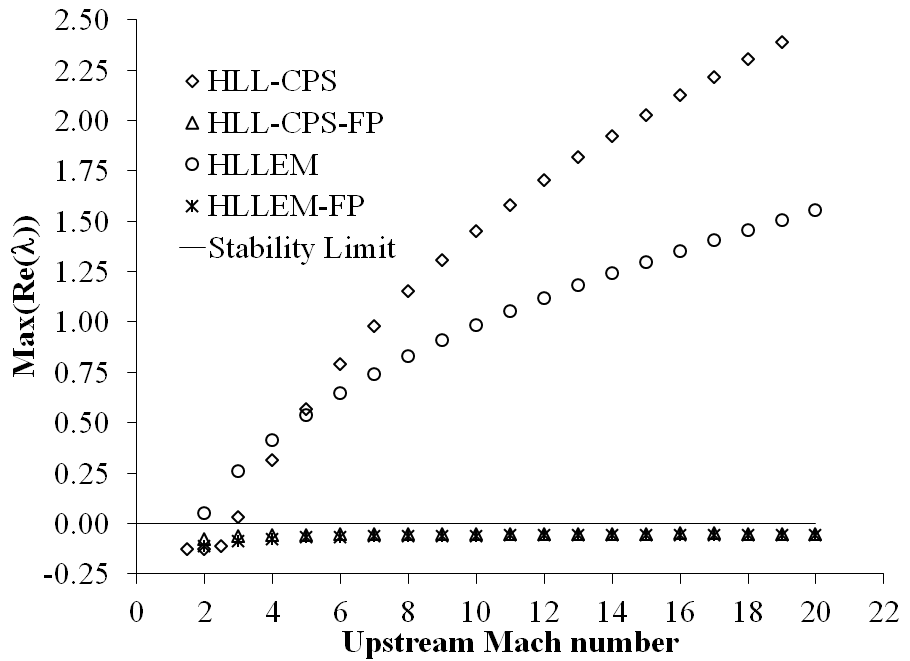} 
	\caption{Plot of maximum real eigenvalue vs upstream Mach number for the proposed HLLEM-FP and HLL-CPS-FP schemes}
	\label{max-real-ev-2p}
	\end{center}
\end{figure}

A random perturbation of $10^{-6}$ is introduced to all the conserved variables of all the grid cells. The matrix stability analysis can be carried out with either a `thin shock' comprising only upstream and downstream values or with a `shock structure' comprising an intermediate state. In the present work, the analysis is carried out for a `thin' shock. The plot of the eigenvalues of $\mathbf{S}$ in the complex plane for a typical upstream Mach number of 7.0 is shown in Fig. \ref{hllemfp-eigenvalue-m7} for the proposed HLLEM-FP and HLL-CPS-FP schemes. It can be seen from the figure that all the real eigenvalues are negative. The plot of the maximum real eigenvalues for different upstream Mach numbers is shown in Fig. \ref{max-real-ev-2p}. It can be seen from the figure that the maximum real eigenvalue of the proposed HLLEM and HLL-CPS schemes are less than zero for all the upstream Mach numbers considered in the analysis. Hence, the proposed HLLEM-FP and HLL-CPS-FP schemes are expected to be stable for high-speed flows.
\section {Results}
The performance of the proposed HLL-CPS-FP and HLLEM-FP schemes is evaluated by solving the problems involving planar shock, double Mach reflection, forward-facing step, blunt body, expansion corner, low-speed cylinder and flat plate boundary layer. The performance of the proposed schemes is compared with the original HLL-CPS and HLLEM schemes. First-order results are shown for the inviscid, high-speed test cases since the numerical shock instabilities are prominently visible in first-order schemes. However, second-order results are shown for a few test cases to demonstrate the ability of the proposed schemes to resolve the flow features. The results for the viscous test cases are also second-order accurate.
\subsection{Planar shock problem}
The problem consists of a Mach 6.0 shock wave propagating down in a rectangular channel. The domain consists of 800$\times{}$20 cells. The centreline in the y direction (11$^{th}$ grid line) is perturbed to promote odd-even decoupling along the length of the shock as	
$y_{i,11}=\left\lbrace\begin{array}{c}\ y_{i,11}+0.001 \hspace{2mm}\text{if  i is even}\\ 
 y_{i,11}-0.001 \hspace{2mm}\text{if  i is odd}\end{array}\right.$. 
 
\begin{figure}[H]
	\begin{center}
		\includegraphics[width=220pt]{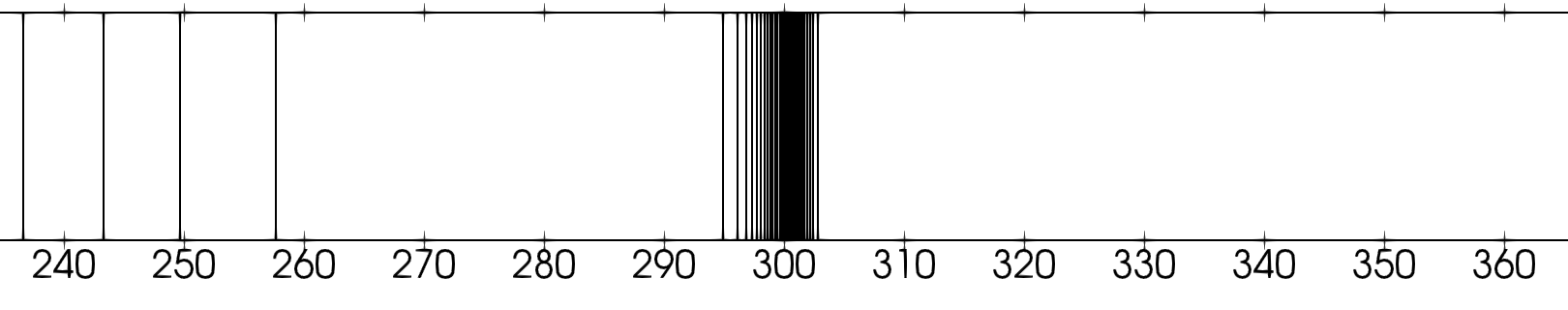} 
		\includegraphics[width=220pt]{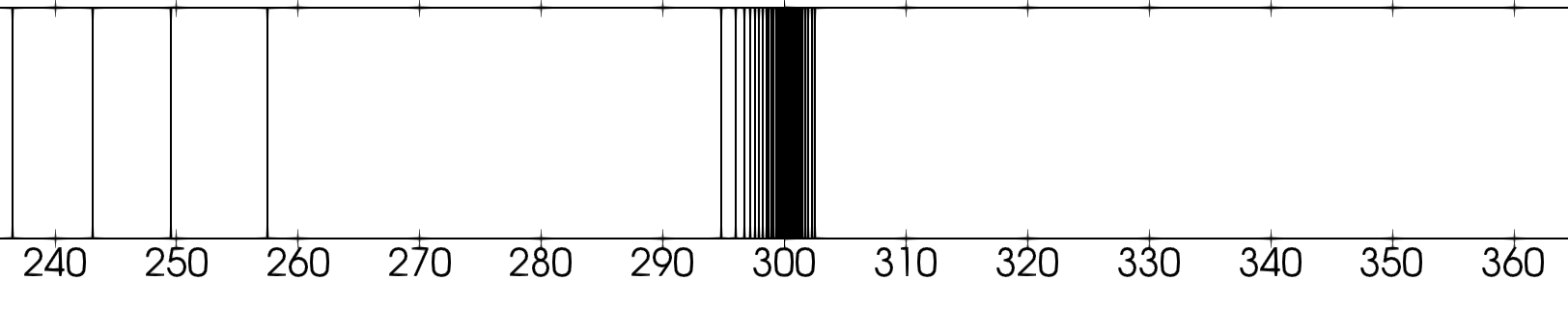}\\ 
		(a) HLL-CPS Scheme \hspace{5cm} (b) HLL-CPS-FP Scheme\\ 
		\includegraphics[width=220pt]{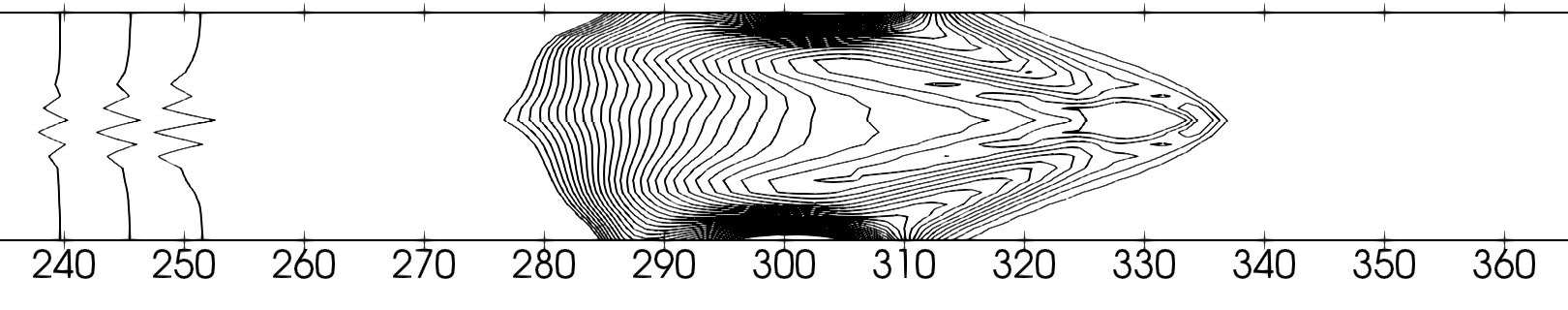} 
		\includegraphics[width=220pt]{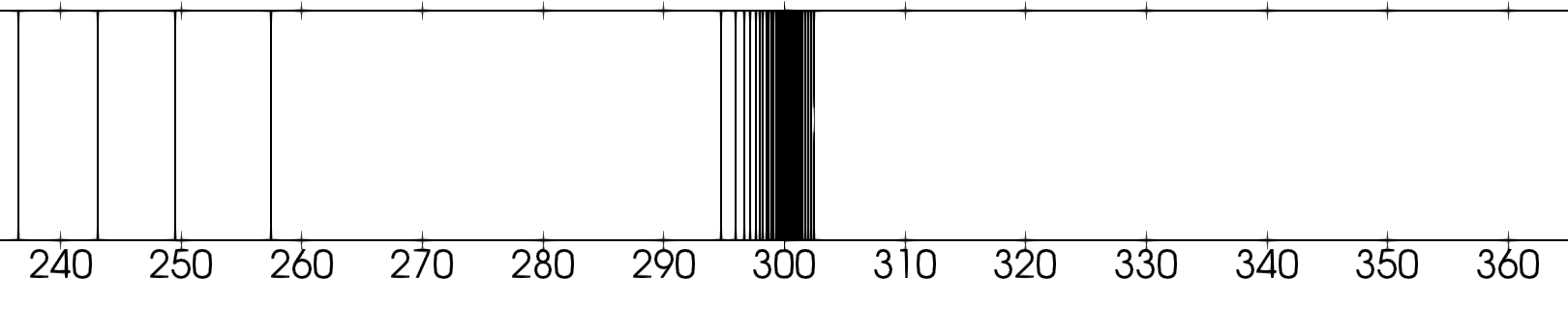} \\ 
		(c) HLLEM Scheme \hspace{5cm} (d)  HLLEM-FP Scheme 
		\caption{Density contour plots for the $M_{\infty}=6$ Planar shock problem computed by the HLL-CPS, HLL-CPS-FP, HLLEM and HLLEM-FP schemes. The results are shown at time t=55 units.}
	\label{planar-shock-hllcps-hllem}
	\end{center}
\end{figure}

The domain is initialized with $\rho=1.4,\;p=1.0,\;u=0$, and $v=0$. Post-shock values are imposed at the inlet and zero gradient boundary conditions are imposed at the exit. Solid wall boundary conditions are imposed at the top and bottom. The density contour plots for the original and proposed HLL-CPS and HLLEM schemes at time $t=55$ are shown in Fig.\ref{planar-shock-hllcps-hllem} and 30 contour levels ranging from 16 to 7.0 are drawn. It can be seen from the figure that odd-even decoupling is absent in the proposed HLL-CPS-FP scheme, just like the original HLL-CPS scheme. It can be seen from the figure that the original HLLEM scheme show odd-even decoupling while the proposed HLLEM-FP scheme is free from the odd-even decoupling problem.
\subsection{Double Mach reflection problem}
\begin{figure}[H]
	\begin{center}
		\includegraphics[width=220pt]{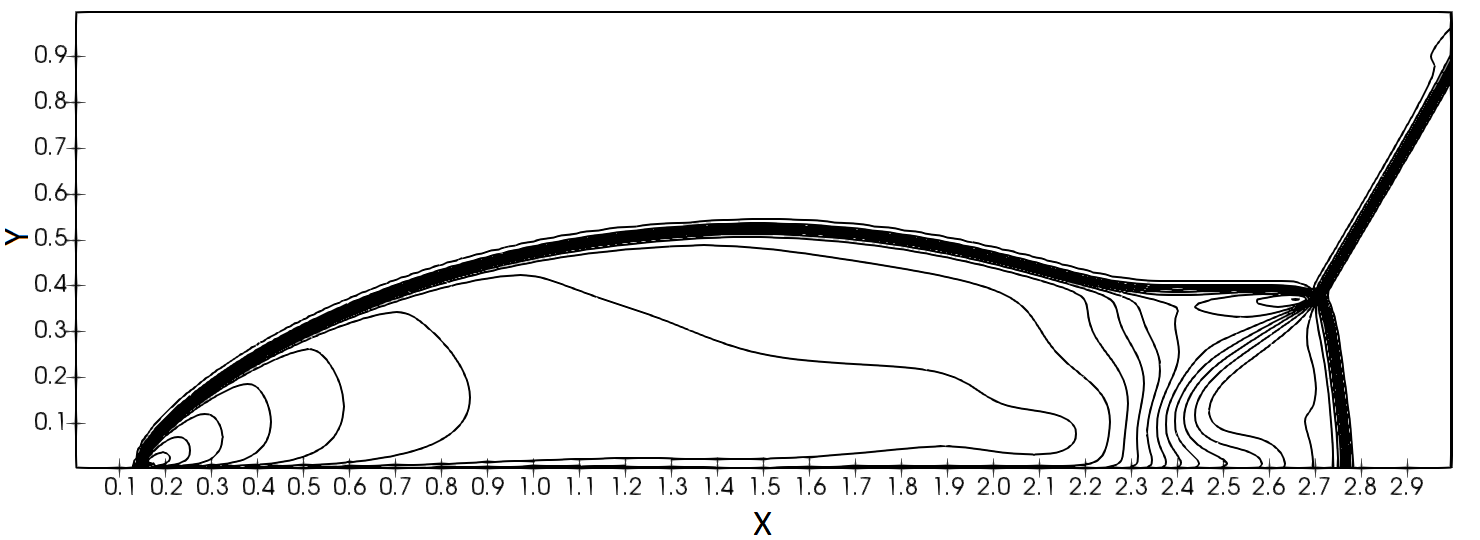} \includegraphics[width=220pt]{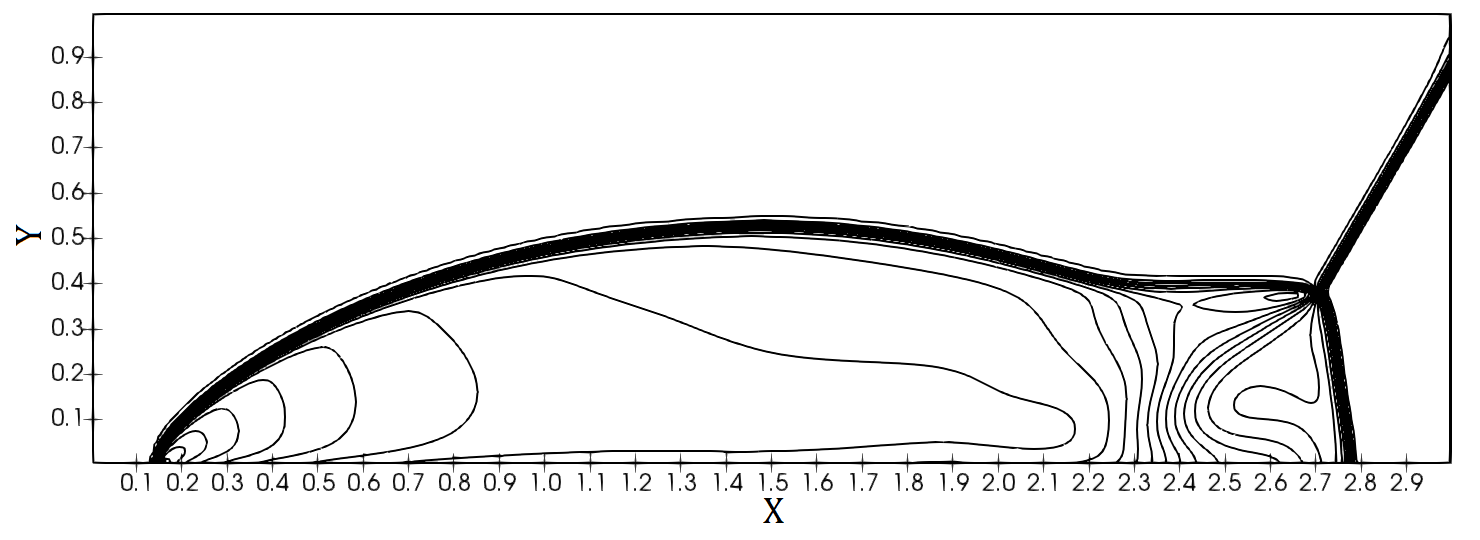}\\ 
		(a) HLL-CPS \hspace{5cm} (b) HLL-CPS-FP \\
		\includegraphics[width=220pt]{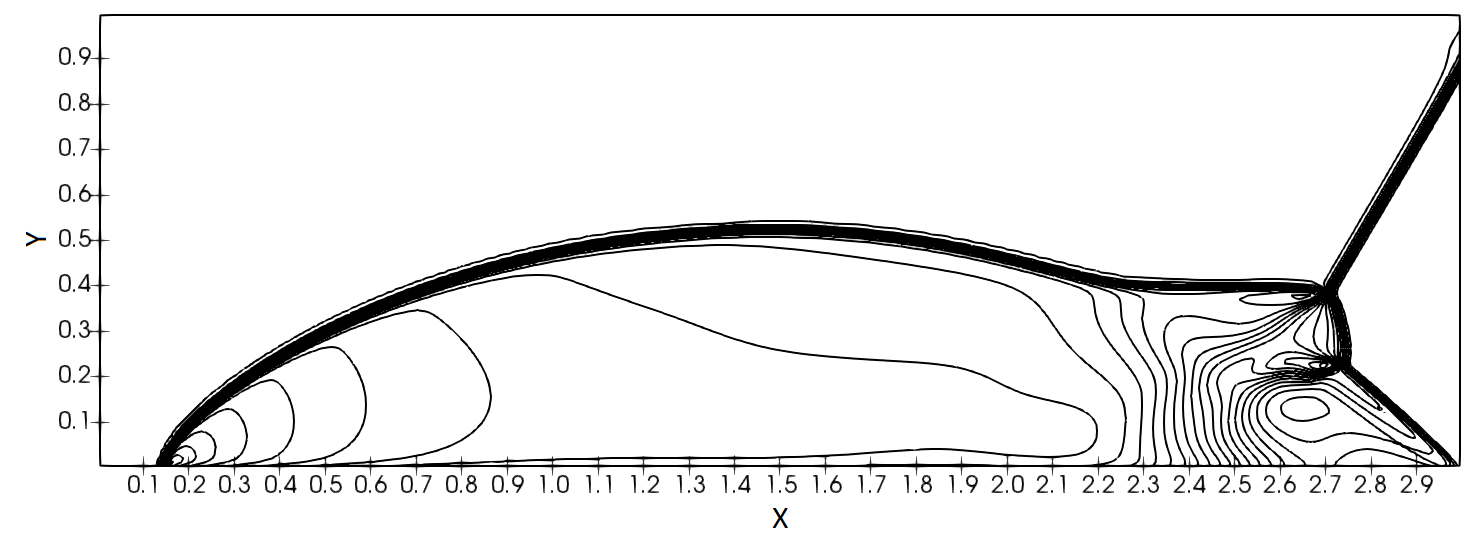} \includegraphics[width=220pt]{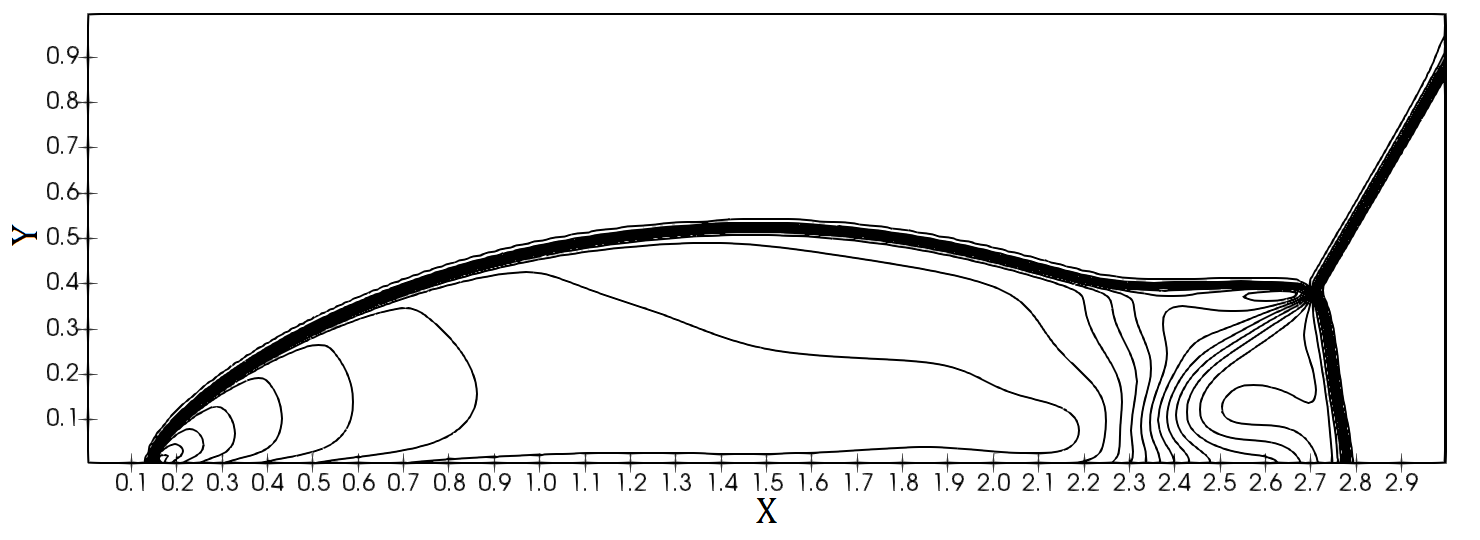}\\ 
		(c) HLLEM  \hspace{5cm} (d) HLLEM-FP  \\
		\caption{Density contours for the $M_{\infty}=10$ double Mach reflection problem computed by the HLL-CPS, HLL-CPS-FP, HLLEM and HLLEM-FP schemes. The results are shown at time $t=0.020026$ units.}
		\label{dmr-hllcps-hllem}
	\end{center}
\end{figure}

The domain is four units long and one unit wide. The domain is divided into $480\times120$ cells. A Mach 10 shock initially making a 60-degree angle at $x=1/6$ with the bottom reflective wall is made to propagate through the domain. The domain ahead of shock is initialized to pre-shock value ($\rho=1.4$, $p=1$, $u=0$, $v=0)$ and domain behind shock is assigned post-shock values. The inlet boundary condition is set to post-shock values and the outlet boundary is set to zero gradients. The top boundary is set to simulate actual shock movement. At the bottom, the post-shock boundary condition is set up to x=1/6 and the reflective wall boundary condition is set thereafter. Density contours are shown in Fig. \ref{dmr-hllcps-hllem} for the original and proposed HLL-CPS and HLLEM schemes at time $t=2.00260\times{}10^{-1}$ units. Thirty (30) contour levels ranging from 2.0 to 21.5 are shown. It can be seen from the figure that the proposed HLL-CPS-FP scheme is free from kinked Mach stem, just like the original HLL-CPS \cite{mandal}. A severe kinked Mach stem is observed in the original HLLEM scheme, while the proposed HLLEM-FP scheme is free from the kinked Mach stem problem.
\subsection{Forward-facing step problem}
\begin{figure}[H]
	\begin{center}
		\includegraphics[width=220pt]{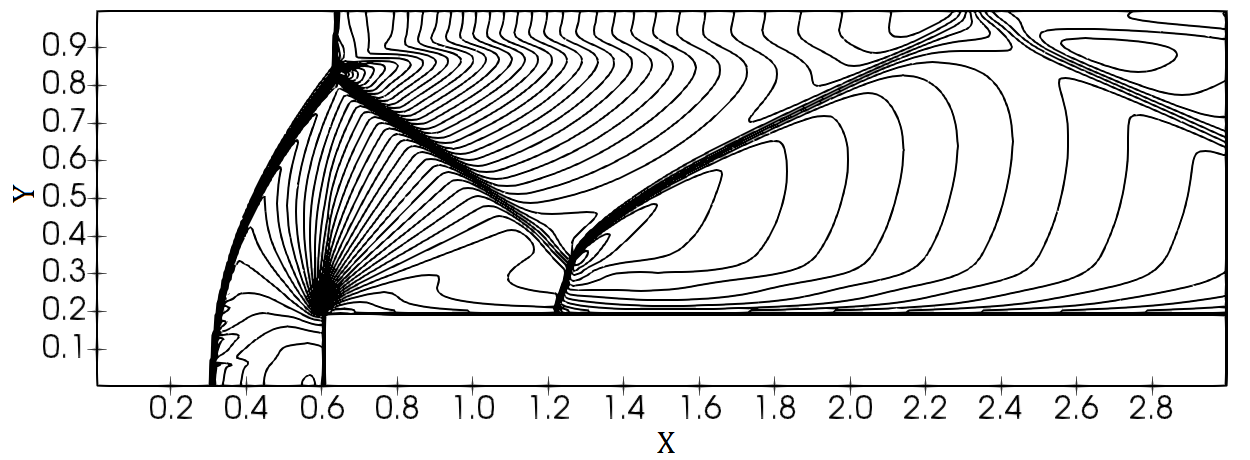} \includegraphics[width=220pt]{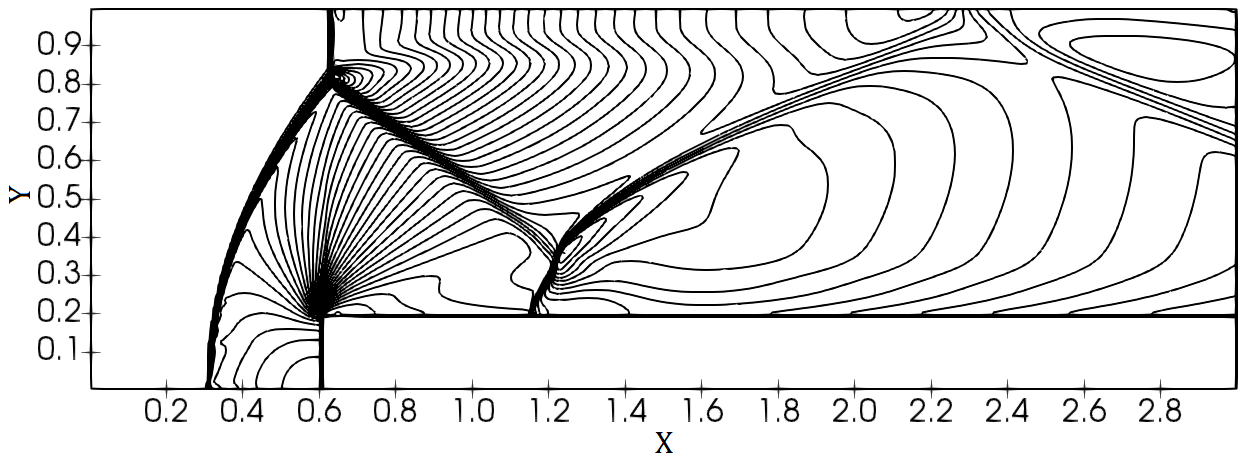}\\
		(a) HLL-CPS \hspace{5cm} (b) HLL-CPS-FP\\
		\includegraphics[width=220pt]{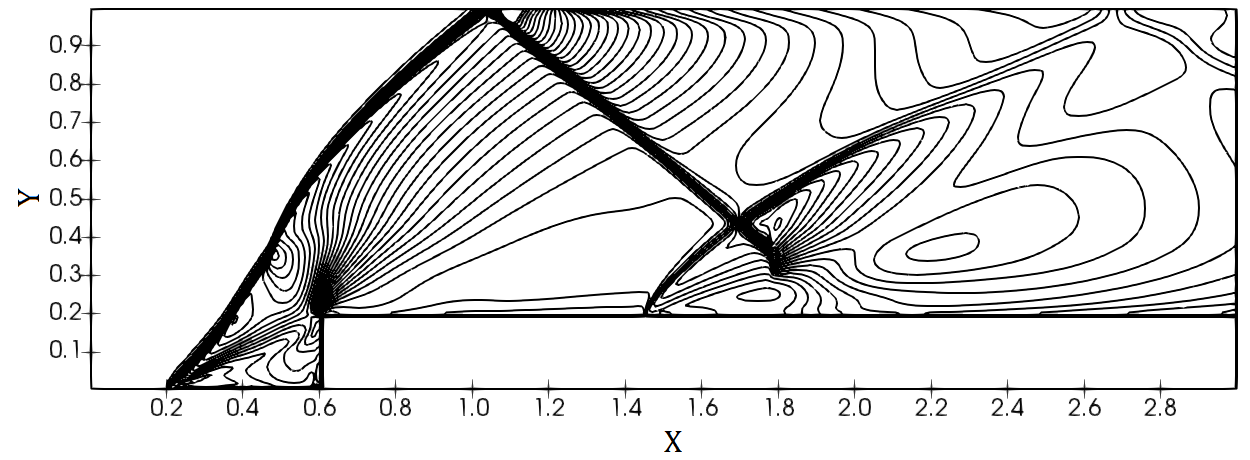} \includegraphics[width=220pt]{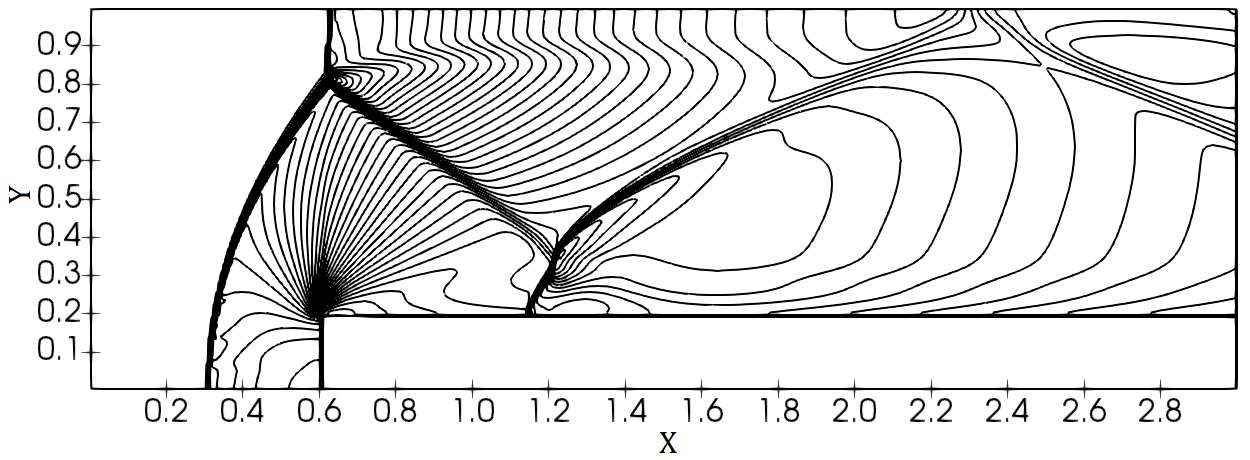} \\
		(c) HLLEM  \hspace{5cm} (d)  HLLEM-FP 
		\caption{Density contour for the $M_{\infty}=3$ flow over a forward-facing step computed by the first-order HLL-CPS, HLL-CPS-FP, HLLEM and HLLEM-FP schemes. The results are shown at time $t=4$ units.}
		\label{ffs-hllcps-hllem-o1}
	\end{center}
\end{figure}

\begin{figure}[H]
	\begin{center}
		\includegraphics[width=220pt]{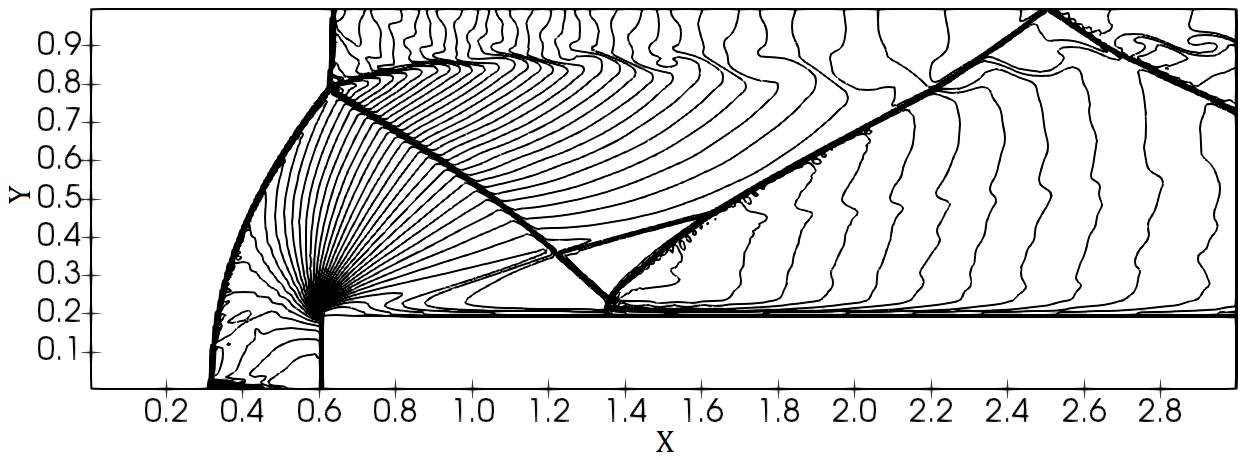} \includegraphics[width=220pt]{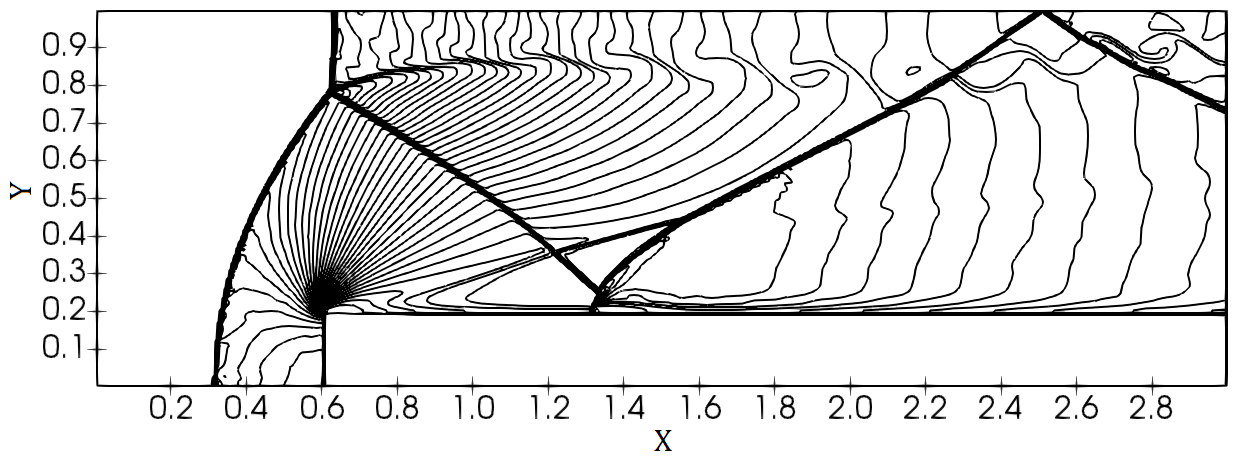}\\
		(a) HLL-CPS \hspace{5cm} (b) HLL-CPS-FP\\
		\includegraphics[width=220pt]{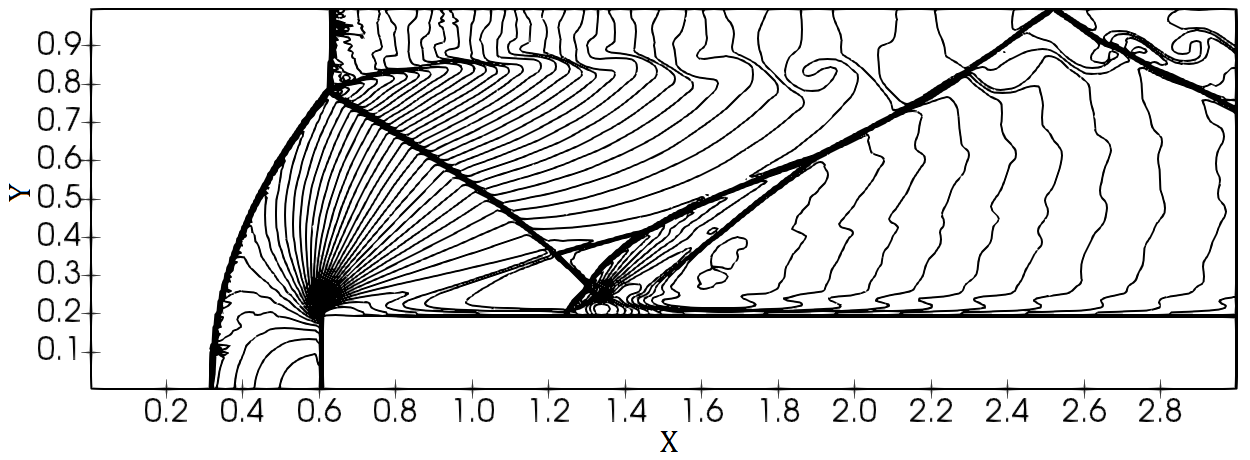} \includegraphics[width=220pt]{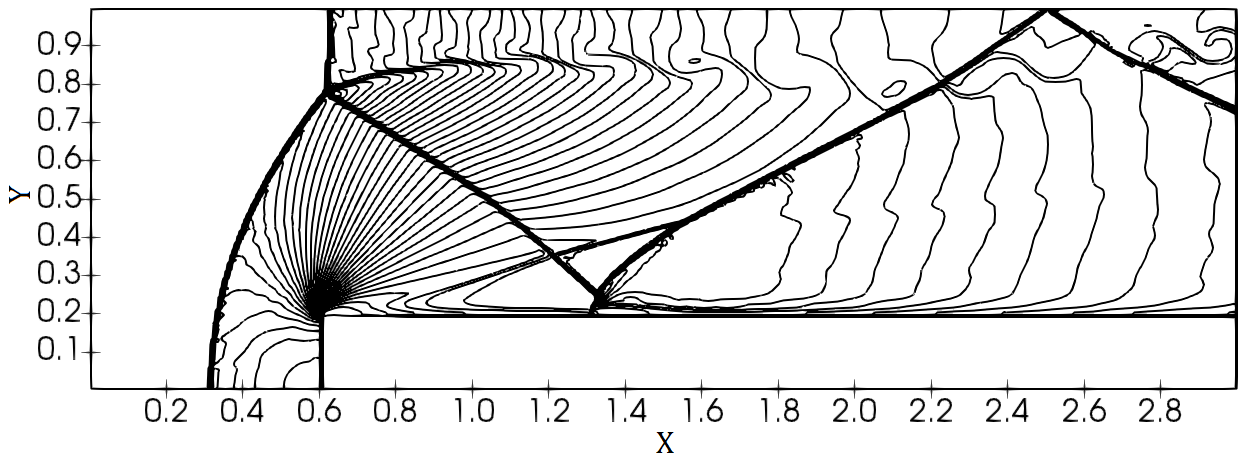} \\
		(c) HLLEM  \hspace{5cm} (d)  HLLEM-FP 
		\caption{Density contour for the $M_{\infty}=3$ flow over a forward-facing step computed by the second-order HLL-CPS, HLL-CPS-FP, HLLEM and HLLEM-FP schemes. The results are shown at time $t=4$ units.}
		\label{ffs-hllcps-hllem-o2}
	\end{center}
\end{figure}

The problem consists of a Mach 3 flow over a forward-facing step. The geometry consists of a step located 0.6 units downstream of the inlet and 0.2 units high. A mesh of $480\times160$ cells is used for a domain of 3 units long and 1 unit high. The complete domain is initialized with $\rho=1.4$, $p=1$, $u=3$, and $v=0$. The inlet boundary is set to free-stream conditions, while the outlet boundary has zero gradients. At the top and bottom, reflective wall boundary conditions are set. Density contour plots for the original and proposed first-order HLL-CPS and HLLEM schemes are shown in Fig. \ref{ffs-hllcps-hllem-o1} at time $t=4.0$ units. A total of 45 contours from 0.2 to 7.0 are shown in the figure. It can be seen from the figure that both the proposed HLL-CPS schemes capture the primary and reflected shocks without any numerical instability while a mild instability is observed in the primary bow ahead of the step in the HLL-CPS scheme. The primary bow shock is replaced by an oblique shock in the original HLLEM scheme. The HLLEM scheme fails to capture the triple point and the primary reflected shock is also not well captured. The proposed HLLEM-FP scheme captures the primary and reflected shocks without numerical instabilities.
It is felt that the first-order schemes are unable to resolve all the flow features of the problem and hence, second-order computations are also carried out. The popular MUSCL reconstruction method \cite{van-leer-muscl} along with the van Leer limiter is used for achieving second-order accuracy in space and time integration is performed using the two-stage Runge-Kutta method of Gottlieb and Shu \cite{gott}. Density contour plots for the original and proposed second-order schemes are shown in Fig. \ref{ffs-hllcps-hllem-o2} at time $t=4.0$ units, and a total of 45 contours from 0.2 to 7.0 are shown in the figure. It can be seen from the figure that the original and proposed schemes are able to resolve the contact discontinuity originating from the triple point near the top wall. A mild instability is observed in the primary bow shock ahead of the step in the HLL-CPS scheme. The HLLEM scheme exhibit mild instability in the primary bow shock ahead of the step and in the normal shock near the top wall. The HLLEM scheme is also unable to resolve the secondary reflected shock around the step region. The proposed HLL-CPS-FP and HLLEM-FP schemes, on the other hand, do not exhibit any numerical shock instability and are also able to resolve the contact discontinuity originating at the triple point near the top wall. It is thus demonstrated that the proposed HLL-CPS-FP and HLLEM-FP schemes are free from numerical shock instabilities and are capable of resolving contact discontinuities at the same time.

\subsection{Blunt body problem}
\begin{figure}[H]
	\begin{center}
	\includegraphics[width=100pt]{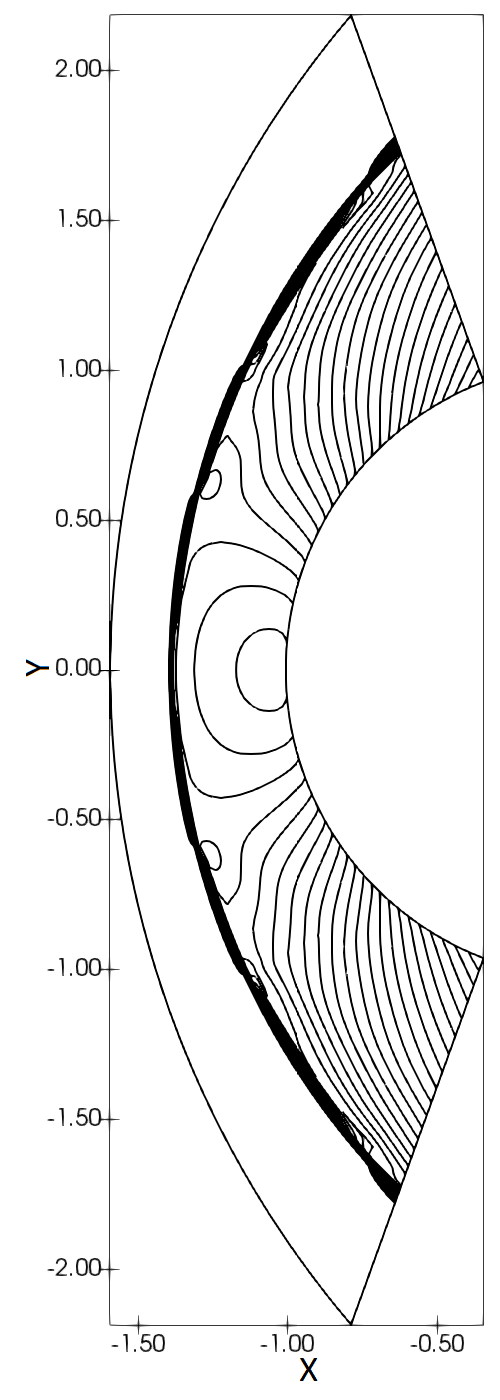}  \includegraphics[width=100pt]{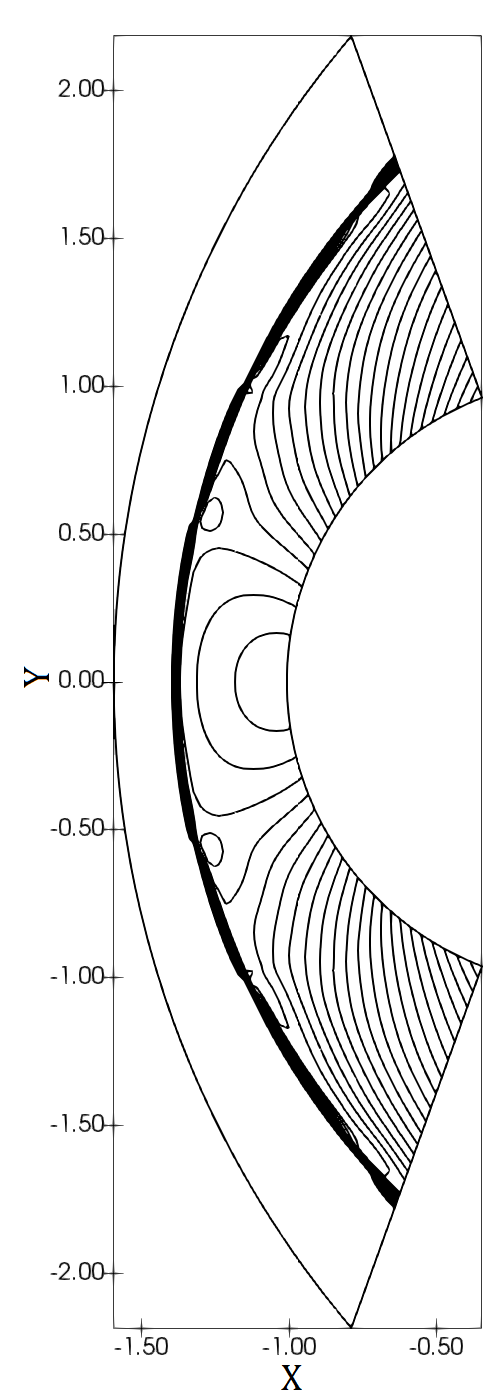}
	\includegraphics[width=100pt]{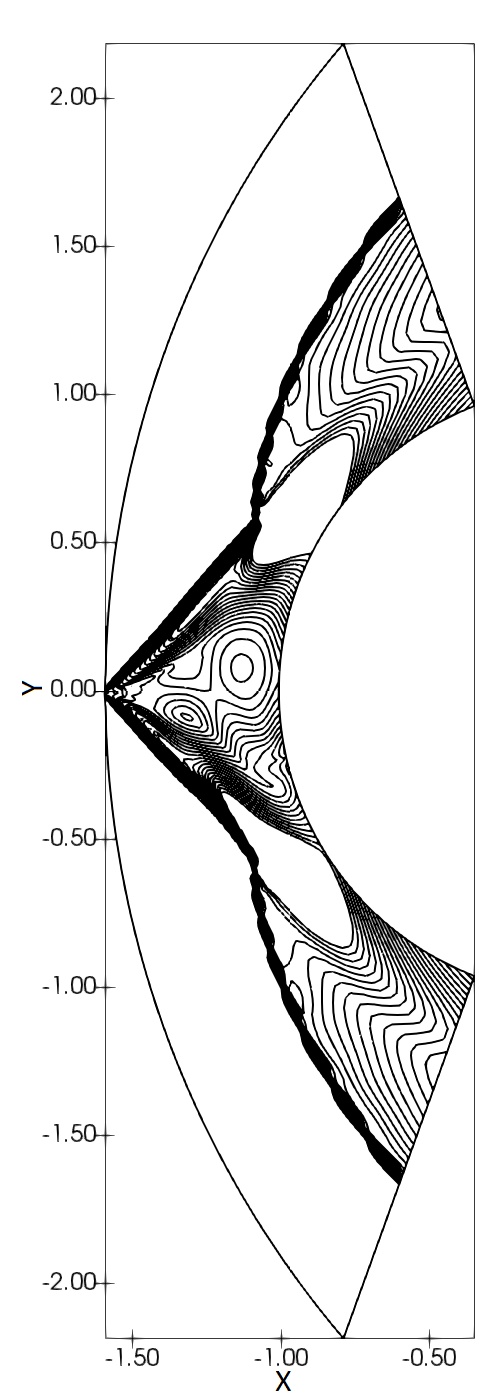} \includegraphics[width=100pt]{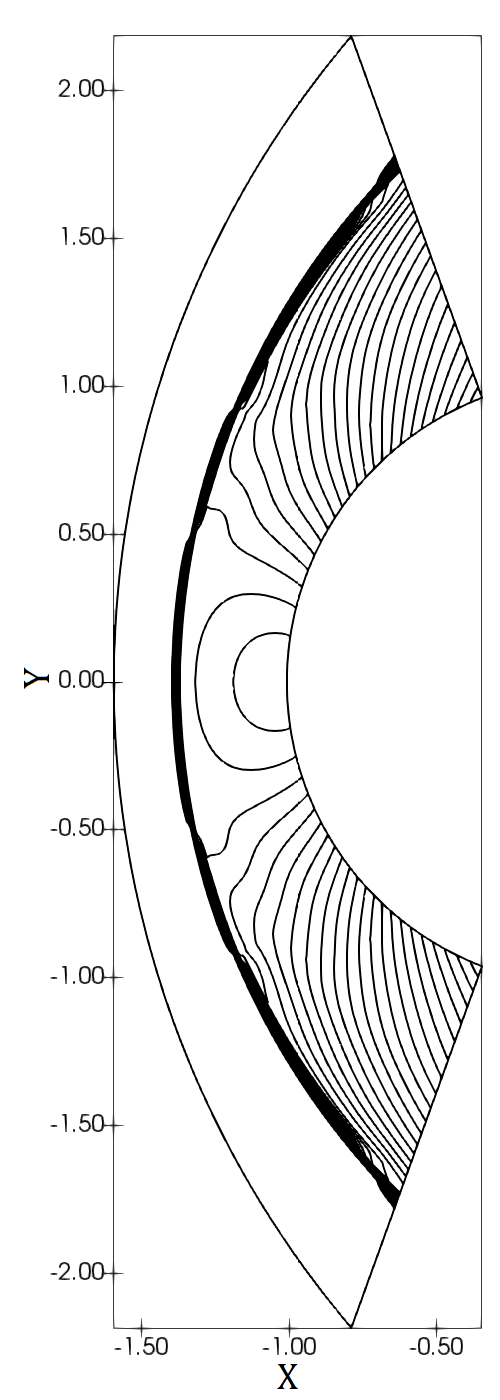} \\
	(a) HLL-CPS  \hspace{1cm} (b)  HLL-CPS-FP  \hspace{1cm}  (c) HLLEM  \hspace{1cm} (d) HLLEM-FP  
	\caption{Density contour for $M_{\infty}=20$ flow over a blunt body computed by the HLL-CPS, HLL-CPS-FP, HLLEM and HLLEM-FP schemes. The results are shown after 100,000 iterations}
		\label{carbuncle-hll-xfp}
	\end{center}
\end{figure}

\begin{figure}[H]
	\begin{center}
		\includegraphics[width=250pt]{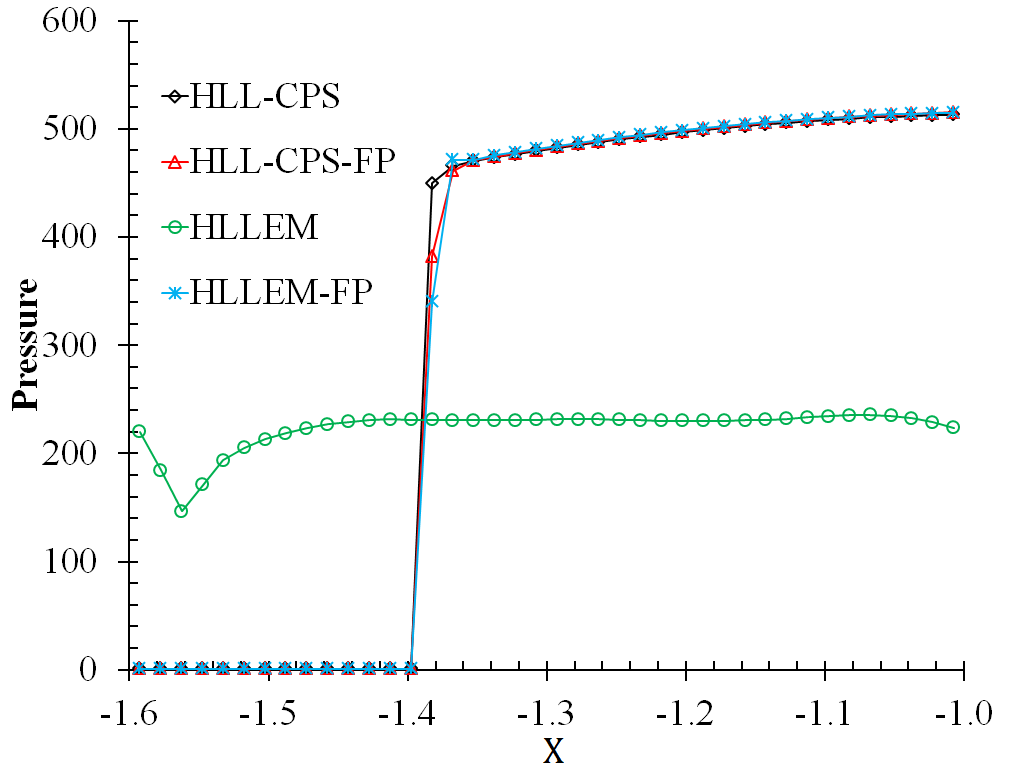} 
		\caption{Static pressure plot along the centreline $(y=0)$ for $M_{\infty}=20$ flow over a blunt body}
		\label{cp-carbuncle}
	\end{center}
\end{figure}

\begin{figure}[H]
	\begin{center}
	\includegraphics[width=225pt]{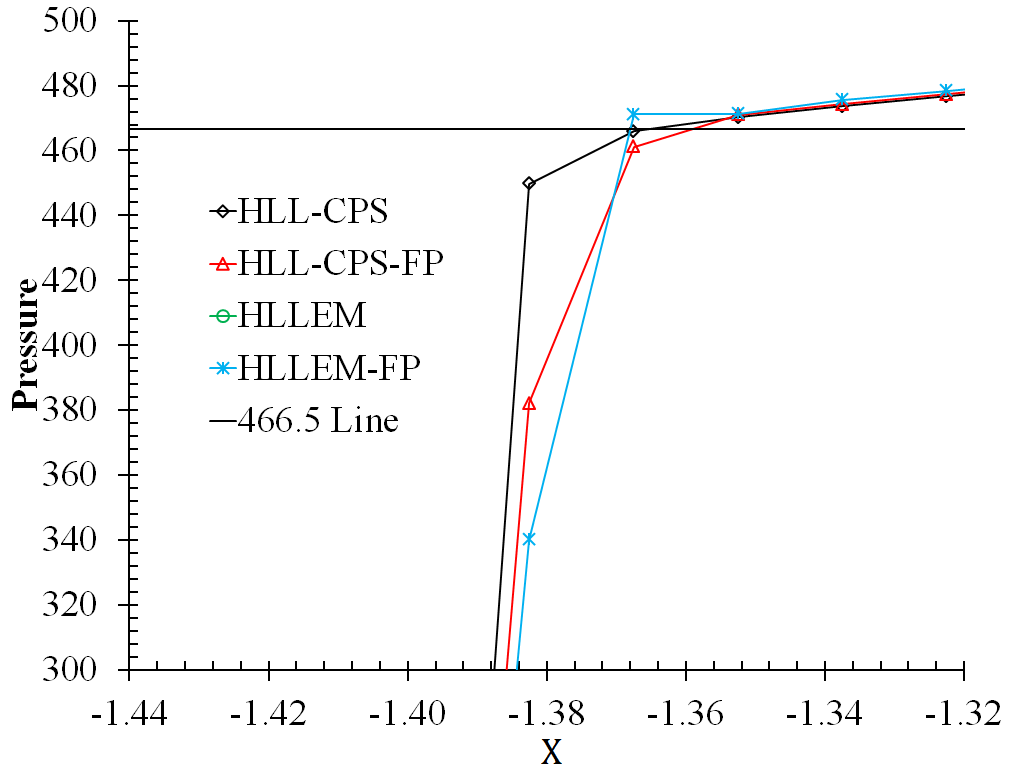} \includegraphics[width=225pt]{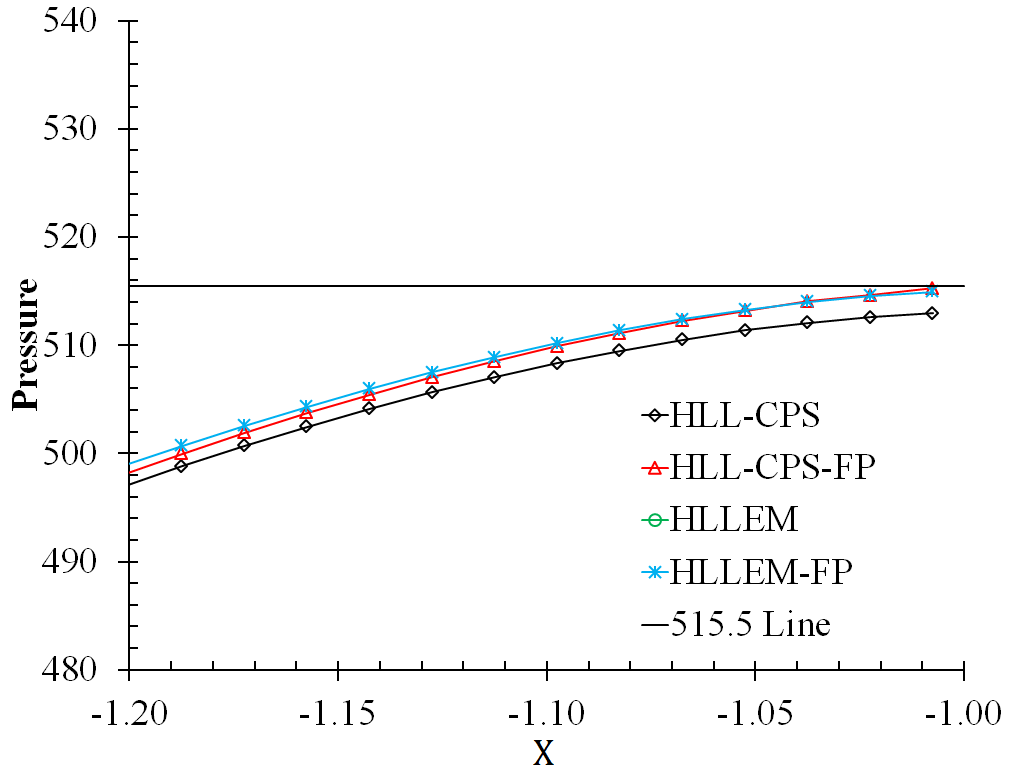}  \\
	(a) Near Shock Location \hspace{1cm} (b) Ahead of the  Blunt Body
	\caption{Zoomed-up view of static pressure plot along the centreline for $M_{\infty}=20$ flow over a blunt body}
	\label{cp-zoom-carbuncle}
	\end{center}
\end{figure}

\begin{figure}[H]
	\begin{center}
		\includegraphics[width=220pt]{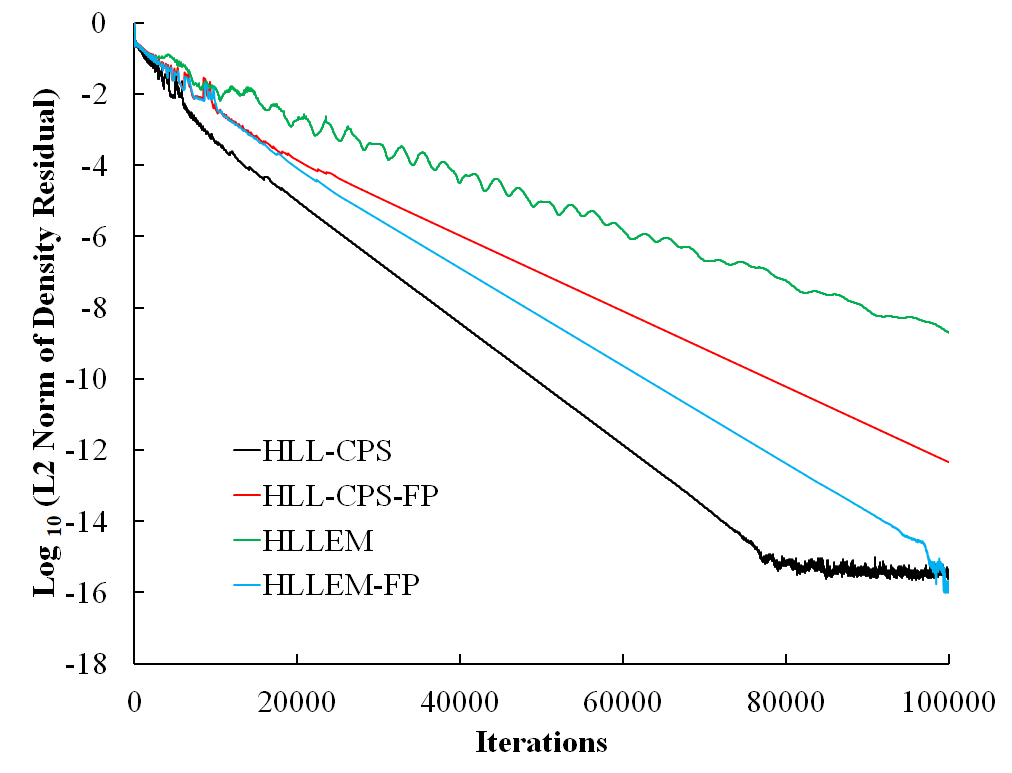} 
		\caption{Convergence history of L2 norm of density residuals for $M_{\infty}=20$ flow over a blunt body computed by the HLL-CPS, HLLEM, HLL-CPS-FP and HLLEM-FP schemes}
		\label{convergence-carbuncle-hllx-fp}
	\end{center}
\end{figure}

Hypersonic flow over a blunt body is a typical problem to assess the performance of a scheme for carbuncle instability. A free-stream Mach number of 20 is used in the computation. The grid for the blunt body had a size of $40\times320$ cells. The domain is initialized with values of $\rho=1.4$, $p=1$, $u=20$, and $v=0$. The inlet boundary condition is set to free-stream value. The solid wall boundary condition is applied to the blunt body. The computations are carried out for 100,000 iterations. The density contour plot for the original and the modified HLL-CPS and HLLEM schemes are shown in Fig. \ref{carbuncle-hll-xfp}, and a total of 27 density contours from 2.0 to 8.7 are drawn. It can be seen from the figure that the original and the modified HLL-CPS schemes are free from the carbuncle phenomenon. It can be seen from the figure that the original HLLEM scheme exhibits a severe carbuncle phenomenon, while the proposed HLLEM-FP scheme is free from the carbuncle problem of the original scheme. The static pressure plot of various schemes along the centreline of the blunt body is shown in Fig. \ref{cp-carbuncle}. It can be seen from the figure that due to the severe carbuncle phenomenon, the static pressure profile of the HLLEM scheme is drastically different from the other schemes which are free from the carbuncle phenomenon. The zoomed-up view of the centreline static pressure plot near the shock location and ahead of the blunt body is shown in Fig. \ref{cp-zoom-carbuncle}. It can be seen from the figure that the post-shock static pressure value of the HLL-CPS, HLL-CPS-FP and HLLEM-FP are very close to the analytical value of 466.5 Pa. The stagnation pressure of the HLL-CPS-FP and HLLEM-FP are about 515 Pa and are very close to the analytical value of 515.5 Pa while the stagnation pressure of the HLLE and HLL-CPS schemes are about 512 Pa. Thus, the stagnation pressure of all the new schemes are closer to the analytical value. Therefore, it is felt that the low Mach corrections to the HLL-CPS and HLLEM schemes may be reducing the pressure losses in the subsonic region downstream of the shock and thus contributing to the improved stagnation pressure values in the new schemes. The post-shock pressure and stagnation pressure of the various schemes are presented in Table \ref{shock-hllx-fp} and the improved stagnation pressure of the proposed schemes can be observed in the table. It is observed that the stagnation pressure of the HLLC-ADC, HLLC-SWM, HLLEM-ADC, and HLLEM-SWM schemes of Simon and Mandal \cite{san1, san2, san3} which do not involve low Mach corrections are around 509 Pa which lends credibility to the above observation. The convergence history of the original and modified HLL-CPS and HLLEM schemes are shown in Fig. \ref{convergence-carbuncle-hllx-fp}. It can be seen from the figure that the original HLLEM scheme converges to about $10^{-8}$ while the modified HLLEM-FP scheme converges to machine accuracy within 100,000 iterations. The original HLL-CPS scheme also converges to about $10^{-16}$ while the HLL-CPS-FP scheme converges to about $10^{-12}$. The convergence histories of the proposed HLLEM-FP and HLL-CPS-FP schemes are better than the HLLEM-ADC scheme of Simon and Mandal \cite{san3} and can be considered to be satisfactory.
\begin{table}[H]
\begin{center}
\begin{doublespace}
\caption{Post Shock and Stagnation Pressure values of various schemes}
\label{shock-hllx-fp}
\begin{tabular}{|c|c|c|c|}
\hline Sl No & Scheme & Post Shock Pressure (Pa) & Stagnation Pressure (Pa) \\
\hline 1 & Analytical &466.5 & 515.5 \\
\hline 2 & HLLE &465.84 & 512.23 \\
\hline 3 & HLL-CPS & 465.99 & 512.94\\
\hline 4 &HLL-CPS-FP &  461.13& 515.25\\
\hline 5 & HLLEM & 228.24 & 220.61 \\
\hline 6 &HLLEM-FP & 471.22 & 514.95 \\
\hline
\end{tabular}
\end{doublespace}
\end{center}
\end{table}
\subsection{Supersonic corner problem}
\begin{figure}[H]
	\begin{center}
		\includegraphics[width=175pt]{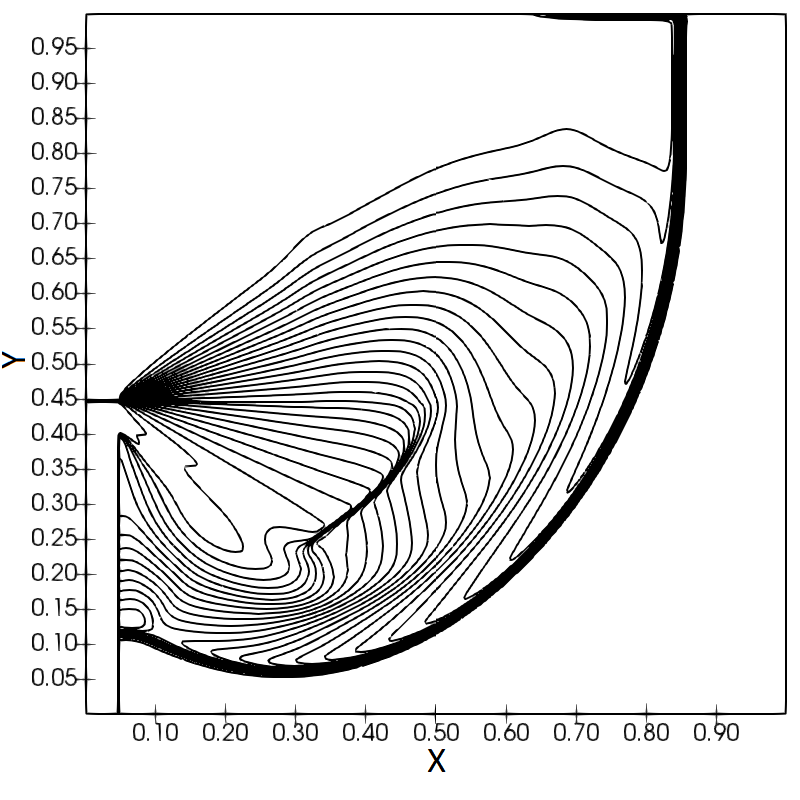}  \includegraphics[width=175pt]{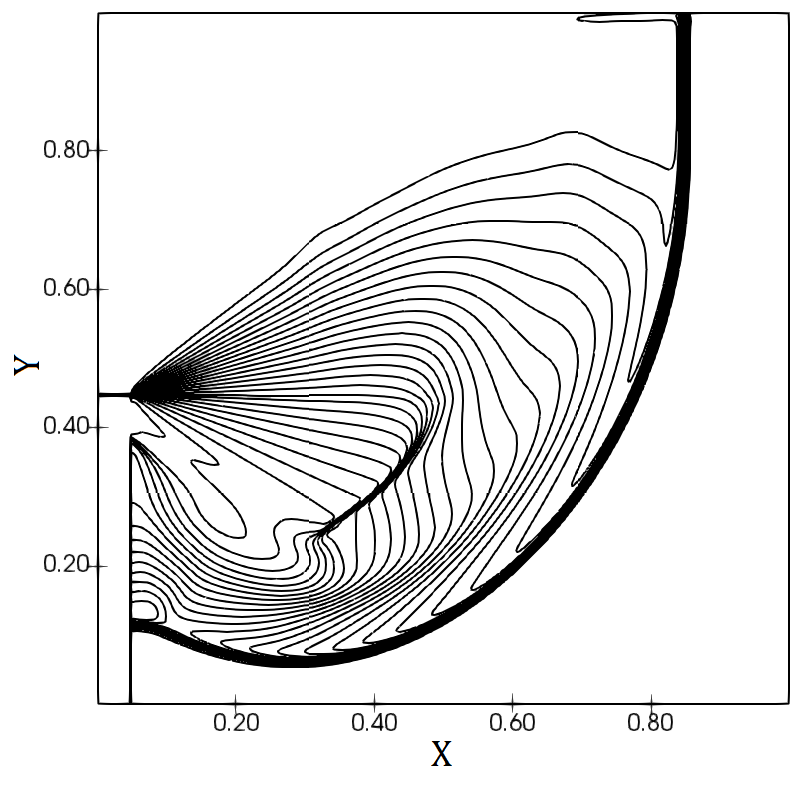} \\
		(a) HLL-CPS Scheme \hspace{3cm} (b) HLL-CPS-FP Scheme \\ 
		\includegraphics[width=175pt]{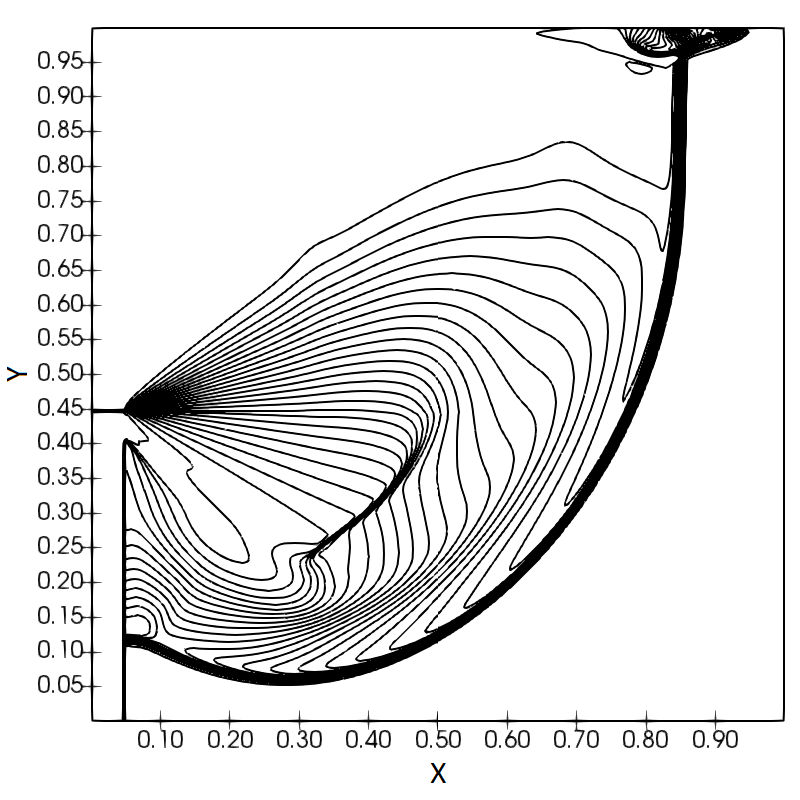}  \includegraphics[width=175pt]{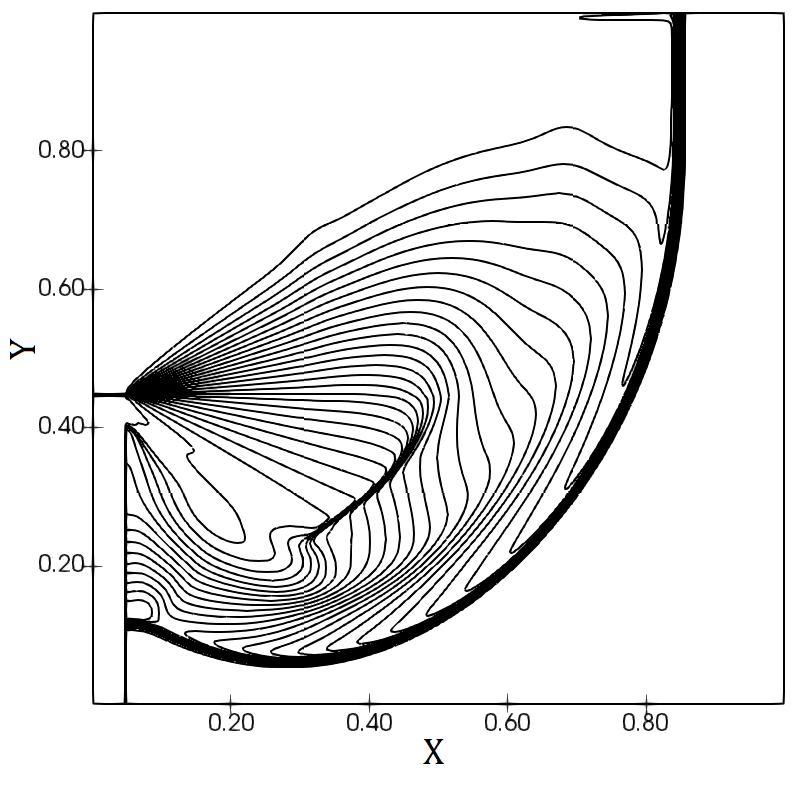} \\ 
		(c) HLLEM Scheme \hspace{3cm} (d)  HLLEM-FP Scheme \\ 
		\caption{Density contour for $M_{\infty}=5.09$ normal shock diffraction around a $90^0$ corner computed by the HLL-CPS, HLL-CPS-FP, HLLEM and HLLEM-FP schemes. The results are shown for time t=0.1561.}
		\label{scr-hll-xfp}
	\end{center}
\end{figure}

The problem consists of a sudden expansion of a Mach 5.09 normal shock around a 90-degree corner. The domain is a square of one unit and is divided into $400\times400$ cells. The corner is located at $x=0.05$ and $y=0.45$. The initial normal shock is located at $x=0.05$. The domain to the right of shock is assigned initially with pre-shock conditions of $\rho{}=1.4,\;p=1.0,\;u=0.0,$ and $v=0.0$. The domain to the left of the shock is assigned post-shock conditions. The inlet boundary is supersonic, the outlet boundary has zero gradients and the bottom boundary behind the corner uses extrapolated values. Reflective wall boundary conditions are imposed on the corner. A CFL value of 0.8 is used and the density contour plots are generated for time $t=0.1561$ units. A total of 30 density contours ranging from 0 to 7.1 is shown in Fig. \ref{scr-hll-xfp} for the various HLL-CPS and HLLEM schemes. The figure shows a benign perturbation at the top just ahead of the shock for HLL-CPS schemes. The results of the HLL-CPS-FP scheme show significant improvement over the original HLL-CPS scheme. It can be seen from the figure that severe numerical shock instability is present in the original HLLEM scheme, while the figure that the proposed HLLEM-FP scheme is free from the numerical shock instability problem of the original scheme. 
\subsection{Flow around a cylinder at low Mach numbers}
\begin{figure}[H]
	\begin{center}
		\includegraphics[width=150pt]{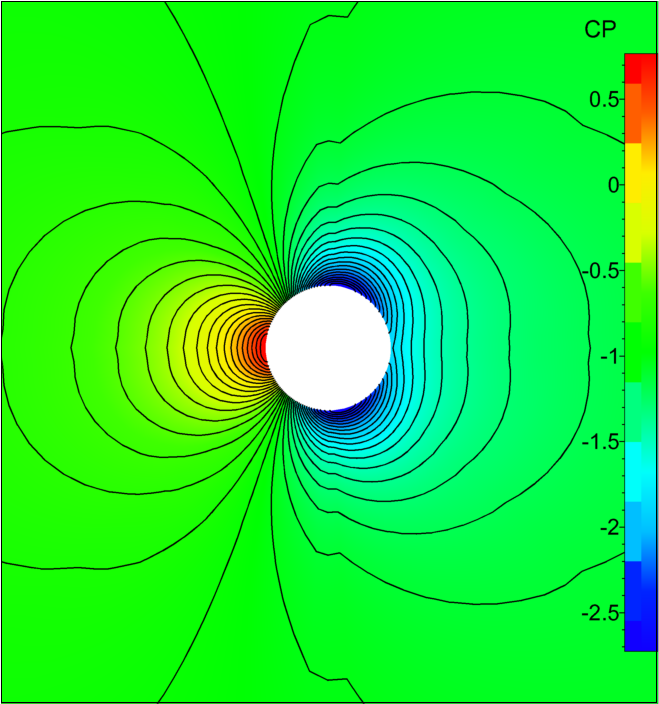} \includegraphics[width=150pt]{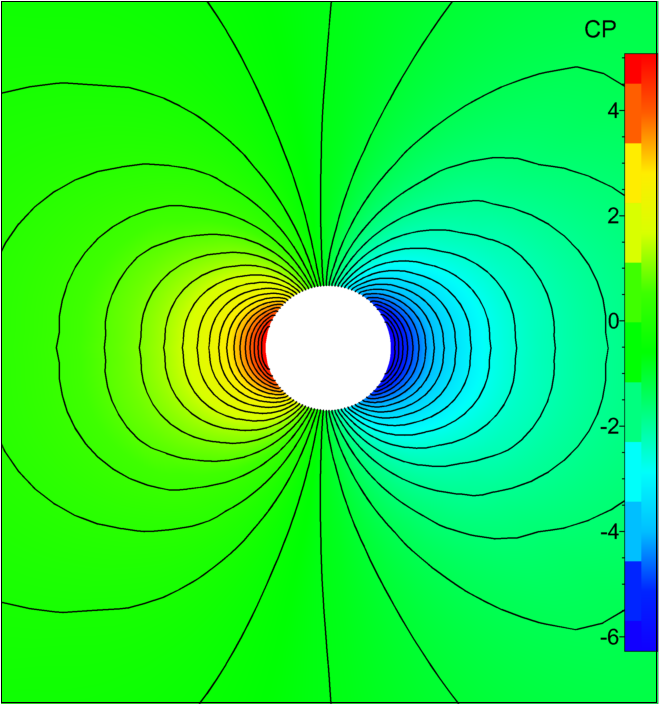} \includegraphics[width=150pt]{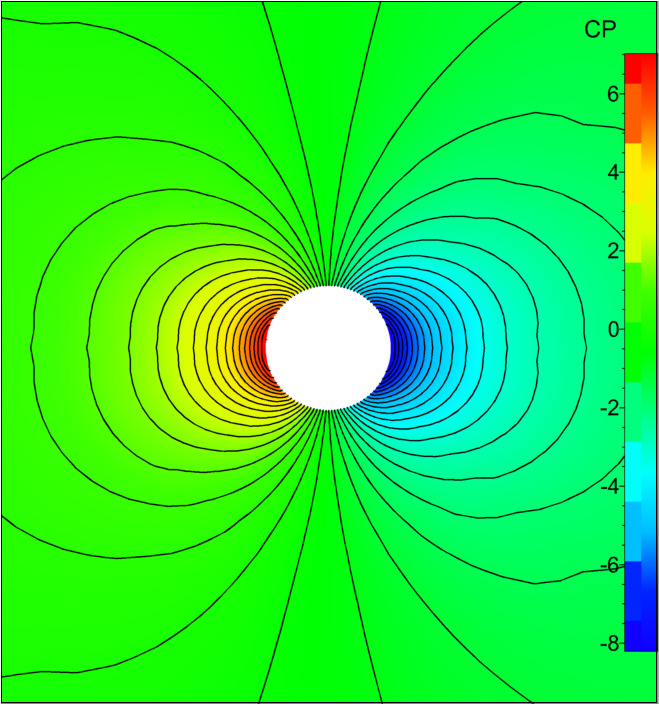}\\
  		(a) $M_{\infty}=0.10$ \hspace{3cm}  (b) $M_{\infty}=0.01$ \hspace{3cm} 	(c) $M_{\infty}=0.001$
\caption{Static pressure coefficient contour plot for flow around a cylinder computed by the HLL-CPS scheme for (a) $M_{\infty}=0.10$, (b) $M_{\infty}=0.01$ and (c) $M_{\infty}=0.001$}
\label{cyl-hllcps}
	\end{center}
\end{figure}

\begin{figure}[H]
	\begin{center}
		\includegraphics[width=150pt]{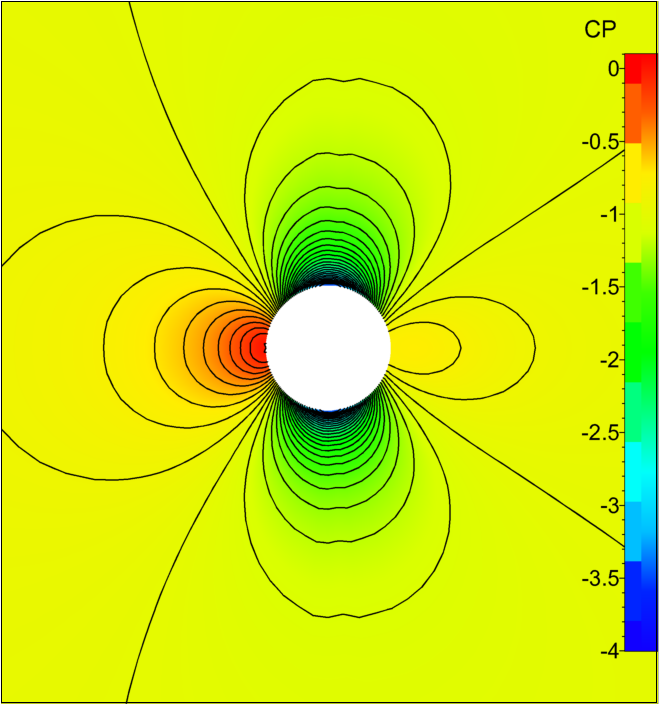} \includegraphics[width=150pt]{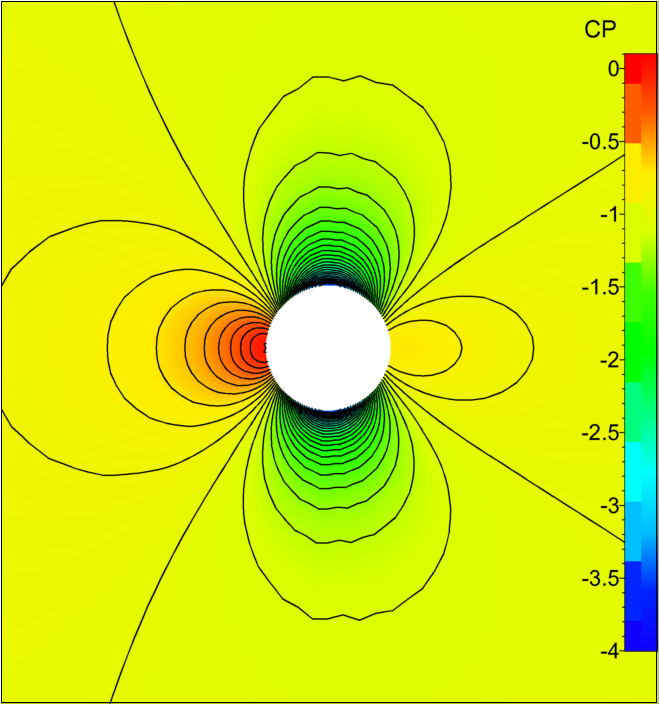} \includegraphics[width=150pt]{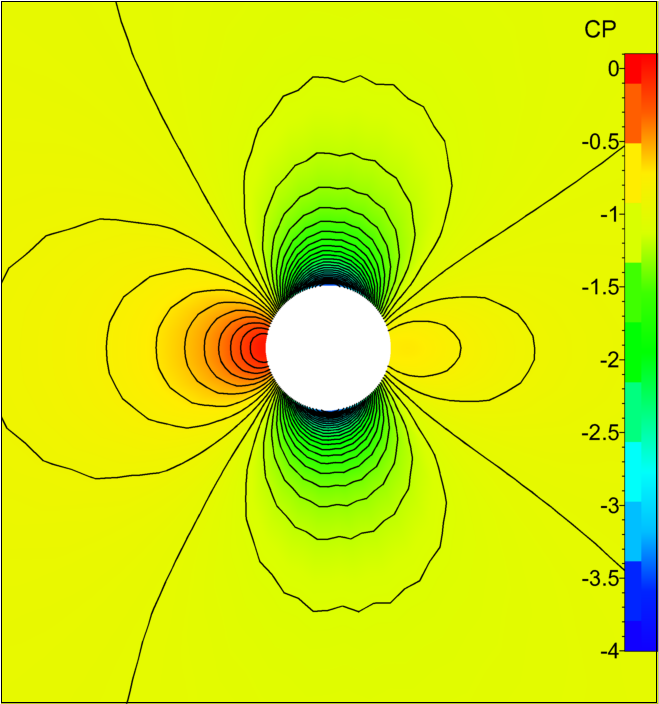}\\
  		(a) $M_{\infty}=0.10$ \hspace{3cm}  (b) $M_{\infty}=0.01$ \hspace{3cm} (c) $M_{\infty}=0.001$
\caption{Static pressure coefficient contour plot for flow around a cylinder computed by the HLL-CPS-FP scheme for (a) $M_{\infty}=0.10$, (b) $M_{\infty}=0.01$ and (c) $M_{\infty}=0.001$}
\label{cyl-hllcps-fp}
	\end{center}
\end{figure}

\begin{figure}[H]
	\begin{center}	
		\includegraphics[width=150pt]{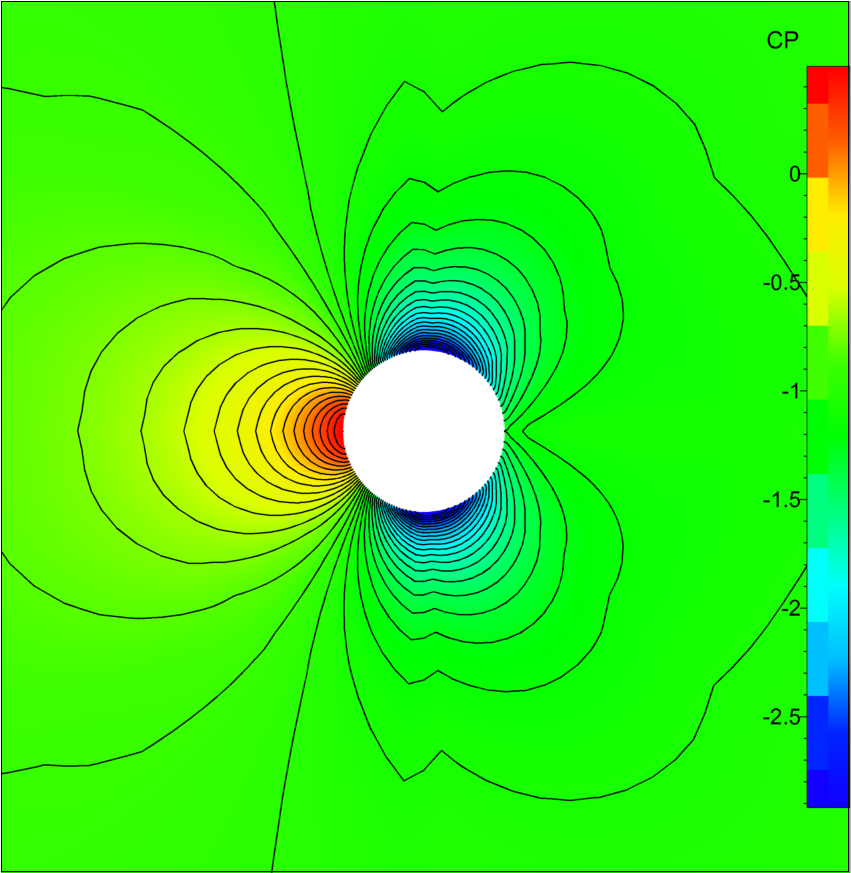} 
		\includegraphics[width=150pt]{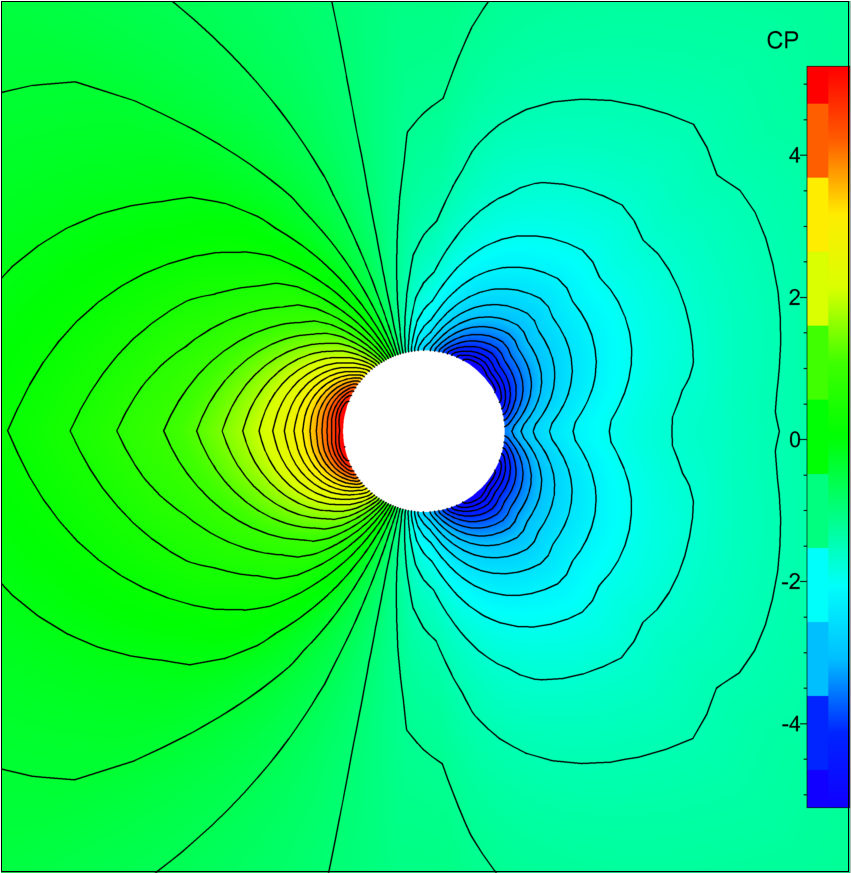} 
		\includegraphics[width=150pt]{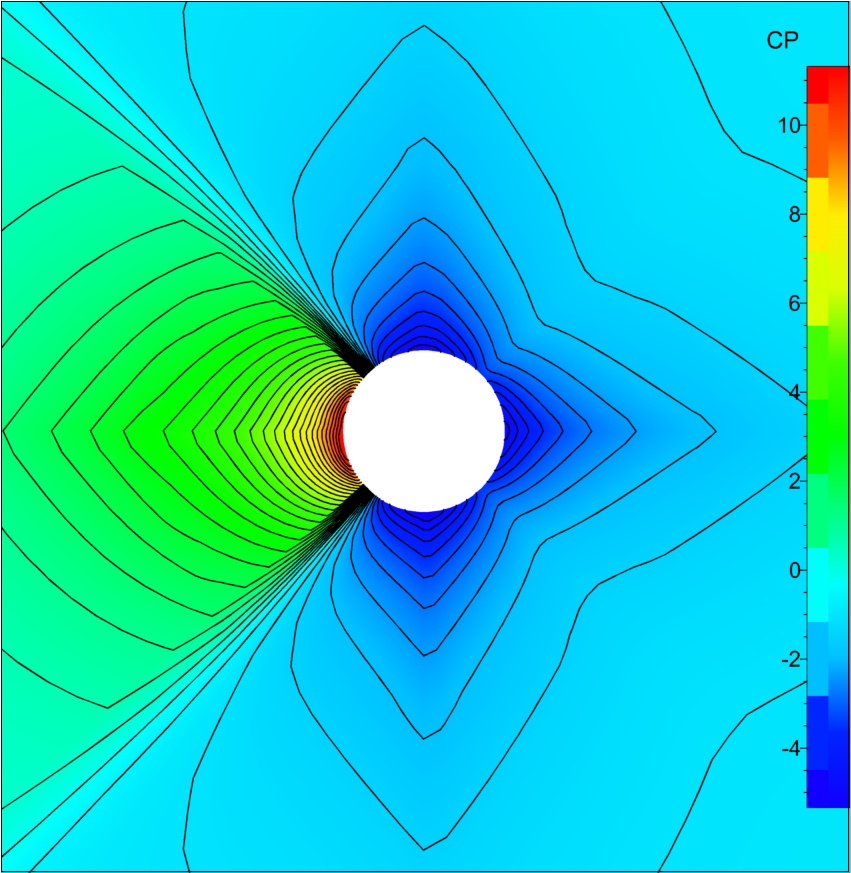}\\
		(a) $M_{\infty}=0.10$ \hspace{3cm} (b) $M_{\infty}=0.01$\hspace{3cm} (c) $M_{\infty}=0.001$
		\caption{Static pressure coefficient contour plot for flow around a cylinder computed by the HLLEM scheme for (a) $M_{\infty}=0.10$, (b) $M_{\infty}=0.01$ and (c) $M_{\infty}=0.001$ }
		\label{cyl-hllem}
	\end{center}
\end{figure}

\begin{figure}[H]
	\begin{center}	
		\includegraphics[width=150pt]{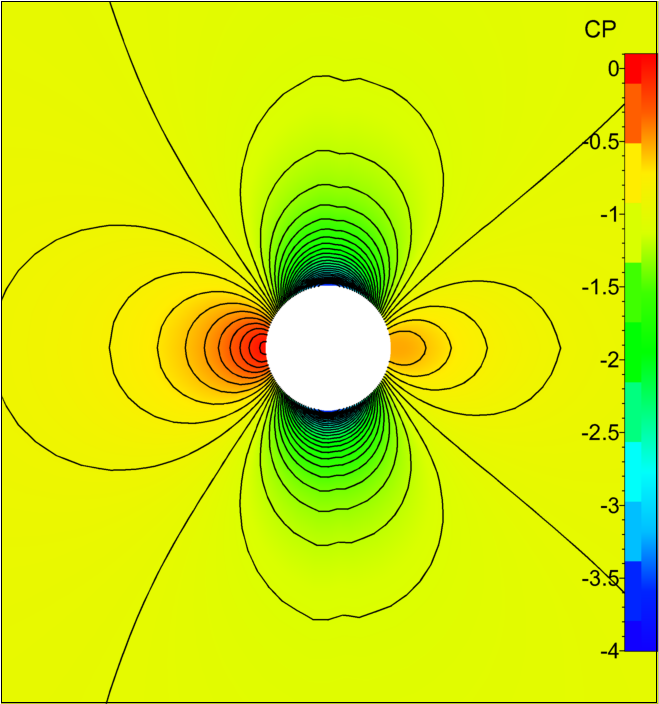} 
		\includegraphics[width=150pt]{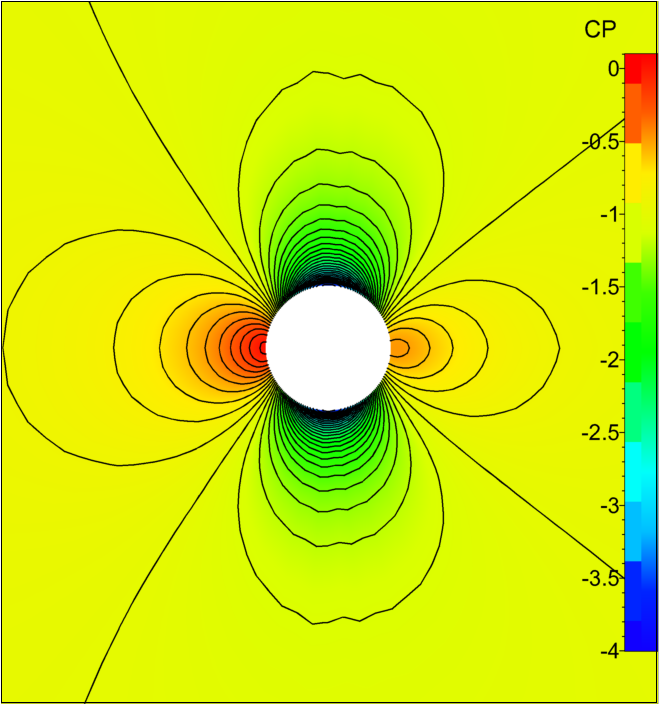} 
		\includegraphics[width=150pt]{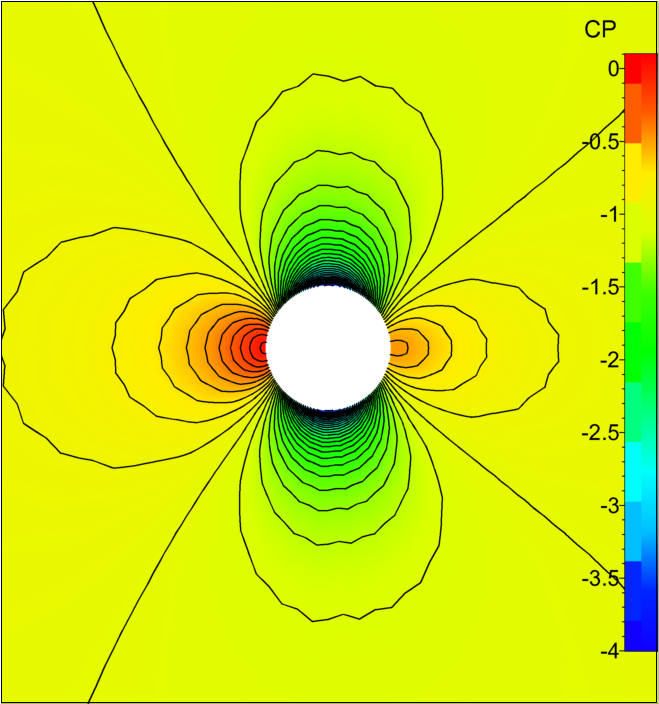}\\
		(a) $M_{\infty}=0.10$ \hspace{3cm} (b) $M_{\infty}=0.01$\hspace{3cm} (c) $M_{\infty}=0.001$
		\caption{Static pressure coefficient contour plot for flow around a cylinder computed by the proposed HLLEM-FP scheme for (a) $M_{\infty}=0.10$, (b) $M_{\infty}=0.01$ and (c) $M_{\infty}=0.001$}
		\label{cyl-hllem-fp}
	\end{center}
\end{figure}

\begin{figure}[H]
	\begin{center}	
		\includegraphics[width=220pt]{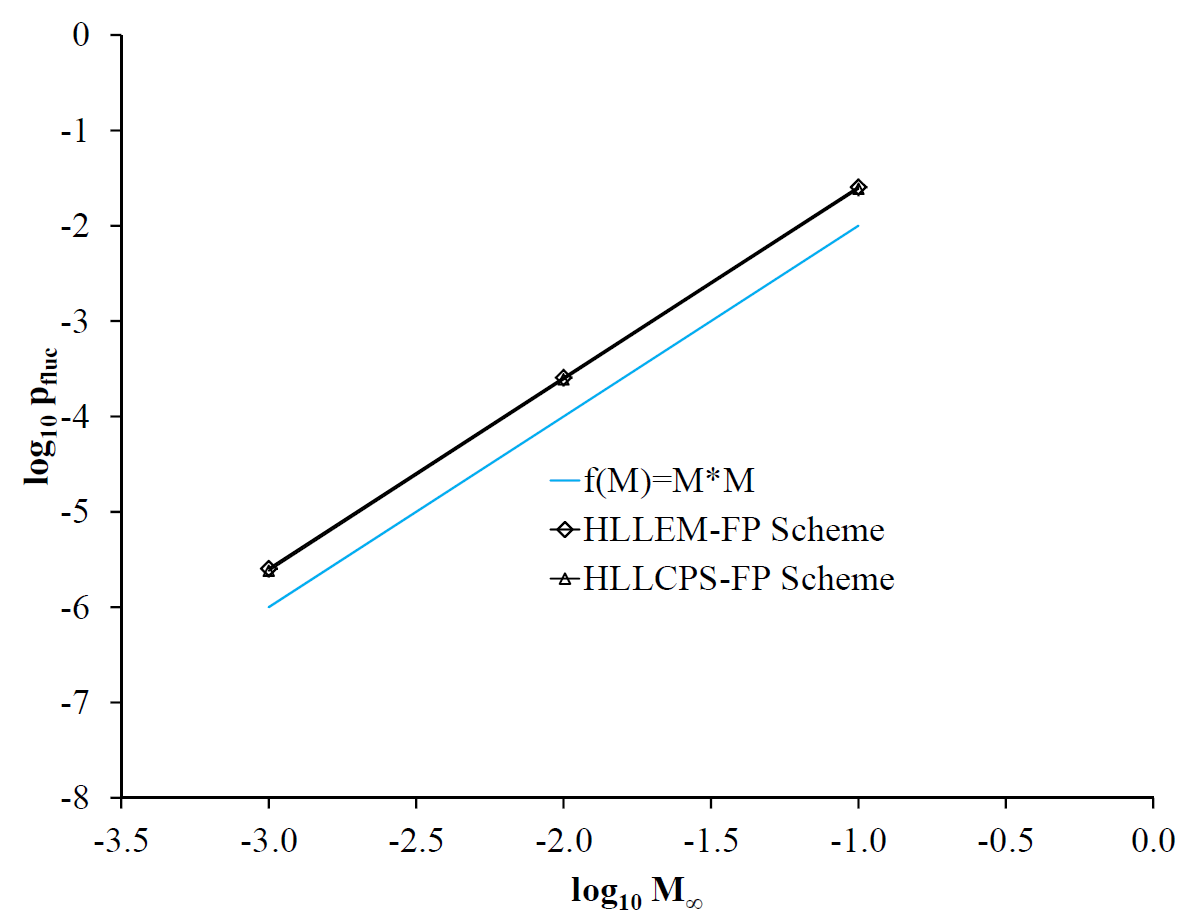} 
		\caption{Maximum pressure fluctuation $p_{fluc}=(p_{max}-p_{min})/p_{max}$ against inflow Mach number for the flow around a cylinder obtained with the proposed HLL-CPS-FP and HLLEM-FP schemes}
		\label{cylinder-pfluc-hll-xfp}
	\end{center}
\end{figure}

The flow around a cylinder at low Mach numbers is an important test case to determine the accuracy and numerical dissipation of a numerical scheme. First-order, inviscid computations are carried out over a cylinder using the original and proposed HLL-CPS and HLLEM schemes. The computations are carried out with a CFL number of 0.80. Structured grids of size $49\times37\times5$, $97\times37\times5$, and $97\times73\times5$ are used and results are shown for a grid size of $97\times73\times5$. The static pressure coefficient plots of the HLL-CPS scheme for Mach 0.10, 0.01 and 0.001 are shown in Fig. \ref{cyl-hllcps} and a total of 41 contour lines are drawn. The static pressure coefficient is defined as $c_p=(p-p_{\infty})/\frac{1}{2}\rho_{\infty}u_{\infty}^2$ where $p_{\infty}$, $\rho_{\infty}$ and $u_{\infty}$ are the free-stream pressure, density and velocity respectively. Due to the increase in numerical viscosity with reduction in Mach number in the HLL-CPS scheme, the pressure contours resemble creeping Stokes flow. The static pressure coefficient plot shows the highest numerical viscosity for the Mach 0.001 case, indicating that numerical viscosity increases with a decrease in Mach number for the HLL-CPS scheme. The static pressure coefficient plots for the proposed HLL-CPS-FP scheme for Mach 0.1, 0.01 and 0.001 are shown in Fig. \ref{cyl-hllcps-fp} and total of 41 contour lines from -4.0 to 0.0 are drawn. Static pressure coefficient plot of the proposed scheme resembles potential flow even at Mach 0.001. A weak checkerboard is seen for the modified HLL-CPS-FP schemes at Mach 0.001. The improved accuracy of the proposed HLL-CPS scheme can be attributed to the reduction in numerical dissipation through an increase in anti-diffusion coefficients $\delta_3$ and $\delta_n$.
The static pressure coefficient contour plots of the HLLEM and the HLLEM-FP scheme are shown in Fig. \ref{cyl-hllem} and \ref{cyl-hllem-fp} respectively, and a total of 41 contour lines are drawn. The static pressure coefficient plots of the HLLEM scheme do not resemble potential flow but are closer to Stokes flow. This can be attributed to the high numerical dissipation in the normal momentum flux of the scheme at low Mach numbers. On the other hand, the static pressure coefficient plots of the proposed HLLEM-FP scheme, shown in Fig. \ref{cyl-hllem-fp}, resemble potential flow for all three Mach numbers. This is due to the reduction in numerical dissipation in the normal momentum flux through the addition of anti-diffusion terms in the scheme. 
The maximum pressure fluctuation obtained with the proposed HLL-CPS-FP and HLLEM-FP schemes for different inflow Mach numbers is plotted in Fig. \ref{cylinder-pfluc-hll-xfp}. It can be seen from the figure that the pressure fluctuation is of the order of Mach number squared ($\mathcal{O} (M^2)$) which is consistent with the pressure fluctuation in the continuous Euler equations.
\subsection{Flat plate boundary layer profile}
\begin{figure}[H]
	\begin{center}	
		\includegraphics[width=220pt]{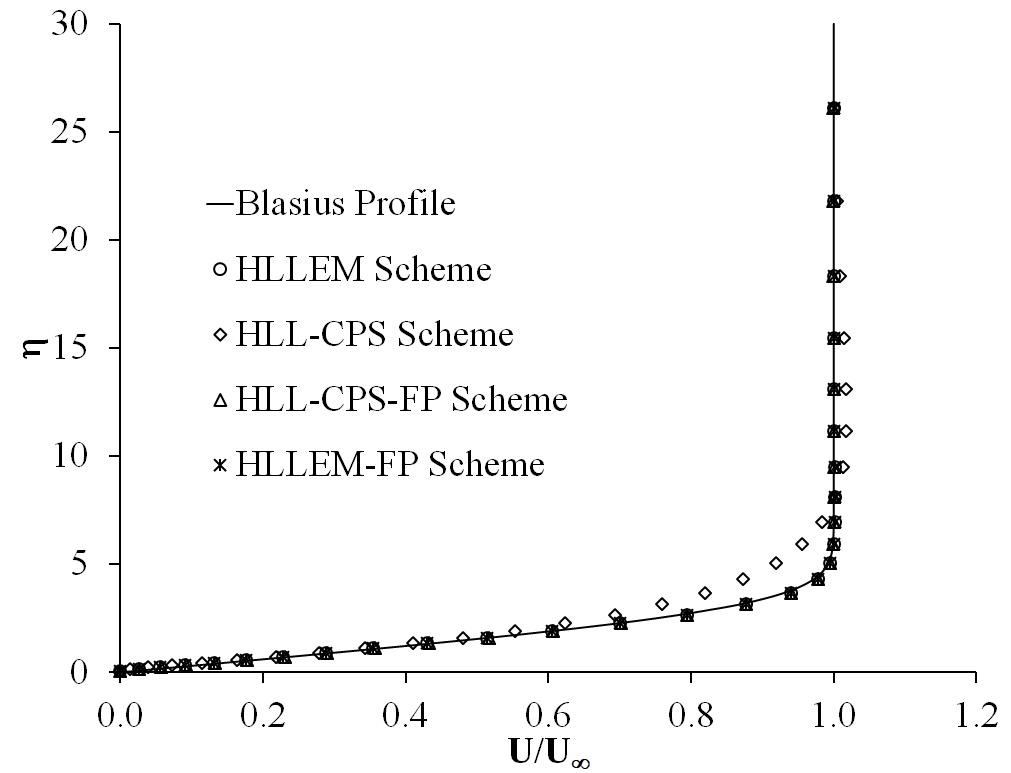}
		\caption{Boundary layer profile for $M_{\infty}=0.20$ laminar flow over a flat plate computed by the HLL-CPS-FP and HLLEM-FP Scheme }
		\label{blprofile_hllem_hllcps_fp}
	\end{center}
\end{figure}

Computations are carried out with the proposed schemes over a flat plat to demonstrate the boundary layer resolving capability of the modified schemes. The computations are carried out for a free-stream Mach number of 0.20 and Reynolds number of about 10,000 like in section \ref{bl-profile}. The third-order computations are carried out with the MUSCL approach \cite{van-leer-muscl} without any limiters. The time integration is carried out by the two-step Runge-Kutta method of Gottlieb and Shu \cite{gott}. The boundary layer profiles are shown after 50,000 iterations with a CFL number of 0.50. The grid had a size of $81 \times 33$. The boundary layer profiles of the proposed HLL-CPS-FP and HLLEM-FP schemes are shown in Fig. \ref{blprofile_hllem_hllcps_fp}. The normalized velocity $U/U_{\infty}$ is plotted against the Blasius parameter $\eta=y\sqrt{U/(\nu{}x)}$ in the figure. It can be seen from the figure that both the HLL-CPS-FP and HLLEM-FP schemes are capable of accurate boundary layer resolution, like the original HLLEM scheme.
\section{Conclusion}

This study compares the HLLEM and HLL-CPS approximate Riemann solvers and proposes enhancements for their performance in all-Mach number flows. The improvements involve modifying the anti-diffusion coefficients associated with contact and shear waves, as well as the anti-diffusion coefficient related to the velocity jump in the face normal direction.

The robustness and stability of the enhanced HLL-CPS and HLLEM schemes are validated through linear perturbation and matrix stability analyses. For high-speed flow problems, the new HLL-CPS scheme exhibits greater robustness compared to the original version, successfully resolving a flat plate boundary layer—a task that the original scheme struggled with. Additionally, the revised HLLEM scheme addresses numerical instabilities such as odd-even decoupling, kinked Mach stems, and the carbuncle phenomenon, which affected the original scheme.

In low Mach number scenarios, specifically for the case of a cylinder, both the modified HLL-CPS and HLLEM schemes outperform their original counterparts, producing results that closely resemble potential flow even at a very low Mach number of 0.001. Overall, the proposed HLL-CPS and HLLEM schemes demonstrate the capability to deliver accurate results across a wide range of Mach numbers.

\end{doublespace}

\textbf{Appendix I: Asymptotic Analysis of  HLL-CPS-FP Scheme for Low Mach Number Flows}\\
The HLL-CPS-FP scheme for the discrete Euler equations shown in (\ref{fv},\ref{hll-cps-flux2},\ref{bdq-new}) can be written in non-dimensional form, after rotation, as 
\begin{equation}
\begin{split}
& |\bar{\Omega}|_i\dfrac{\partial{}\bar{U}_i}{\partial{}\bar{t}}+  \sum_{l\epsilon\upsilon(i)}\left[\bar{u}_{n,il}\left[\begin{array}{c} \bar{\rho}\\ \bar{\rho}\bar{u}\\ \bar{\rho}\bar{v} \\ \bar{\rho}\bar{e}\end{array}\right]_K+\left[\begin{array}{c} 0 \\ \dfrac{\bar{p}}{M_*^2}n_x\\  \dfrac{\bar{p}}{M_*^2}n_y \\ \bar{p}\bar{u_n}\end{array}\right]_{il}+ \left(\dfrac{\bar{u}_nM_*}{2\bar{a}}\right)_{il}\left[\begin{array}{c} 0 \\ \left(\dfrac{n_x}{M_*^2}\right)_{il}\Delta_{il}(\bar{p})\\  \left(\dfrac{n_y}{M_*^2}\right)_{il} \Delta_{il}(\bar{p})\\ \Delta_{il}(\bar{p}\bar{u}_n)\end{array}\right]\right]\Delta{}s_{il}-  \\ &  \sum_{l\epsilon\upsilon(\mathbf{i})}\left(\dfrac{\bar{u}_n^2M_*}{2\bar{a}}-\dfrac{\bar{a}}{2M_*}\right)_{il}\left[\begin{array}{c} \Delta_{il}(\bar{\rho})-\left(\dfrac{\bar{a}}{\bar{u}_nM_*+\bar{a}}\right)_{il}(\Delta{}_{il}\bar{\rho}-\dfrac{\Delta_{il}\bar{p}}{\bar{a}^2_{il}}) \\ \Delta_{il}(\bar{\rho}\bar{u})-\left(\dfrac{\bar{a}}{\bar{u}_nM_*+\bar{a}}\right)_{il}\left(\bar{u}_{il}(\Delta{}_{il}\bar{\rho}-\dfrac{\Delta_{il}\bar{p}}{\bar{a}^2_{il}}) -(\bar{\rho}n_y)_{il}\Delta{}_{il}\bar{u}_t\right) \\ \Delta_{il}(\bar{\rho}\bar{v})-\left(\dfrac{\bar{a}}{\bar{u}_nM_*+\bar{a}}\right)_{il}\left(\bar{v}_{il}(\Delta{}_{il}\bar{\rho}-\dfrac{\Delta_{il}\bar{p}}{\bar{a}^2_{il}})+(\bar{\rho}n_x)_{il}\Delta{}_{il}\bar{u}_t\right) \\ \Delta_{il}(\bar{\rho}\bar{e})- \left(\dfrac{\bar{a}}{\bar{u}_nM_*+\bar{a}}\right)_{il}\left(\dfrac{\bar{q}^2_{il}}{2}(\Delta{}_{il}\bar{\rho}-\dfrac{\Delta_{il}\bar{p}}{\bar{a}^2_{il}}) +(\bar{\rho}\bar{u}_t)_{il}\Delta{}_{il}\bar{u}_t\right)\\\end{array}\right]\Delta{}s_{il} \\ & 
+\sum_{l\epsilon\upsilon(\mathbf{i})}\left(\dfrac{\bar{u}_n^2M_*}{2\bar{a}}-\dfrac{\bar{a}}{2M_*}\right)_{il}\left[\begin{array}{c} 0 \\ \bar{\delta}_n(\rho{}n_x)_{il}\Delta_{il}\bar{u}_n \\  \bar{\delta}_n(\rho{}n_y)_{il}\Delta_{il}\bar{u}_n \\ \bar{\delta}_n(\rho{}u_n)_{il}\Delta_{il}\bar{u}_n\end{array}\right]\Delta{}s_{il}=0
\end{split}
\end{equation}
The normal anti-diffusion coefficient in non-dimensional terms is $\bar{\delta}_n=1-M_*\bar{f(M)}$.\\
Expanding the terms asymptotically and arranging the terms in the equal power of $M_*$, we obtain 
\begin{enumerate}
	\item{Order of $M_*^0$}
	\begin{enumerate}
		\item{Continuity Equation}
		\begin{equation}
		|\bar{\Omega}|\dfrac{\partial\bar{\rho}_{0i}}{\partial{}\bar{t}}+\sum_{l\epsilon\upsilon(i)}\left(\bar{u}_{n0,il}\bar{\rho}_{0K}+\dfrac{1}{2\bar{a}_{0il}}\Delta_{il}\bar{p}_1\right)\Delta{}s_{il}=0
		\end{equation}
		\item{ x-momentum Equation}
		\begin{equation}
		\begin{split}
		& |\bar{\Omega}|\dfrac{\partial(\bar{\rho}_0\bar{u}_0)_i}{\partial{}\bar{t}}+ \sum_{l\epsilon\upsilon(i)}\left(\bar{u}_{n0,il}(\bar{\rho}_0\bar{u}_0)_K+\left(\bar{p}_2n_x\right)_{il}+\left(\dfrac{\bar{u}_{n0}}{2\bar{a}_0}n_x\right)_{il}\Delta_{il}\bar{p}_1 \right)\Delta{}s_{il}\\ & +\sum_{l\epsilon\upsilon(i)}\left(\dfrac{1}{2\bar{a}_{0il}}\Delta_{il}(\bar{p}_1\bar{u}_0+\bar{p}_0\bar{u}_1)\right)\Delta{}s_{il}=0
		\end{split}
		\end{equation}
		\item{ y-momentum Equation}
		\begin{equation}
		\begin{split}
		& |\bar{\Omega}|\dfrac{\partial(\bar{\rho}_0\bar{v}_0)_i}{\partial{}\bar{t}}+ \sum_{l\epsilon\upsilon(i)}\left(\bar{u}_{n0,il}(\bar{\rho}_0\bar{v}_0)_K+\left(\bar{p}_2n_y\right)_{il}+\left(\dfrac{\bar{u}_{n0}}{2\bar{a}_0}n_y\right)_{il}\Delta_{il}\bar{p}_1 \right)\Delta{}s_{il}\\ & +\sum_{l\epsilon\upsilon(i)}\left(\dfrac{1}{2\bar{a}_{0il}}\Delta_{il}(\bar{p}_1\bar{v}_0+\bar{p}_0\bar{v}_1)\right)\Delta{}s_{il}=0
		\end{split}
		\end{equation}
	\end{enumerate}
	\item{Order of $M_*^{-1}$} 
	\begin{enumerate}
		\item{Continuity equation}
		\begin{equation}\label{deltap0}
		\sum_{l\epsilon\upsilon(i)}\left(\dfrac{1}{2\bar{a}_{il}}\Delta_{il}\bar{p}_0\right)\Delta{}s_{il}=0
		\end{equation}
		From the above equation, we obtain $\sum_{l\epsilon\upsilon(i)}\Delta_{il}\bar{p}_0=0$ and hence $\bar{p}_0=$ constant for all $\mathbf{i}$.
		
		\item{x-momentum equation}
		\begin{equation}
		\begin{split}
		& \sum_{l\epsilon\upsilon(\mathbf{i})}\left((\bar{p}_1n_x)_{il}+\dfrac{\bar{u}_{n0}n_x}{2\bar{a}_0}_{il}\Delta_{il}\bar{p}_0\right)\Delta{}s_{il}+ \\ &
		\sum_{l\epsilon\upsilon(\mathbf{i})}\dfrac{a_{0il}}{2}\left(\Delta_{il}(\bar{\rho}_0\bar{u}_0)-\bar{u}_0(\Delta_{il}\bar{\rho}_0-\dfrac{\Delta_{il}\bar{p}_0}{\bar{a}^2_0})+(\bar{\rho}_0n_y)_{il}\Delta_{il}u_{t0}-(\bar{\rho}_0n_x)_{il}\Delta_{il}u_{n0}\right)\Delta{}s_{il}=0
		\end{split}
		\end{equation}
		Since $\sum_{l\epsilon\upsilon(\mathbf{i})}\Delta_{il}\bar{p}_0=0$, we obtain
		\begin{equation}
		\sum_{l\epsilon\upsilon(\mathbf{i})}(\bar{p}_1n_x)_{il}+\dfrac{\bar{a}_{0il}}{2}\left(\bar{\rho}_{0il}\Delta_{il}\bar{u}_0+\bar{\rho}_{0il}n_y\Delta_{il}\bar{u}_{t0}-\bar{\rho}_{0il}n_x\Delta_{il}\bar{u}_{n0}\right)=0
		\end{equation}
		
		\begin{equation} \label{result-om1-hllcpsfp-x}
		\sum_{l\epsilon\upsilon(\mathbf{i})}(\bar{p}_1n_x)_{il}=0
		\end{equation}
		
		\item{ y-momentum equation}
		\begin{equation}
		\begin{split}
		& \sum_{l\epsilon\upsilon(\mathbf{i})}\left((\bar{p}_1n_y)_{il}+\dfrac{\bar{u}_{n0}n_y}{2\bar{a}_0}_{il}\Delta_{il}\bar{p}_0+\dfrac{\bar{u}_{0il}}{2\bar{a}_{0il}}\Delta_{il}(\bar{p}_0)\right)\Delta{}s_{il}+ \\ &
		\sum_{l\epsilon\upsilon(\mathbf{i})}\dfrac{a_{0il}}{2}\left(\Delta_{il}(\bar{\rho}_0\bar{v}_0)-\bar{v}_0(\Delta_{il}\bar{\rho}_0-\dfrac{\Delta_{il}\bar{p}_0}{\bar{a}^2_{0il}})-(\bar{\rho}_0n_x)_{il}\Delta_{il}u_{t0}-(\bar{\rho}_0n_y)_{il}\Delta_{il}u_{n0}\right)\Delta{}s_{il}=0
		\end{split}
		\end{equation}
		Since $\sum_{l\epsilon\upsilon(\mathbf{i})}\Delta_{il}\bar{p}_0=0$, we obtain
		\begin{equation} 
		\sum_{l\epsilon\upsilon(\mathbf{i})}(\bar{p}_1n_y)_{il}+\dfrac{\bar{a}_{0il}}{2}\left(\bar{\rho}_{0il}\Delta_{il}\bar{v}_0-\bar{\rho}_{0il}n_x\Delta_{il}\bar{u}_{t0}-\bar{\rho}_{0il}n_y\Delta_{il}\bar{u}_{n0}\right)=0
		\end{equation}
		
		\begin{equation} \label{result-om1-hllcpsfp-y}
		\sum_{l\epsilon\upsilon(\mathbf{i})}(\bar{p}_1n_y)_{il}=0
		\end{equation}
		
	\end{enumerate}
\end{enumerate}

Equations (\ref{result-om1-hllcpsfp-x}, \ref{result-om1-hllcpsfp-y}) imply that 
\begin{equation}
\bar{p}_{1,i-1,j}-\bar{p}_{1,i+1,j}=0, \hspace{1cm} \bar{p}_{1,i,j-1}-\bar{p}_{1,i,j+1}=0, 
\end{equation}
This implies that $\bar{p}_1=$constant for all $\mathbf{i}$. Therefore, the HLL-CPS-FP scheme permits pressure fluctuation of the type $p(x,t)=p_0(t)+M_*^2p_2(x,t)$ and hence, the scheme will be able to resolve the flow features at low Mach numbers. \\

\textbf{Appendix II: Asymptotic Analysis of  HLLEM Scheme for Low Mach Number Flows}\\
The HLLEM scheme shown in equations (\ref{fv}, \ref{hll-type-flux},\ref{bdq-roe-hllem}) can be written in non-dimensional form as 
\begin{equation}
\begin{split}
& |\bar{\Omega}|\dfrac{\partial{}\bar{U}_i}{\partial{}\bar{t}}+  \sum_{l\epsilon\upsilon(\mathbf{i})}\left[\begin{array}{c} \bar{\rho}\bar{u}_n \\ \bar{\rho}\bar{u}\bar{u}_n+\dfrac{\bar{p}}{M_*^2}n_x\\  \bar{\rho}\bar{v}\bar{u}_n+\dfrac{\bar{p}}{M_*^2}n_y \\ (\bar{\rho}\bar{e}+\bar{p})\bar{u_n}\end{array}\right]_{il}+ \left(\dfrac{\bar{u}_nM_*}{2\bar{a}}\right)_{il}\left[\begin{array}{c} \Delta_{il}(\bar{\rho}\bar{u}_n) \\ \Delta_{il}(\bar{\rho}\bar{u}\bar{u}_n)+\left(\dfrac{n_x}{M_*^2}\right)_{il}\Delta_{il}(\bar{p})\\ \Delta_{il}(\bar{\rho}\bar{v}\bar{u}_n)+ \left(\dfrac{n_y}{M_*^2}\right)_{il} \Delta_{il}(\bar{p})\\ \Delta_{il}(\bar{\rho}\bar{e}\bar{u}_n+\bar{p}\bar{u}_n)\end{array}\right]\Delta{}s_{il}-  \\ & \sum_{l\epsilon\upsilon(\mathbf{i})}\left(\dfrac{\bar{u}_n^2M_*}{2\bar{a}}-\dfrac{\bar{a}}{2M_*}\right)_{il}\left[\begin{array}{c} \Delta_{il}(\bar{\rho})-\left(\dfrac{\bar{a}}{\bar{u}_nM_*+\bar{a}}\right)_{il}(\Delta{}_{il}\bar{\rho}-\dfrac{\Delta_{il}\bar{p}}{\bar{a}^2_{il}}) \\ \Delta_{il}(\bar{\rho}\bar{u})-\left(\dfrac{\bar{a}}{\bar{u}_nM_*+\bar{a}}\right)_{il}\left(\bar{u}_{il}(\Delta{}_{il}\bar{\rho}-\dfrac{\Delta_{il}\bar{p}}{\bar{a}^2_{il}}) -(\bar{\rho}n_y)_{il}\Delta{}_{il}\bar{u}_t\right) \\ \Delta_{il}(\bar{\rho}\bar{v})-\left(\dfrac{\bar{a}}{\bar{u}_nM_*+\bar{a}}\right)_{il}\left(\bar{v}_{il}(\Delta{}_{il}\bar{\rho}-\dfrac{\Delta_{il}\bar{p}}{\bar{a}^2_{il}})+(\bar{\rho}n_x)_{il}\Delta{}_{il}\bar{u}_t\right) \\ \Delta_{il}(\bar{\rho}\bar{e})- \left(\dfrac{\bar{a}}{\bar{u}_nM_*+\bar{a}}\right)_{il}\left(\dfrac{\bar{q}^2_{il}}{2}(\Delta{}_{il}\bar{\rho}-\dfrac{\Delta_{il}\bar{p}}{\bar{a}^2_{il}}) +(\bar{\rho}\bar{u}_t)_{il}\Delta{}_{il}\bar{u}_t\right)\\\end{array}\right]\Delta{}s_{il}=0 
\end{split}
\end{equation}

The flow variables are expanded using the following asymptotic expansions 
\begin{equation} \label{expansion}
\bar{\rho}=\bar{\rho_0}+M_*\bar{\rho_1}+M_*^2\bar{\rho_2}, \hspace{1cm} \bar{u}_n=\bar{u}_{n0}+M_*\bar{u}_{n1}+M_*^2\bar{u}_{n2}, \hspace{1cm} \bar{p}=\bar{p}_0+M_*\bar{p}_1+M_*^2\bar{p}_2 
\end{equation}
As per the rules of Rieper\cite{rieper2}, only the leading terms are considered for the terms appearing in the denominator. Expanding the terms asymptotically using equation (\ref{expansion}) and arranging the terms in the equal power of $M_*$, we obtain 
\begin{enumerate}
\item{Order of $M_*^0$}
\begin{enumerate}
\item{Continuity Equation}
\begin{equation}
|\bar{\Omega}|\dfrac{\partial\bar{\rho}_{0i}}{\partial{}\bar{t}}+\sum_{l\epsilon\upsilon(\mathbf{i})}\left((\bar{\rho}\bar{u}_{n0})_{il}+\dfrac{\bar{u}_{n0il}}{2}\Delta_{il}\bar{\rho}_0-\dfrac{\bar{u}_{n0il}}{2\bar{a}^2}\Delta_{il}\bar{p}_0+\dfrac{1}{2\bar{a}_{0il}}\Delta_{il}\bar{p}_1\right)\Delta{}s_{il}=0
\end{equation}
\item{x-Momentum Equation}
\begin{equation}
\begin{split}
& |\bar{\Omega}|\dfrac{\partial(\bar{\rho}_0\bar{u}_0)_i}{\partial{}\bar{t}}+\sum_{l\epsilon\upsilon(\mathbf{i})}\left((\bar{\rho}_{0}\bar{u}_0\bar{u}_{n0})_{il}+(\bar{p}_2n_x)_{il}+\dfrac{n_x}{2\bar{a}_0}_{il}(\bar{u}_{n1}\Delta_{il}\bar{p}_0+\bar{u}_{n0}\Delta_{il}\bar{p}_1)\right)\Delta{}s_{il} \\ &
+\sum_{l\epsilon\upsilon(\mathbf{i})}\dfrac{\bar{a}_{0il}}{2}(\bar{\rho}_0\Delta_{il}\bar{u}_1+\bar{u}_1\Delta_{il}\bar{\rho}_0+\bar{\rho}_1\Delta_{il}\bar{u}_0+\bar{u}_0\Delta_{il}\bar{\rho}_1)\Delta{}s_{il}\\ & 
+\sum_{l\epsilon\upsilon(\mathbf{i})}\dfrac{\bar{u}_{n0}}{2}\left(\bar{u}_0(\Delta_{il}\bar{\rho}_0-\dfrac{\Delta_{il}p_0}{\bar{a}^2_0})-(\bar{\rho}_0n_y)_{il}\Delta{}_{il}\bar{u}_{t0}\right) \Delta{}s_{il}\\&
-\sum_{l\epsilon\upsilon(\mathbf{i})}\dfrac{\bar{a}_{0il}}{2}\left(\bar{u}_{0il}\Delta_{il}\bar{\rho}_1+\bar{u}_{1il}\Delta_{il}\bar{\rho}_0-\bar{u}_{0il}\dfrac{\Delta_{il}\bar{p}_1}{\bar{a}^2}-\bar{u}_{1il}\dfrac{\Delta_{il}\bar{p}_0}{\bar{a}^2}\right) \Delta{}s_{il}\\ &
+\sum_{l\epsilon\upsilon(\mathbf{i})}\dfrac{\bar{a}_{0il}}{2}((\bar{\rho}_0n_y)_{il}\Delta{}_{il}\bar{u}_{t1}+(\bar{\rho}_1n_y)_{il}\Delta{}_{il}\bar{u}_{t0})\Delta{}s_{il}=0
\end{split}
\end{equation}
\item{y-Momentum Equation}
\begin{equation}
\begin{split}
& |\bar{\Omega}|\dfrac{\partial(\bar{\rho}_0\bar{v}_0)_i}{\partial{}\bar{t}}+\sum_{l\epsilon\upsilon(\mathbf{i})}\left((\bar{\rho}_{0}\bar{v}_0\bar{u}_{n0})_{il}+(\bar{p}_2n_y)_{il}+\left(\dfrac{n_y}{2\bar{a}_0}\right)_{il}(\bar{u}_{n1}\Delta_{il}\bar{p}_0+\bar{u}_{n0}\Delta_{il}\bar{p}_1)\right)\Delta{}s_{il} \\ &
+\sum_{l\epsilon\upsilon(\mathbf{i})}\dfrac{\bar{a}_{0il}}{2}(\bar{\rho}_0\Delta_{il}\bar{v}_1+\bar{v}_1\Delta_{il}\bar{\rho}_0+\bar{\rho}_1\Delta_{il}\bar{v}_0+\bar{v}_0\Delta_{il}\bar{\rho}_1)\Delta{}s_{il}\\ & 
+\sum_{l\epsilon\upsilon(\mathbf{i})}\dfrac{\bar{u}_{n0}}{2}\left(\bar{v}_0\Delta_{il}\bar{\rho}_0-\dfrac{v_0\Delta_{il}p_0}{\bar{a}^2_0}+(\bar{\rho}_0n_x)_{il}\Delta{}_{il}\bar{u}_{t0}\right) \Delta{}s_{il}\\&
-\sum_{l\epsilon\upsilon(\mathbf{i})}\dfrac{\bar{a}_{0il}}{2}\left(\bar{v}_{0il}\Delta_{il}\bar{\rho}_1+\bar{v}_{1il}\Delta_{il}\bar{\rho}_0-\bar{v}_{0il}\dfrac{\Delta_{il}\bar{p}_1}{\bar{a}^2}-\bar{v}_{1il}\dfrac{\Delta_{il}\bar{p}_0}{\bar{a}^2}\right) \Delta{}s_{il}\\ &
-\sum_{l\epsilon\upsilon(\mathbf{i})}\dfrac{\bar{a}_{0il}}{2}((\bar{\rho}_0n_x)_{il}\Delta{}_{il}\bar{u}_{t1}+(\bar{\rho}_1n_x)_{il}\Delta{}_{il}\bar{u}_{t0})\Delta{}s_{il}=0
\end{split}
\end{equation}
\end{enumerate}
\item{Order of $M_*^{-1}$ }
\begin{enumerate}
\item{Continuity Equation}
\begin{equation}
\sum_{l\epsilon\upsilon(\mathbf{i})}\dfrac{1}{2\bar{a}_{0il}}\Delta_{il}\bar{p}_0\Delta{}s_{il}=0
\end{equation}
Therefore,  $\sum_{l\epsilon\upsilon(\mathbf{i})}\Delta_{il}\bar{p}_0=0$ and hence $\bar{p}_0=$ constant for all $\mathbf{i}$
\item{x-Momentum Equation}
\begin{equation}
\begin{split}
& \sum_{l\epsilon\upsilon(\mathbf{i})}\left((\bar{p}_1n_x)_{il}+\dfrac{\bar{u}_{n0}n_x}{2\bar{a}_0}_{il}\Delta_{il}\bar{p}_0\right)\Delta{}s_{il} \\ &
+ \sum_{l\epsilon\upsilon(\mathbf{i})}\dfrac{a_{0il}}{2}\left(\Delta_{il}(\bar{\rho}_0\bar{u}_0)-\bar{u}_0(\Delta_{il}\bar{\rho}_0-\dfrac{\Delta_{il}\bar{p}_0}{\bar{a}^2_0})+(\bar{\rho}_0n_y)_{il}\Delta_{il}u_{t0}\right)\Delta{}s_{il}=0
\end{split}
\end{equation}
Since $\sum_{l\epsilon\upsilon(\mathbf{i})}\Delta_{il}\bar{p}_0=0$, we obtain
\begin{equation}\label{result-om1-hllem-x}
\sum_{l\epsilon\upsilon(\mathbf{i})}(\bar{p}_1n_x)_{il}+\dfrac{\bar{a}_{0il}}{2}\left(\bar{\rho}_{0il}\Delta_{il}\bar{u}_0+\bar{\rho}_{0il}n_y\Delta_{il}\bar{u}_{t0}\right)=0
\end{equation}
\item{y-Momentum Equation}
\begin{equation}
\begin{split}
& \sum_{l\epsilon\upsilon(\mathbf{i})}\left((\bar{p}_1n_y)_{il}+\dfrac{\bar{u}_{n0}n_y}{2\bar{a}_0}_{il}\Delta_{il}\bar{p}_0+\dfrac{\bar{u}_{0il}}{2\bar{a}_{0il}}\Delta_{il}(\bar{p}_0)\right)\Delta{}s_{il} \\ &
+ \sum_{l\epsilon\upsilon(\mathbf{i})}\dfrac{a_{0il}}{2}\left(\Delta_{il}(\bar{\rho}_0\bar{v}_0)-\bar{v}_0(\Delta_{il}\bar{\rho}_0-\dfrac{\Delta_{il}\bar{p}_0}{\bar{a}^2_{0il}})-(\bar{\rho}_0n_x)_{il}\Delta_{il}u_{t0}\right)\Delta{}s_{il}=0
\end{split}
\end{equation}
Since $\sum_{l\epsilon\upsilon(\mathbf{i})}\Delta_{il}\bar{p}_0=0$, we obtain
\begin{equation} \label{result-om1-hllem-y}
\sum_{l\epsilon\upsilon(\mathbf{i})}(\bar{p}_1n_y)_{il}+\dfrac{\bar{a}_{0il}}{2}\left(\bar{\rho}_{0il}\Delta_{il}\bar{v}_0-\bar{\rho}_{0il}n_x\Delta_{il}\bar{u}_{t0}\right)=0
\end{equation}
\end{enumerate}
\end{enumerate}

Equations (\ref{result-om1-hllem-x}) and (\ref{result-om1-hllem-y}) imply that  

\begin{equation} \label{hllem-p1}
\sum_{l\epsilon\upsilon(\mathbf{i})}\bar{p}_1 \neq 0
\end{equation}

Therefore, the  HLLEM scheme for the discrete Euler equations supports pressure fluctuation of the type $p(x,t)=p_0(t)+M_*p_1(x,t)$ and hence shall not be able to resolve flow features at low Mach numbers.

The HLLEM-FP scheme shown in equations (\ref{fv}, \ref{hll-type-flux},\ref{bdq-new}) can be written in non-dimensional form, after rotation,  as 

\begin{equation}
\begin{split}
& |\bar{\Omega}|\dfrac{\partial{}\bar{U}_i}{\partial{}\bar{t}}+  \sum_{l\epsilon\upsilon(\mathbf{i})}\left[\begin{array}{c} \bar{\rho}\bar{u}_n \\ \bar{\rho}\bar{u}\bar{u}_n+\dfrac{\bar{p}}{M_*^2}n_x\\  \bar{\rho}\bar{v}\bar{u}_n+\dfrac{\bar{p}}{M_*^2}n_y \\ (\bar{\rho}\bar{e}+\bar{p})\bar{u_n}\end{array}\right]_{il}+ \left(\dfrac{\bar{u}_nM_*}{2\bar{a}}\right)_{il}\left[\begin{array}{c} \Delta_{il}(\bar{\rho}\bar{u}_n) \\ \Delta_{il}(\bar{\rho}\bar{u}\bar{u}_n)+\left(\dfrac{n_x}{M_*^2}\right)_{il}\Delta_{il}(\bar{p})\\ \Delta_{il}(\bar{\rho}\bar{v}\bar{u}_n)+ \left(\dfrac{n_y}{M_*^2}\right)_{il} \Delta_{il}(\bar{p})\\ \Delta_{il}(\bar{\rho}\bar{e}\bar{u}_n+\bar{p}\bar{u}_n)\end{array}\right]\Delta{}s_{il}-  
\\ &  \sum_{l\epsilon\upsilon(\mathbf{i})}\left(\dfrac{\bar{u}_n^2M_*}{2\bar{a}}-\dfrac{\bar{a}}{2M_*}\right)_{il}\left[\begin{array}{c} \Delta_{il}(\bar{\rho})-\left(\dfrac{\bar{a}}{\bar{u}_nM_*+\bar{a}}\right)_{il}(\Delta{}_{il}\bar{\rho}-\dfrac{\Delta_{il}\bar{p}}{\bar{a}^2_{il}}) \\ \Delta_{il}(\bar{\rho}\bar{u})-\left(\dfrac{\bar{a}}{\bar{u}_nM_*+\bar{a}}\right)_{il}\left(\bar{u}_{il}(\Delta{}_{il}\bar{\rho}-\dfrac{\Delta_{il}\bar{p}}{\bar{a}^2_{il}}) -(\bar{\rho}n_y)_{il}\Delta{}_{il}\bar{u}_t\right) \\ \Delta_{il}(\bar{\rho}\bar{v})-\left(\dfrac{\bar{a}}{\bar{u}_nM_*+\bar{a}}\right)_{il}\left(\bar{v}_{il}(\Delta{}_{il}\bar{\rho}-\dfrac{\Delta_{il}\bar{p}}{\bar{a}^2_{il}})+(\bar{\rho}n_x)_{il}\Delta{}_{il}\bar{u}_t\right) \\ \Delta_{il}(\bar{\rho}\bar{e})- \left(\dfrac{\bar{a}}{\bar{u}_nM_*+\bar{a}}\right)_{il}\left(\dfrac{\bar{q}^2_{il}}{2}(\Delta{}_{il}\bar{\rho}-\dfrac{\Delta_{il}\bar{p}}{\bar{a}^2_{il}}) +(\bar{\rho}\bar{u}_t)_{il}\Delta{}_{il}\bar{u}_t\right)\\\end{array}\right]\Delta{}s_{il} \\ & 
+\sum_{l\epsilon\upsilon(\mathbf{i})}\left(\dfrac{\bar{u}_n^2M_*}{2\bar{a}}-\dfrac{\bar{a}}{2M_*}\right)_{il}\left[\begin{array}{c} 0 \\ \bar{\delta}_n(\rho{}n_x)_{il}\Delta_{il}\bar{u}_n \\  \bar{\delta}_n(\rho{}n_y)_{il}\Delta_{il}\bar{u}_n \\ \bar{\delta}_n(\rho{}u_n)_{il}\Delta_{il}\bar{u}_n\end{array}\right]\Delta{}s_{il}=0
\end{split}
\end{equation}
The normal anti-diffusion coefficient in non-dimensional terms is $\bar{\delta}_n=1-M_*\bar{f(M)}$.

Expanding the terms asymptotically using equation (\ref{expansion}) and arranging the terms in the equal power of $M_*$, we obtain 
\begin{enumerate}
\item{Order of $M_*^{-1}$ }
\begin{enumerate}
\item{Continuity Equation}
\begin{equation}
\sum_{l\epsilon\upsilon(\mathbf{i})}\dfrac{1}{2\bar{a}_{0il}}\Delta_{il}\bar{p}_0\Delta{}s_{il}=0
\end{equation}
Therefore,  $\sum_{l\epsilon\upsilon(\mathbf{i})}\Delta_{il}\bar{p}_0=0$ and hence $\bar{p}_0=Cte \forall i$ where $Cte$ is a constant.
\item{x-Momentum Equation}
\begin{equation}
\begin{split}
& \sum_{l\epsilon\upsilon(\mathbf{i})}\left((\bar{p}_1n_x)_{il}+\dfrac{\bar{u}_{n0}n_x}{2\bar{a}_0}_{il}\Delta_{il}\bar{p}_0\right)\Delta{}s_{il}+ \\ &
\sum_{l\epsilon\upsilon(\mathbf{i})}\dfrac{a_{0il}}{2}\left(\Delta_{il}(\bar{\rho}_0\bar{u}_0)-\bar{u}_0(\Delta_{il}\bar{\rho}_0-\dfrac{\Delta_{il}\bar{p}_0}{\bar{a}^2_0})+(\bar{\rho}_0n_y)_{il}\Delta_{il}u_{t0}-(\bar{\rho}_0n_x)_{il}\Delta_{il}u_{n0}\right)\Delta{}s_{il}=0
\end{split}
\end{equation}
Since $\sum_{l\epsilon\upsilon(\mathbf{i})}\Delta_{il}\bar{p}_0=0$, we obtain
\begin{equation}
\sum_{l\epsilon\upsilon(\mathbf{i})}(\bar{p}_1n_x)_{il}+\dfrac{\bar{a}_{0il}}{2}\left(\bar{\rho}_{0il}\Delta_{il}\bar{u}_0+\bar{\rho}_{0il}n_y\Delta_{il}\bar{u}_{t0}-\bar{\rho}_{0il}n_x\Delta_{il}\bar{u}_{n0}\right)=0
\end{equation}

\begin{equation} \label{result-om1-hllemfp-x}
\sum_{l\epsilon\upsilon(\mathbf{i})}(\bar{p}_1n_x)_{il}=0
\end{equation}

\item{y-Momentum Equation}
\begin{equation}
\begin{split}
& \sum_{l\epsilon\upsilon(\mathbf{i})}\left((\bar{p}_1n_y)_{il}+\dfrac{\bar{u}_{n0}n_y}{2\bar{a}_0}_{il}\Delta_{il}\bar{p}_0+\dfrac{\bar{u}_{0il}}{2\bar{a}_{0il}}\Delta_{il}(\bar{p}_0)\right)\Delta{}s_{il}+ \\ &
\sum_{l\epsilon\upsilon(\mathbf{i})}\dfrac{a_{0il}}{2}\left(\Delta_{il}(\bar{\rho}_0\bar{v}_0)-\bar{v}_0(\Delta_{il}\bar{\rho}_0-\dfrac{\Delta_{il}\bar{p}_0}{\bar{a}^2_{0il}})-(\bar{\rho}_0n_x)_{il}\Delta_{il}u_{t0}-(\bar{\rho}_0n_y)_{il}\Delta_{il}u_{n0}\right)\Delta{}s_{il}=0
\end{split}
\end{equation}
Since $\sum_{l\epsilon\upsilon(\mathbf{i})}\Delta_{il}\bar{p}_0=0$, we obtain
\begin{equation} 
\sum_{l\epsilon\upsilon(\mathbf{i})}(\bar{p}_1n_y)_{il}+\dfrac{\bar{a}_{0il}}{2}\left(\bar{\rho}_{0il}\Delta_{il}\bar{v}_0-\bar{\rho}_{0il}n_x\Delta_{il}\bar{u}_{t0}-\bar{\rho}_{0il}n_y\Delta_{il}\bar{u}_{n0}\right)=0
\end{equation}

\begin{equation} \label{result-om1-hllemfp-y}
\sum_{l\epsilon\upsilon(\mathbf{i})}(\bar{p}_1n_y)_{il}=0
\end{equation}

Equations (\ref{result-om1-hllemfp-x}, \ref{result-om1-hllemfp-y}) imply that 

\begin{equation}
\bar{p}_{1,i-1,j}-\bar{p}_{1,i+1,j}=0, \hspace{1cm} \bar{p}_{1,i,j-1}-\bar{p}_{1,i,j+1}=0, 
\end{equation}
This implies that $\bar{p}_1=$constant for all $\mathbf{i}$

\end{enumerate}
\end{enumerate}

Therefore, the  HLLEM-FP scheme for the discrete Euler equations supports pressure fluctuation of the type $p(x,t)=p_0(t)+M_*^2p_1(x,t)$ and hence is consistent with the continuous Euler equations.


\begin{thebibliography}{9}
	\bibitem{godu} S. R. Godunov, A difference Scheme for Numerical Solution of Discontinuous Solution of Hydrodynamic Equations, Math. Sbornik, 47, 271-306, translated US Joint Publ. Res. Service, JPRS 7226, 1969.
	\bibitem{toro-book} E. F. Toro, Riemann Solvers and Numerical Methods for fluid dynamics, Third Edition, Springer-Verlag, Berlin Heidelberg, 2009.
	\bibitem{hll} A. Harten, P. D. Lax, B. van Leer, On Upstream Differencing and Godunov-type Methods for Hyperbolic Conservation Laws, SIAM Review 1983, 25, 35-61.
	\bibitem{einf1} B. Einfeldt, On Godunov-type Methods for Gas Dynamics, SIAM J. Numerical Analysis, Vol 25, No 2, April 1988.
	\bibitem{roe} P. L. Roe, Approximate Riemann solvers, parameter vectors and difference schemes, J. Comput. Phys. 43 (1981) 357-372.
	\bibitem{osher} S. Osher, F. Solomon, Upwind Difference Schemes for Hyperbolic System of Conservation Laws, Mathematics of Computation, Vol. 38, Num. 158 (1982) 339-374. 
	\bibitem{einf2} B. Einfeldt, C. D. Munz, P. L. Roe, On Godunov type methods near low densities , J. Comput. Phys. 92 (1991) 273-295.
	\bibitem{toro1} E. F. Toro, M. Spruce, W. Speares, Restoration of the contact surface in the HLL-Riemann solver, Shock Waves (1994) 4:25-34.
	\bibitem{batt} P. Batten, N. Clarke, C. Lambert, D. M. Causon, On the choice of wavespeeds for the HLLC Riemann Solver, SIAM J. Sci. Comput. Vol. 18, No 6, pp 1553-1570, November 1997.
	\bibitem{mandal} J. C. Mandal , V. Panwar, Robust HLL type Riemann Solver Capable of Resolving Contact Discontinuity, Comput.  Fluids 63 (2012) 148-164.
	\bibitem{kita1} K. Kitamura, A Further Survey of Shock Capturing Methods on Hypersonic Heating Issues, AIAA-2013-2698.
	\bibitem{kita-slau2} K. Kitamura, Assessment of SLAU2 and other flux functions with slope limiters in hypersonic shock-interaction heating, computers and Fluids 129(2016) 134-145.
	\bibitem{peery} K. M. Peery, S. T. Imlay, Blunt Body Flow Simulations, AIAA Paper-88-2924, 1988.
	\bibitem{liou1}  M-S. Liou, Mass Flux schemes and connection to shock instability, J. Comput. Phys. 160 (2000) 623-648.
	\bibitem{kim}  S-S. Kim, C. Kim, O-H. Rho, S. K. Hong, Cures for the shock instability: Development of a shock stable Roe scheme, J. Comput. Phys. 185 (2003) 342-374.
	\bibitem{quirk} J. J. Quirk, A contribution to the great Riemann solver debate, Int. J. Numer. Meth. Fluids 18 (1994) 555-574.
	\bibitem{pand} M. Pandolfi, D. D'Ambrosio, Numerical Instabilities in Upwind Methods: Analysis and Cures for the ``Carbuncle'' Phenomenon, J. Comput. Phys. 166 (2001) 271-301.
	\bibitem{gress1} J. Gressier, J-M. Moschetta, Robustness versus accuracy in shock-wave computations, Int. J. Numer. Methods Fluids 33 (2000) 313-332.
	\bibitem{dumb-matrix} M. Dumbser, J-M. Moschetta, J. Gressier, A matrix stability analysis of the carbuncle problem, J. Comput. Phys. 197 (2004) 647-670.
	\bibitem{ren} Y-X. Ren, A robust shock-capturing scheme based on rotated Riemann solvers, Comput. Fluids 32 (2003) 1379-1403.
	\bibitem{nishi} H. Nishikawa, K. Kitamura, Very simple, carbuncle-free, boundary-layer-resolving, rotated-hybrid Riemann solvers, J. Comput. Phys. 227 (2008) 2560-2581.
	\bibitem{park} S. H. Park, J. H. Kwon, On the Dissipation mechanism of Godunov-type schemes J. Comput. Phys. 188 (2003) 524-542.
	\bibitem{shen}  Z. Shen, W. Yan, G. Yuan, A robust HLLC-type Riemann solver for strong shock, J. Comput. Phys.  309 (2016) 185-206.
	\bibitem{xie-hllem} W. Xie, W. Li, H. Li, Z. Tian and S. Pan, On numerical instabilities of Godunov-type scheme for strong shock, J. Comput. Phys., 350 (2017) 607-637.
	\bibitem{chen-roe} S-S Chen, C. Yen, B-X Lin and Y-S. Li, A new Robust Carbuncle-free Roe scheme for Strong Shock, J. Sci. Comput. 2018(77) 1250-1277.
	\bibitem{chen-all-scheme}  S-S Chen, C. Yan, B-X Lin, L-Y Liu, J. Yu, Affordable shock-stable item for Godunov-type schemes against carbuncle phenomenon, J. Comput. Phys., 373 (2018) 662-672.
	\bibitem{san1} S. Simon, J. C. Mandal, A simple cure for numerical shock instability in the HLLC Riemann solver, J. Comput. Phys. 378 (2018) 477-496.
	\bibitem{san2} S. Simon, J. C. Mandal, A cure for numerical shock instability in HLLC Riemann solver using anti-diffusion control, Comput. Fluids 174 (2018) 144-166.
	\bibitem{san3} S. Simon, J. C. Mandal, Strategies to cure numerical shock instability in the HLLEM Riemann solver, Int. J. Numer. Meth. Fluids 89 (2019) 533-569.
	\bibitem{kemm} F. Kemm, Heuristic and Numerical Considerations for the Carbuncle Phenomenon, Applied Mathematics and Computation 320 (2018) 596-613.	
	\bibitem{fleis1} N. Fleischmann, S. Adami,  X. Y. Hu, N. A. Adams, A low dissipation method to cure the grid-aligned shock instability, J. Comput. Phys. 410 (2020) 109004.
	\bibitem{fleis2} N. Fleischmann, S. Adami, N. A. Adams, A shock-stable modification of  the HLLC Riemann solver with reduced numerical dissipation, J. Comput. Phys. 423 (2020) 109762.
	\bibitem{guillard1} H. Guillard, C. Viozat, On the Behaviour of Upwind Schemes in the Low Mach Limit, Comput. Fluids 28 (1999) 63-86.
	\bibitem{rieper1} F. Rieper, On the Behaviour of Numerical Scheme in the Low Mach Number Regime, PhD Dissertation, July 2008.
	\bibitem{rieper2} F. Rieper, A low-Mach number fix for Roe's Approximate Riemann Solver, J. Comput. Phys. 230 (2011) 5263-5287.
	\bibitem{park-pre} S. H. Park, J. E. Lee, J. H. Kwon, Preconditioned HLLE Method for Flows at All Mach Numbers, AIAA Journal, 2006 44:11 pp. 2645-2653.
	\bibitem{li-roe} X-S. Li, C-W. Gu, An All-Speed Roe-Type scheme and its Asymptotic Analysis of Low Mach Number behaviour, J. Comput. Phys. 227 (2008) 5144-5159.
	\bibitem{th1} B. Thornber, A. Mosedale, D. Drikakis, D. Youngs, R. Williams, An Improved Reconstruction Method for Compressible Flows with Low Mach Number Features, J. Comput. Phys 227 (2008) 4873-4894.
	\bibitem{th3} A. Garcia-Ucida Juarez, A. Raimo, E. Shapiro and B. Thornber, Steady Turbulent Flow Computations Using a Low Mach Fully Compressible Scheme, AIAA Journal 2014 52:11 pp. 2559-2575.
	\bibitem{ossw} K. Osswald, A. Siegmund, P. Birken, V. Hannemann, A. Meister, $L^2$Roe: A low-dissipation version of Roe's approximate Riemann solver for low Mach numbers, Int. J. Numer. Meth. Fluids 81 (2016) 71-86.
	\bibitem{della1} S. Dellacherie, Analysis of Godunov type schemes applied to compressible Euler system at low Mach number, J. Comput. Phys. 229 (2010) 978-1016.
	\bibitem{della2} S. Dellacherie, J. Jung, P. Omnes, P-A. Raviart, Construction of modified Godunov-type schemes accurate at any Mach number for the compressible Euler system, Mathematical Models and Methods in Applied Sciences, Vol 26, No 13 (2016) 2525-2615.
	\bibitem{yu-hllem2} H. Yu, Z. Tian, F. Yang and H. Li, On Asymptotic Behaviour of HLL-type Schemes at Low Mach numbers, Mathematical Problems in Engineering, Volume 2020, Article ID 7451240.
	\bibitem {xie-roe} W. Xie, Y. Zhang, Q. Chang, H. Li, Towards an Accurate and Robust Roe-Type Scheme for All Mach Number Flows, Adv. Appl. Math. Mech. Vol 11, No 1, pp.132-167 (2019). 
	\bibitem{shima} E. Shima, K. Kitamura, Parameter-Free Simple Low-Dissipation AUSM-Family Scheme for All Speeds, AIAA Journal, 2011 49:8 pp. 1693-1709.
	\bibitem{chen-ausm} S-S. Chen, F-J. Cai, H-C Xue, N. Wang, C. Yan, An improved AUSM-family scheme with robustness and accuracy for all Mach number flows, Applied Mathematical Modeling, 77 (2020) 1065-1081.
	\bibitem{xie-hllc} W. Xie, Y. Zhang, J. Lai, H. Li, An accurate and robust HLLC-type Riemann solver for the compressible Euler system at various Mach numbers, Int. J. Numer. Meth. Fluids, 89 (2019) 430-463.
	\bibitem{chen-hllc} S-S. Chen, B. Lin, Y Li, C. Yan, HLLC+: LOW-MACH SHOCK-STABLE HLLC-TYPE RIEMANN SOLVER FOR ALL-SPEED FLOWS, SIAM J. Sci. Comput., 42 (2020) B921-B950.
	\bibitem{qu-hllem} F. Qu, J. Chen, D. Sun, J. Bai, C. Yun, A new all-speed flux scheme for the Euler equations, Computers and Mathematics with Applications, 77 (2019) 1216-1231.
	\bibitem{dumbb} M. Dumbser, D. S. Balsara, A new efficient formulation of the HLLEM Riemann solver for general conservation and non-conservation hyperbolic systems, J. Comput. Phys. 304 (2016) 275-319.	
	\bibitem{zha1} G-C. Zha, E. Bilgen, Numerical Solutions of Euler equations by using a new flux vector splitting scheme, Int. J. Numer. Meth. Fluids 17 (1993) 115-144.
	\bibitem{ausm1} M-S. Liou, C.J. Steffen Jr., A New Flux Splitting Scheme, J. Comput. Phys. 107 (1993) 23-39.
	\bibitem{tramel} R. W. Tramel, R. H. Nichols, P. G. Buning, Addition of Improved Shock-Capturing Schemes to OVERFLOW 2.1, AIAA-2009-3988.
	\bibitem{van-leer-muscl} B. van Leer, Towards the ultimate conservative difference scheme, A second-order sequel to Godunov's Method, J. Comput. Phys. 32 (1979) 101-136.
	\bibitem{gott} S. Gottlieb and C-W. Shu, Total Variation Diminishing Runge-Kutta Schemes, Mathematics of Computation, 67 (1998) 73-85.
	\bibitem{tan} S. Tan, L. Hu and H. Yuan, Development of a Shock Stable and Contact-Preserving Scheme for Multidimensional Euler Equations, AIAA Journal 2022 (60) 5232-5248.
	\bibitem{Hu} L. Hu, L. Yuan and  K. Zhao, Development of accurate and robust genuinely two-dimensional HLL-type Riemann solver for compressible flows, Comput. Fluids 213 (2020) 104719.
	\bibitem{chen-aupm} S-S. Chen, C. Yan, K. Zhong, H-C. Xue and E-L. Li, A novel flux splitting scheme with robustness and low dissipation for hypersonic heating prediction, International Journal of Heat and Mass Transfer 127 (2018) 126-137. 
\end{thebibliography}
\end{document}